\documentclass[a4paper]{amsart}

\usepackage[T1]{fontenc}	
\usepackage{graphicx}
\usepackage{amsmath,amsthm, amssymb}
\usepackage{upgreek}
\usepackage{mathrsfs}
\usepackage{pdfsync}
\usepackage[a4paper,top=3.5cm,bottom=3cm,left=3.5cm,right=3.5cm, bindingoffset=0mm]{geometry}

\usepackage[usenames,dvipsnames]{xcolor}
\usepackage[inline,shortlabels]{enumitem}
\usepackage[all]{xy}									
\usepackage{dsfont}									

\usepackage{hyperref} 
\usepackage{mathtools}

\allowdisplaybreaks

\def\@fnsymbol#1{\ensuremath{\ifcase#1\or *
\or
   \mathsection\or \mathparagraph\or \|\or **\or \dagger\dagger
   \or \ddagger\ddagger \else\@ctrerr\fi}}

\newcommand{\boldBeta}{\boldsymbol{\mathsf{B}}}
\newcommand{\boldGamma}{{\boldsymbol\Gamma}}
\newcommand{\boldDelta}{{\boldsymbol\Delta}}
\newcommand{\boldLambda}{{\boldsymbol\Lambda}}

\newcommand{\boldTau}{\mbfT}
\newcommand{\boldPhi}{{\boldsymbol\Phi}}

\newcommand{\boldUpsilon}{{\boldsymbol\Upsilon}}

\newcommand{\Beta}{\mathsf{B}}
\newcommand{\Eta}{\boldsymbol\eta}
\newcommand{\Mu}{\mu}

\newcommand{\boldgamma}{{\boldsymbol \gamma}}

\usepackage{xparse}

\ExplSyntaxOn
\NewDocumentCommand{\makeabbrev}{mmm}
 {
  \yoruk_makeabbrev:nnn { #1 } { #2 } { #3 }
 }

\cs_new_protected:Npn \yoruk_makeabbrev:nnn #1 #2 #3
 {
  \clist_map_inline:nn { #3 }
   {
    \cs_new_protected:cpn { #2 } { #1 { ##1 } }
   }
 }
 \ExplSyntaxOff

\makeabbrev{\textbf}{tbf#1}{a,b,c,d,e,f,g,h,i,j,k,l,m,n,o,p,q,r,s,t,u,v,w,x,y,z,A,B,C,D,E,F,G,H,I,J,K,L,M,N,O,P,Q,R,S,T,U,V,W,X,Y,Z}

\makeabbrev{\textbf}{bf#1}{a,b,c,d,e,f,g,h,i,j,k,l,m,n,o,p,q,r,s,t,u,v,w,x,y,z,A,B,C,D,E,F,G,H,I,J,K,L,M,N,O,P,Q,R,S,T,U,V,W,X,Y,Z}

\makeabbrev{\textsf}{tsf#1}{a,b,c,d,e,f,g,h,i,j,k,l,m,n,o,p,q,r,s,t,u,v,w,x,y,z,A,B,C,D,E,F,G,H,I,J,K,L,M,N,O,P,Q,R,S,T,U,V,W,X,Y,Z}

\makeabbrev{\mathsf}{mss#1}{a,b,c,d,e,f,g,h,i,j,k,l,m,n,o,p,q,r,s,t,u,v,w,x,y,z,A,B,C,D,E,F,G,H,I,J,K,L,M,N,O,P,Q,R,S,T,U,V,W,X,Y,Z}

\makeabbrev{\mathfrak}{mf#1}{a,b,c,d,e,f,g,h,i,j,k,l,m,n,o,p,q,r,s,t,u,v,w,x,y,z,A,B,C,D,E,F,G,H,I,J,K,L,M,N,O,P,Q,R,S,T,U,V,W,X,Y,Z}

\makeabbrev{\mathrm}{mrm#1}{a,b,c,d,e,f,g,h,i,j,k,l,m,n,o,p,q,r,s,t,u,v,w,x,y,z,A,B,C,D,E,F,G,H,I,J,K,L,M,N,O,P,Q,R,S,T,U,V,W,X,Y,Z}

\makeabbrev{\mathbf}{mbf#1}{a,b,c,d,e,f,g,h,i,j,k,l,m,n,o,p,q,r,s,t,u,v,w,x,y,z,A,B,C,D,E,F,G,H,I,J,K,L,M,N,O,P,Q,R,S,T,U,V,W,X,Y,Z}

\makeabbrev{\mathcal}{mc#1}{A,B,C,D,E,F,G,H,I,J,K,L,M,N,O,P,Q,R,S,T,U,V,W,X,Y,Z}

\makeabbrev{\mathbb}{mbb#1}{A,B,C,D,E,F,G,H,I,J,K,L,M,N,O,P,Q,R,S,T,U,V,W,X,Y,Z}

\makeabbrev{\mathscr}{ms#1}{A,B,C,D,E,F,G,H,I,J,K,L,M,N,O,P,Q,R,S,T,U,V,W,X,Y,Z}

\makeabbrev{\mathrm}{#1}{
Id,id,ran,rk,diag,stab,ann,conv,pr,ev,tr,End,Hom,sgn,im,op,can,fin,ext,red,tot,
rot,usc,lsc,Lip,lip,bLip,osc,AC,loc,spec,
supp,Opt,Adm,Cpl,Geo,GeoOpt,GeoAdm,GeoCpl,reg,
bd,co,Ric,Exp,dExp,dist,seg,Seg,cut,fcut,Cut,SDiff,Iso,Isom,diam,cl,Homeo,Diff,Der,vol,dvol,inj,relint,
var,law,Var,Poi,Gam,pa,so,iso,fs,inv,pqi,mix,
}

\makeabbrev{\mathsf}{#1}{CD,BE,RCD,MCP,Ent,wMTW,MTW,Ch}

\newcommand{\KonvRe}{{\mathsc{af}}}
\newcommand{\vReStu}{{\mathsc{wd}}}
\newcommand{\CoalFrag}{{\mathsc{cf}}}

\newcommand{\volm}{\boldsymbol{\mathsf{m}}}
\newcommand{\hvolm}{\widehat{\boldsymbol{\mathsf{m}}}}
\newcommand{\Heat}{\boldsymbol{\mathsf{H}}}
\newcommand{\dirE}{\boldsymbol{\mathsf{E}}}
\newcommand{\heat}{\boldsymbol{\mathsf{h}}}
\newcommand{\BM}{\boldsymbol{\mathsf{W}}}
\newcommand{\TT}{\boldsymbol{\mathsf{T}}}
\newcommand{\Gen}{\boldsymbol{\mathsf{L}}}

\newcommand{\eps}{\varepsilon}

\renewcommand{\div}{\mathrm{div}}
\newcommand{\defeq}{\eqqcolon}
\renewcommand{\complement}{\mathrm{c}}

\newcommand{\acts}{\,\raisebox{\depth}{\scalebox{1}[-1]{$\circlearrowleft$}}\,}

\newcommand{\mathsc}[1]{\text{\textsc{#1}}}

\newcommand{\emparg}{{\,\cdot\,}}

\newcommand{\Ppaiso}{\andi{\msP^\pa}}

\newcommand{\Ppa}{\msP^\pa}
\newcommand{\Ppafs}{\msP^{\pa,\fs}}
\newcommand{\Pfs}{\msP^\fs}
\newcommand{\Ppaisofs}{\andi{\msP^{\pa,\fs}}}

\newcommand{\slo}[1]{\abs{D#1}}

\DeclareMathOperator{\grad}{\boldsymbol\nabla\!}

\newcommand{\Proj}[2]{{{\rm{pr}}^{#1}_{#2}}}

\newcommand{\forallae}[1]{{\textrm{\,for ${#1}$-a.e.\,}}}

\newcommand{\Bo}[1]{\mcB_{#1}}
\newcommand{\T}{\tau}

\newcommand{\Leb}{{\mathrm{Leb}}}

\renewcommand{\Cap}{\mathrm{cap}}

\newcommand{\dom}[1]{\msD(#1)}

\DeclareMathOperator{\eqdef}{\coloneqq}

\let\epsilon\varepsilon

\newcommand{\longrar}{\longrightarrow}
\newcommand{\rar}{\rightarrow}

\newcommand{\nlim}{\lim_{n}}								
\newcommand{\mlim}{\lim_{m}}

\newcommand{\nliminf}{\liminf_{n }}

\newcommand{\nlimsup}{\limsup_{n }\,}

\newcommand{\diff}{\mathop{}\!\mathrm{d}}						
	
\DeclareMathOperator{\orb}{orb}

\newcommand{\tabs}[1]{\big\lvert#1\big\rvert}	
\newcommand{\abs}[1]{\left\lvert#1\right\rvert}						
\newcommand{\ttabs}[1]{\lvert#1\rvert}						
\newcommand{\tnorm}[1]{\big\lVert#1\big\rVert}					
\newcommand{\norm}[1]{\left\lVert#1\right\rVert}					
\newcommand{\set}[1]{\left\{#1\right\}}							
\newcommand{\tset}[1]{\big\{#1\big\}}							
\newcommand{\floor}[1]{\left\lfloor#1\right\rfloor}					
\newcommand{\tonde}[1]{\left(#1\right)}							
\newcommand{\ttonde}[1]{\big({#1}\big)}
\newcommand{\quadre}[1]{\left[#1\right]}							
\newcommand{\tquadre}[1]{\big[#1\big]}							
\newcommand{\tscalar}[2]{\big\langle #1 \, \big |\, #2\big\rangle}		
\newcommand{\scalar}[2]{\left\langle #1 \,\middle |\, #2\right\rangle}		
\newcommand{\qvar}[1]{\quadre{#1}}								
	
\newcommand{\seq}[1]{\tonde{#1}}								
\newcommand{\tseq}[1]{{\big(#1\big)}}
\newcommand{\ttseq}[1]{{(#1)}}
\newcommand{\Cb}{\mcC_b}									
\newcommand{\Cc}{\mcC_c}									
\newcommand{\Cz}{\mcC_0}									

\newcommand{\Mp}{\mathscr M^+}
\newcommand{\Mb}{\mathscr M_b}
\newcommand{\Mbp}{\mathscr M_b^+}

\newcommand{\pfwd}{\sharp}

\DeclareMathOperator*{\essinf}{essinf}

\newcommand{\quotient}[2]{\left.\raisebox{.0em}{$#1\!$}\middle/\raisebox{-.0em}{$#2$}\right.}

\DeclareMathOperator{\car}{\mathds 1}

\DeclareMathOperator{\emp}{\varnothing} 
\DeclareMathOperator{\N}{{\mathbb N}}
\newcommand{\R}{{\mathbb R}}
\DeclareMathOperator{\Q}{{\mathbb Q}}
\DeclareMathOperator{\Z}{{\mathbb Z}}

\DeclareMathOperator{\cont}{cont} 

\newcommand{\restr}{\big\lvert}	
\newcommand{\trid}{{\star}}
\newcommand{\Fl}{\Psi}
\newcommand{\fl}{\uppsi}

\newcommand{\iref}[1]{\ref{#1}}

\newcommand{\comma}{\,\,\mathrm{,}\;\,}
\newcommand{\semicolon}{\,\,\mathrm{;}\;\,}
\newcommand{\fstop}{\,\,\mathrm{.}}

\DeclareMathOperator{\zero}{{\mathbf 0}}
\DeclareMathOperator{\uno}{{\mathbf 1}}

\DeclareMathOperator{\bexp}{\mathbf{exp}}

\DeclareMathOperator{\inter}{int}

\newcommand{\n}[1]{{\overline{#1}}}

\newcommand{\gscal}[2]{\scalar{#1}{#2}_\mssg}
\newcommand{\tgscal}[2]{\tscalar{#1}{#2}_\mssg}
\newcommand{\hTF}[2]{\widehat\mfZ^{#1}_{#2}}
\newcommand{\hTFW}[2]{\widehat\mfW^{#1}_{#2}}
\newcommand{\TF}[1]{\mfZ^{#1}}
\newcommand{\TFL}[2]{\mfL^{#1,#2}}
\newcommand{\hTFB}[1]{\mfB^{#1}}
\newcommand{\Cyl}[1]{\mcF\mcC^{#1}}

\newcommand{\deq}{\overset{\mrmd}{=}}

\let\temp\phi
\let\phi\varphi
\let\varphi\temp

\newcommand{\andi}[1]{#1_\circ}
\newcommand{\Tso}{\andi{\boldTau}}
\newcommand{\To}{\boldTau}
\newcommand{\hM}{\widehat\mbfM}
\newcommand{\andihM}{\andi{\widehat\mbfM}}

\newcommand{\DF}{{\mcD}}

\newcommand{\GP}{{\mcG}}

\newcommand{\PD}{{\Pi}}

\newcommand{\Vect}{\mfX}
\newcommand{\man}{M}
\newcommand{\hman}{\widehat M}

\numberwithin{equation}{section}
\theoremstyle{plain}
\newtheorem{thm}{Theorem}[section]
\newtheorem*{thm*}{Theorem}
\newtheorem*{mthm*}{Main Theorem}

\newtheorem{prop}[thm]{Proposition}
\newtheorem{lem}[thm]{Lemma}
\newtheorem{cor}[thm]{Corollary}

\theoremstyle{definition}
\newtheorem{defs}[thm]{Definition}
\newtheorem*{defs*}{Definition}

\theoremstyle{remark}
\newtheorem{rem}[thm]{Remark}
\newtheorem{ass}{Assumption}

\renewcommand{\paragraph}[1]{\medskip \emph{#1}. \ }

\begin{document}

\title[Dirichlet--Ferguson Diffusions]{The Dirichlet--Ferguson diffusion\\on the space of probability measures\\over a closed Riemannian manifold\textsuperscript{*}
}
\thanks{\textsuperscript{*}Research supported by the Sonderforschungsbereich~1060 and the Hausdorff Center for Mathematics.}

\author[L.\ Dello Schiavo]{Lorenzo Dello Schiavo\textsuperscript{\S}}


\address{
\noindent
Institut f\"ur Angewandte Mathematik
\newline
Rheinische Friedrich-Wilhelms-Universit\"at Bonn
\newline
Endenicher Allee 60
\newline
DE 53115 Bonn
\newline
Germany
}

\thanks{\textsuperscript{\S}I am especially grateful to Prof.s K.-T.~Sturm, E.~W.~Lytvynov, M.~Gordina and M.~Gubinelli for several remarks and comments. I thank Prof.~G.~Last for a preliminary version of~\cite{Las18}.
Finally, I am also grateful to Prof.s~F.~Bassetti and E.~Dolera and to Dr.~C.~Orrieri for discussions about the Dirichlet--Ferguson measure during my stay, in February 2018, at the Department of Mathematics `F.~Casorati' of the Universit\`a degli Studi di Pavia. I~thank the Department for their kind hospitality.
}

\email{delloschiavo@iam.uni-bonn.de}
\date{\today}



\begin{abstract}
We construct a recurrent diffusion process with values in the space of probability measures over an arbitrary closed Riemannian manifold of dimension~$d\geq 2$. The process is associated with the Dirichlet form defined by integration of the Wasserstein gradient w.r.t.~the Dirichlet--Ferguson measure, and is the counterpart on multi-dimensional base spaces to the Modified Massive Arratia Flow over the unit interval described in V.~Konarovskyi, M.-K.\ von~Renesse, \emph{Comm.\ Pure Appl.\ Math.}, 72, 0764--0800 (2019). Together with two different constructions of the process, we discuss its ergodicity, invariant sets, finite-dimensional approximations, and Varadhan short-time asymptotics.
\end{abstract}

\maketitle

\tableofcontents

\section{Introduction}
We provide two constructions of a Markov diffusion~$\Eta_\bullet$ with values in the space of probability measures~$\msP$ over a closed Riemannian manifold~$\man$ of dimension~$d\geq 2$.

On the one hand, combining results of Bendikov--Saloff-Coste~\cite{BenSaC97} and Albeverio--Daletskii--Kondratiev~\cite{AlbDalKon97,AlbDalKon00} about elliptic diffusions on infinite products, we characterize~$\Eta_\bullet$ as the superprocess constituted by up-to-countable independent massive Brownian particles with volatility equal to their inverse mass. Thus, we may regard~$\Eta_\bullet$ as the counterpart over~$\man$ of Konarovskyi's Modified Massive Arratia Flow~\cite{Kon17} over the unit interval. Here, no coalescence occurs by reasons of the dimension of~$\man$.

On the other hand, we show that~$\Eta_\bullet$ is associated with a symmetric regular strongly local Dirichlet form~$\mcE$ on the space of real-valued functions on~$\msP$ square-integrable with respect to the Dirichlet--Ferguson random measure~$\DF$~\cite{Fer73}. The form~$\mcE$ is defined as the closure of the Dirichlet integral induced by~$\DF$ and by the natural gradient of the $L^2$-Wasserstein geometry of~$\msP$, on the algebra of cylinder functions induced by smooth potential energies in the sense of Otto calculus~\cite{Ott01}. Thus, we may regard~$\Eta_\bullet$ as a possible candidate for a ``Brownian motion'' ---~that is, a canonical diffusion process~--- on the $L^2$-Wasserstein space~$\msP_2$.

\smallskip

Let~$(\man,\mssg)$ be a closed Riemannian manifold of dimension~$d\geq 2$, with volume measure~$\mssm$, Riemannian distance~$\mssd$, and Laplace--Beltrami operator~$\Delta^\mssg$. Let~$(\msP_2,W_2)$ be the $L^2$-Wasserstein metric space over~$(\man,\mssd)$, endowed with Otto's metric~$\scalar{\emparg}{\emparg}_{T_\mu\msP_2}$ and with the Dirichlet--Ferguson (probability) measure~$\DF_\mssm$ with intensity measure~$\mssm$. For a suitable algebra~$\hTF{\infty}{0}$ of functions~$u\colon \msP_2\rar \R$, see~Definition~\ref{d:CylFunc} below, we prove the following.

\begin{mthm*}
The symmetric bilinear form~$(\mcE,\hTF{\infty}{0})$ given by
\begin{align*}
\mcE(u,v)\eqdef& \tfrac{1}{2}\int_{\msP_2} \scalar{\grad u(\mu)}{\grad v(\mu)}_{T_\mu\msP_2} \diff\DF_\mssm(\mu)\comma \qquad u,v\in\hTF{\infty}{0}\comma
\end{align*}
is closable. Its closure~$(\mcE,\dom{\mcE})$ is a regular strongly local recurrent (conservative) Dirichlet form with generator the (Friedrichs) extension of the essentially self-adjoint operator~$(\mbfL,\hTF{\infty}{0})$ given by
\begin{align*}
(\mbfL u)(\mu)\eqdef& \frac{1}{2} \int_\man \frac{\Delta^z\restr_{z=x} u(\mu+\mu\!\set{x}\delta_z-\mu\!\set{x}\delta_x)}{(\mu\!\set{x})^2}\diff\mu(x) \quad \forallae{\DF_\mssm}\mu\comma\; u\in\hTF{\infty}{0} \fstop
\end{align*}

Additionally:
\begin{itemize}
\item the energy measure space~$(\mcE,\msP_2,\DF_\mssm)$ has the Rademacher property, that is, for any $W_2$-Lipschitz function~$u\colon\msP_2\rar \R$,
\begin{align*}
u\in\dom{\mcE}\quad \text{and} \quad \scalar{\grad u(\mu)}{\grad u(\mu)}_{T_\mu\msP_2}\leq \Lip[u]^2 \quad\forallae{\DF_\mssm} \mu\semicolon
\end{align*}

\item the Markov kernel
\begin{align*}
p_t(A_1,A_2)\eqdef \int_{A_1} e^{-t\, \mbfL}\car_{A_2}  \diff\DF_\mssm \comma\qquad A_1,A_2\subset \msP_2\comma
\end{align*}
associated to~$(\mcE,\dom{\mcE})$ satisfies the one-sided Varadhan-type upper estimate
\begin{align*}
\lim_{t\downarrow 0} t\, \log p_t(A_1, A_2)\leq -\tfrac{1}{2} \inf_{\mu_i\in A_i} W_2(\mu_1,\mu_2)^2 \semicolon
\end{align*}

\item $(\mcE,\dom{\mcE})$ is properly associated with a $\msP$-valued Markov diffusion with infinite life-time
\begin{align*}
\Eta_\bullet\eqdef \ttonde{\Omega, \mcF,\seq{\mcF_t}_{t\geq 0}, \seq{\eta_t}_{t\geq 0}, \seq{P_\eta}_{\eta\in\msP}} \semicolon
\end{align*}

\item for every initial distribution~$\law(\Eta_0)\ll\DF_\mssm$, for every~$u\in\hTF{2}{0}$, the process
\begin{align*}
M^u_t\eqdef u(\Eta_t)-u(\Eta_0)-\int_0^t (\mbfL u)(\Eta_s) \diff s
\end{align*}
is an $\seq{\mcF_t}_{t\geq 0}$-martingale with quadratic variation process
\begin{align*}
\qvar{M^u}_t=\int_0^t \scalar{\grad u(\Eta_s)}{\grad u(\Eta_s)}_{T_{\Eta_s}\msP_2} \diff s \fstop
\end{align*}
\end{itemize}
\end{mthm*}


\section{Motivations and main results}\label{s:Intro}

In the last two decades, the space~$\msP$ of all Borel probability measures over a closed Riemannian manifold~$(\man,\mssg)$, endowed with the $L^2$-Kantorovich--Rubin\-shtein distance~$W_2$~\eqref{eq:WD}, has proven both a powerful tool and an interesting geometric object in its own right.
Since the fundamental works of Y.~Brenier,~R.~J.~McCann, F.~Otto, C.~Villani and many others, e.g.,~\cite{Bre91, Ott01, McC97,Vil09}, several geometric notions have been introduced, including those of geodesic curve, tangent space~$T_\mu\msP$ at a point~$\mu\in\msP$ and gradient~$\grad u(\mu)$ of a scalar-valued function~$u$ at~$\mu$, e.g.,~\cite{Gig11, Gig12, Lot07}.
Indeed, the metric space~$\msP_2\eqdef (\msP,W_2)$ may ---~to some extent~--- be regarded as a kind of infinite-dimensional Riemannian manifold. Furthermore, provided that~$(\man,\mssg)$ be a closed manifold with non-negative sectional curvature,~$\msP_2$ has non-negative lower curvature bound in the sense of Alexandrov~\cite[Thm.~2.20]{AmbGig11}.

\subsection{Volume measures on~$\msP_2$} The question of the existence of a \emph{Riemannian} volume measure on~$\msP_2$, say~$\dvol_{\msP_2}$, has been insistently posed and remains to date not fully answered.
A first natural requirement for such a measure ---~if any~--- is an \emph{integration-by-parts} formula for the gradient, which would imply the closability of the form
\begin{align}\label{eq:Intro:DirForm}
\mcE(u,v)\eqdef \tfrac{1}{2}\int_{\msP} \scalar{\grad u(\mu)}{\grad v(\mu)}_{T_\mu\msP_2} \dvol_{\msP_2}(\mu)\fstop
\end{align}
In turn, the theory of Dirichlet forms would then grant the existence of a diffusion process associated to~$\mcE$ and thus deserving the name of \emph{Brownian motion} on~$\msP_2$.

Further requirements are the validity of a Rademacher-type Theorem, i.e. the $\dvol_{\msP_2}$-a.e. differentiability of $W_2$-Lipschitz functions, which motivated the work~\cite{LzDS19b}, and of its converse, the Sobolev-to-Lipschitz property. Together, these properties would grant the identification of~$W_2$ with the intrinsic distance induced by~$\mcE$.

\subsection{Diffusions processes on~$\msP$}
In the case when~$\man=\mbbS^1$, the unit sphere, or~$\man=I$, the closed unit interval, M.-K.~von~Renesse and K.-T.~Sturm proposed the \emph{entropic measure}~$\mbbP^\beta$~\cite[Dfn.~3.3]{vReStu09} as a candidate for~$\dvol_{\msP_2}$, and constructed the associated \emph{Wasserstein diffusion}~$\Mu^\vReStu_\bullet$. Whereas the construction of the entropic measure in the case when~$\man$ is an arbitrary closed Riemannian manifold was subsequently achieved by K.-T.~Sturm in~\cite{Stu11}, many of its properties, and in particular the closability of the associated form~\eqref{eq:Intro:DirForm}, remain unknown.

Similar constructions to the Wasserstein diffusion ---~up to now confined to one-dimension\-al base spaces~--- include J.~Shao's \emph{Dirichlet--Wasserstein diffusion}~\cite{Sha11}, when~$\man=\mbbS^1$ or~$I$; V.~Konarovskyi \emph{modified massive Arratia flow}~$\Mu^\KonvRe_\bullet$~\cite{Kon17,KonvRe18}, when~$\man=I$; and V.~Konarovskyi and M.-K. von Renesse's \emph{coalescing-fragment\-ating Wasserstein dynamics}~$\Mu^{\CoalFrag}_\bullet$ \cite{KonvRe17}, when~$\man=\R$.
Finally, it is worth mentioning two constructions in the case~$\man=\R^d$, namely the superprocesses of stochastic flows introduced by Z.-M.~Ma and K.-N.~Xiang in~\cite{MaXia01}, and the recent work~\cite{ChoGan17} by Y.~T.~Chow and W.~Gangbo, concerned with a stochastic process
on~$\msP_2$ \guillemotleft modeled after Brownian motion\guillemotright~and generated by a \guillemotleft partial Laplacian\guillemotright.

\subsection{A canonical process}
 If not otherwise stated, we shall assume the following.

\begin{ass}[Riemannian manifolds]\label{ass:RM}
By a \emph{Riemannian manifold} we shall mean any closed (i.e. compact, without boundary) connected oriented smooth Riemannian manifold~$(\man,\mssg)$ with (smooth) Riemannian metric~$\mssg$, intrinsic distance~$\mssd_\mssg$, Borel $\sigma$-algebra~$\Bo{\mssg}$, (Borel) volume measure~$\mssm$, \emph{normalized} volume measure~$\n\mssm$, and heat kernel~$\mssh_t(x,\diff y)$. If not otherwise stated, we shall assume that~$d\eqdef \dim \man\geq 2$.
\end{ass}

We construct a stochastic diffusion process
\begin{align}\label{eq:Process}
\Eta_\bullet\eqdef\tseq{\Omega, \mcF,\seq{\mcF_t}_{t\geq 0}, \seq{\Eta_t}_{t\geq 0}, \seq{P_\eta}_{\eta\in\msP_2}}
\end{align}
with state space~$\msP_2$, modeled after the Brownian motion on~$\man$. By this we mean that~$\Eta_\bullet$ enjoys the following property: Let~$\seq{\Eta_t^{\eta_0}}_t$ denote a stochastic path of~$\Eta_\bullet$ starting at~$\eta_0$. If
\begin{align*}
\eta_0\eqdef(1-r)\delta_{x_0}+r\delta_{y_0}\comma\qquad x_0, y_0\in \man\comma r\in I\comma
\end{align*}
then
\begin{align}\label{eq:ass:A}
\Eta_t^{\eta_0}(\omega)=(1-r)\delta_{x_{t/(1-r)}(\omega)}+r\delta_{y_{t/r}(\omega)} \comma
\end{align}
where~$x_\bullet$ and~$y_\bullet$ are \emph{independent} Brownian motions on~$(\man,\mssg)$ respectively starting at~$x_0$,~$y_0$.
For~$r\in\set{0,1}$, Equation~\eqref{eq:ass:A} entails that~$\Eta_\bullet$ respects the Dirac embedding, that is
\begin{align*}
\Eta_0=\delta_{x_0}\quad\implies\quad \Eta_t(\omega)=\delta_{x_t(\omega)}
\end{align*}
for some Brownian motion~$x_\bullet$ starting at~$x_0$. This is a natural requirement, since
\begin{align*}
\delta\colon (\man,\mssd_\mssg)\rar (\msP_2,W_2)\comma \quad x\mapsto \delta_x
\end{align*}
is an isometric embedding.
If~$r\in (0,1)$, then~\eqref{eq:ass:A} and its straightforward $n$-points generalization may be easily interpreted in terms of particle systems. Indeed,~$\Eta_\bullet$ as in~\eqref{eq:ass:A} describes the evolution of the two massive particles~$(x_0,1-r)$ and~$(y_0,r)$, and translates into the requirement that the evolution of their positions be independent up to the choice of suitable volatilities, namely the inverse of the mass carried by each atom.

We will provide two different constructions of~$\Eta_\bullet$.

\subsection{Construction via semigroups}
In the following endow~$\mbfI\eqdef [0,1]^{\infty}$ with the product topology and set
\begin{equation}\label{eq:DeltaTau}
\begin{aligned}
\boldDelta\eqdef&\set{\mbfs\eqdef \seq{s_i}_i\in \mbfI : \sum_i^\infty s_i=1}\comma
\\
\To\eqdef&\set{\mbfs\eqdef \seq{s_i}_i \in \boldDelta : s_i \geq s_{i+1} \geq 0 \text{~~for all~}i} \fstop
\end{aligned}
\end{equation}

Let~$\Tso\subset \To$ be defined similarly to~$\To$ with~$>$ in place of~$\geq$. For~$\mbfs\in\Tso$  put~$\man_i\eqdef (\man,s_i\mssg)$ and consider the infinite product~$\mbfM=\prod_i^\infty \man_i$, endowed with the product measure
\begin{align*}
\n\volm\eqdef \bigotimes^\infty \n\mssm \fstop
\end{align*}

Letting~$\mbfx=\seq{x_i}_i^\infty$,~$\mbfy\eqdef\seq{y_i}_i^\infty$ and defining the family of heat kernel measures
\begin{align}\label{eq:Intro:heat}
\heat^\mbfs_t(\mbfx,\diff \mbfy)\eqdef \bigotimes_i^\infty \mssh_{t/s_i}(x_i,\diff y_i)\comma \qquad \mbfx\in \mbfM\comma\qquad t>0\comma
\end{align}
the resulting product semigroup~$\seq{\Heat^\mbfs_t}_{t\geq 0}$ given by
\begin{align}\label{eq:Intro:ProdSemigroup}
(\Heat^\mbfs_t u)(\mbfx)\eqdef \int_\mbfM u(\mbfy) \, \heat^\mbfs_t(\mbfx,\diff \mbfy)\comma \qquad u\in L^2(\mbfM, \n\volm)\comma \qquad t>0\comma
\end{align}
is an ergodic Markov semigroup with invariant measure~$\n\volm$. By the general results of A.~Bendi\-kov and L.~Saloff-Coste~\cite{BenSaC97} about infinite-dimensional elliptic diffusions, $\Heat^\mbfs_t$~admits a density, i.e.\ 
\begin{align*}
\heat^\mbfs_t(\mbfx,\diff\mbfy)=\heat^\mbfs_t(\mbfx,\mbfy)\diff\n\volm(\mbfy)\comma \qquad \mbfx,\mbfy\in \mbfM\comma t>0\comma
\end{align*}
additionally continuous and bounded on~$(0,\infty)\times\mbfM^{\times 2}$ for every~$\mbfs\in\Tso$.
We denote by
\begin{align}\label{eq:Intro:BM}
\BM^\mbfs_\bullet\eqdef \tseq{\Omega,\mcF,\seq{\mcF_t}_{t\geq 0},\seq{\BM^\mbfs_t}_{t\geq 0}, \seq{P^\mbfs_\mbfx}_{\mbfx\in \mbfM}}
\end{align}
the associated time-homogeneous recurrent ergodic Markov process with state space~$\mbfM$ and transition kernels~$\seq{\heat^\mbfs_\bullet(\mbfx,\emparg)}_{\mbfx\in\mbfM}$.

\smallskip

Let now~$\mbfP$ be any probability on~$\Tso$, so that~$(\Tso,\mbfP)$ is a standard Borel probability space.
The semigroup defined on~$\hM\eqdef \To\times\mbfM$ as
\begin{equation}\label{eq:IntroHatHeat}
\begin{aligned}
(\widehat\Heat_t v)(\mbfs,\mbfx)\eqdef& \ttonde{(\id \otimes \Heat^\mbfs_t)\, v}(\mbfs,\mbfx)\comma \qquad && (\mbfs,\mbfx)\in \hM\comma t>0\comma
\\
=&\ttonde{\Heat^\mbfs_t\, v(\mbfs, \emparg)}(\mbfx)\comma && v\in L^2(\hM, \mbfP\otimes\n\volm)
\end{aligned}
\end{equation}
is itself a Markov semigroup on~$L^2(\hM, \mbfP\otimes\n\volm)$.
Setting
\begin{align*}
\andi{\mbfM}\eqdef\set{\mbfx\in\mbfM : x_i\neq x_j \text{ for } i\neq j}\comma
\end{align*}
the map
\begin{align}\label{eq:Intro:Phi}
\boldPhi\colon \boldDelta\times\mbfM \longrar \msP_2\comma\qquad \boldPhi(\mbfs,\mbfx)\eqdef \sum_i^\infty s_i\delta_{x_i}
\end{align}
is injective when {restricted} to~$\andihM\eqdef \Tso\times\andi{\mbfM}$.

Say that~$\andi{\mbfM}$ is \emph{$\BM^\mbfs_\bullet$-coexceptional} if~$\BM^\mbfs_t\in\andi{\mbfM}$ for every~$t>0$. We will show in Lemma~\ref{l:Except} that~$\andi{\mbfM}$ is $\BM^\mbfs_\bullet$-coexceptional for every~$\mbfs\in\Tso$.
Therefore,~$\andihM$ is coexceptional for the process~$\widehat\BM_\bullet$ associated to~$\widehat\Heat_t$. Indicate by~$\BM^{\mbfs;\mbfx_0}_t$ the process~$\BM^\mbfs_t$ started at~$\mbfx_0$.
Provided that~$\boldPhi$ is suitably measurable, we may consider the induced stochastic process on~$\msP$, pathwise defined by
\begin{align}\label{eq:Intro:Eta}
\Eta_t^{\eta_0}\eqdef \boldPhi\circ\widehat\BM^{\mbfs,\mbfx_0}_t=\boldPhi\ttonde{\mbfs,\BM^{\mbfs;\mbfx_0}_t} \comma \qquad \eta_0\eqdef \boldPhi\ttonde{\mbfs,\mbfx_0}\comma \qquad t> 0\fstop
\end{align}

By construction,~$\Eta_\bullet$ is a time-homogeneous Markov process with state space~$\msP_2$. However, since~$\boldPhi$ is not continuous, it is not clear at this stage whether~$\Eta_\bullet$ has continuous paths, and its properties may wildly vary, depending on the choice of the law~$\mbfP$ for the starting point~$\mbfs$.

\subsection{A choice for~$\dvol_{\msP_2}$}
Everywhere in the following we let~$\beta>0$ be fixed. For the moment, we shall think of~$\beta$ as the total volume of~$\man$, so that~$\mssm= \beta\n\mssm$. Let
\begin{align*}
\diff\Beta_\beta(r)\eqdef \beta(1-r)^{\beta-1}\diff r
\end{align*}
be the \emph{Beta distribution} on~$I$ with shape parameters~$1$ and~$\beta$. We denote by
\begin{align*}
\boldBeta_\beta\eqdef \bigotimes^\infty \Beta_\beta
\end{align*}
the corresponding product measure on~$\mbfI$.

The push-forward measure of a measure~$\mu$ by measurable map~$T$ is
\begin{align*}
T_\pfwd \mu\eqdef \mu\circ T^{-1}\fstop
\end{align*}

We set~$\Ppaiso=\boldPhi(\andihM)$, the space of \emph{purely atomic} probability measures with \emph{infinite strictly ordered} masses, dense in the compact space~$\msP_2$,~\cite[Thm.~6.18, Rmk.~6.19]{Vil09}. Thus, if we assume~$\mbfP$ to be fully supported on~$\To$, then
\begin{align*}
\mbfQ\eqdef \boldPhi_\pfwd (\mbfP\otimes\n\volm)
\end{align*}
is fully supported on~$\msP_2$.
Under such assumption, the property of~$\mbfQ$ being a `canonical' measure ---~in any suitable sense~--- on~$\msP$ is equivalent to that of~$\mbfP$ being `canonical' on~$\To$. As a candidate for~$\mbfP$ we choose the \emph{Poisson--Dirichlet measure}~$\PD_\beta$ introduced by J.~F.~C.~Kingman in~\cite{Kin75}. We recall its definition following the neat exposition of P.~Donnelly and G.~Grimmet~\cite{DonGri93}.

\begin{defs}[Poisson--Dirichlet measure] For~$\mbfr\eqdef \seq{r_i}_i\in \mbfI$ we denote by~$\boldUpsilon(\mbfr)$ the vector of its entries in non-increasing order, by~$\boldUpsilon\colon \mbfI\rar \To\subset \mbfI$ the reordering map, measurable by~\cite[p.~91]{DonJoy89}. Further, we let~$\boldLambda\colon \mbfI\rar \boldDelta$ be defined by
\begin{equation}\label{eq:LambdaReArr}
\begin{aligned}
\Lambda_1(r_1)\eqdef& r_1\comma\qquad \Lambda_k(r_1,\dotsc, r_k)\eqdef r_k\prod_i^{k-1} (1-r_i)\comma 
\\
\boldLambda(\mbfr)\eqdef& \seq{\Lambda_1(r_1), \Lambda_2(r_1,r_2),\dotsc}\fstop
\end{aligned}
\end{equation}
The \emph{Poisson--Dirichlet measure~$\PD_\beta$ with parameter~$\beta$} on~$\To$, concentrated on~$\Tso$, is
\begin{align*}
\PD_\beta\eqdef (\boldUpsilon\circ\boldLambda)_\pfwd \boldBeta_\beta \fstop
\end{align*}
\end{defs}

For such a choice of~$\mbfP$, the measure~$\mbfQ$ is the \emph{Dirichlet--Ferguson measure}~$\DF_\mssm$ introduced by T.~S.~Ferguson in his seminal work~\cite{Fer73}. The dependence of~$\DF_\mssm$ on~$\beta$ is implicit in the constraint~$\mssm=\beta\n\mssm$. Since in the following~$\DF_\mssm$ will play the role of~$\dvol_{\msP_2}$, we state here some of its several characterizations.

\begin{thm}[A characterization of~$\DF_\mssm$] Let~$\mbfQ$ be a probability measure on~$\msP$. For~$\eta\in\msP$, $x\in \man$ and~$r\in I$, set
\begin{align}\label{eq:Intro:CoCo}
\eta_x\eqdef \eta\!\set{x}\in I \qquad \text{and} \qquad \eta^x_r\eqdef (1-r)\eta+r\delta_x\in\msP \fstop
\end{align}

Then, the following are equivalent:
\begin{itemize}
\item~$\mbfQ$ is the Dirichlet--Ferguson measure~$\DF_\mssm\eqdef \boldPhi_\pfwd (\PD_\beta\otimes\n\volm)$;
\item if~$\eta$ is a $\mbfQ$-distributed $\msP$-valued random field,~$x$ is $\n\mssm$-distributed and~$r$ is $\Beta_\beta$-distributed, then~$\mbfQ$ satisfies \emph{Sethuraman's fixed-point characterization}
\begin{align}\label{eq:Sethuraman}
\eta\deq \eta^x_r \comma
\end{align}
where~$\deq$ denotes equality in law, see~\emph{\cite[Eqn.~(3.2)]{Set94}};
\item~$\mbfQ$ satisfies the \emph{Mecke-type identity} or \emph{Georgii--Nguyen--Zessin formula},
\begin{align}\label{eq:Mecke}
\int_\msP \int_\man u(\eta, x, \eta_x) \diff\eta(x)\diff\mbfQ(\eta) =\int_\msP \int_{\man} \int_I u(\eta^x_r, x, r) \diff \Beta_\beta(r)\diff\n\mssm(x)\diff\mbfQ(\eta)
\end{align}
for any semi-bounded measurable~$u\colon\msP_2\times \man\times I\rar \R$, see~\emph{\cite{LzDSLyt17}}.
\end{itemize}
\end{thm}

We will mostly dwell upon the characterization~\eqref{eq:Mecke}, obtained with E.~W.~Lyt\-vynov in~\cite{LzDSLyt17}, originally proven in the form~\eqref{eq:Sethuraman} by J.~Sethuraman in~\cite{Set94}. 
See also~G.~Last~\cite{Las18} for a similar characterization on more general spaces,~\cite{LzDS19a} for a Fourier transform characterization and T.~J.~Jiang, J.~M.~Dickey, and K.-L.~Kuo's characterization via $c$-transform~\cite{JiaDicKuo04}.

\subsection{Construction via Dirichlet forms theory} By construction, the measure $\boldPhi^{-1}_\pfwd\DF_\mssm=\PD_\beta\otimes\n\volm$ is an invariant measure of~$\widehat\BM_\bullet$. Choosing~$\dvol_{\msP_2}=\DF_\mssm$ in~\eqref{eq:Intro:DirForm}, we will show that the process~$\Eta_\bullet$ in~\eqref{eq:Intro:Eta} is the Markov diffusion associated with the Dirichlet form~$\mcE$.
This requires however some preparations.

\medskip

We shall follow a similar strategy to the one adopted by Yu.~G.~Kondratiev, E.~W.~Lyt\-vynov and A.~M.~Vershik in~\cite{KonLytVer15}, where analogous results are presented for Gibbs measures on the space of non-negative Radon measures over~$\R^d$.
Firstly, let~$\hat f\colon \man\times I\rar \R$ be of the form~$\hat f\eqdef f\otimes \varrho$, where~$f\in\mcC^\infty(\man)$ and~$\varrho\in\mcC^\infty(I)$ is supported in the \emph{open} interval~$(0,1)$. Recalling the notation~$\eta_x\eqdef \eta\!\set{x}$, we further let
\begin{align}\label{eq:Intro:fTrid}
\hat f^\trid(\eta)\eqdef \int_\man f(x)\cdot \varrho(\eta_x) \diff\eta(x)
\end{align}
and consider
\begin{itemize}
\item the algebra~$\hTF{}{0}$ of cylinder functions~$u\colon \msP\rar \R$ of the form
\begin{align*}
u(\eta)=F\ttonde{\hat f_1^\trid(\eta),\dotsc, \hat f_k^\trid(\eta)}\comma
\end{align*}
where~$F\in \mcC^\infty_b(\R^k)$ and~$\hat f_i$ is as in~\eqref{eq:Intro:fTrid} for~$i\leq k$.

\item the algebra~$\hTFB{}$ of cylinder functions induced by measurable potential energies, viz.\
\begin{align*}
u(\eta)=F\ttonde{\eta f_1,\dotsc, \eta f_k}\comma
\end{align*}
where~$F\in \mcC^\infty_b(\R^k)$ and~$f_i\in\mcB_b(\man;\R)$ for~$i\leq k$;

\item the algebra~$\TF{}{}$ of cylinder functions induced by smooth potential energies, defined analogously to~$\hTFB{}$, with the additional requirement that~$f_i\in\mcC^\infty(\man)$ for~$i\leq k$.
\end{itemize}

Let now~$w$ be a smooth vector field and~$\tseq{\fl^{w,t}}_{t\geq 0}$ be the associated flow of diffeomorphisms~\eqref{eq:FlowDef}. For~$\mu\in\msP$ we denote by~$T^\Der_\mu\msP_2$ the completion of the space of all smooth vector fields~$w$ with respect to the pre-Hilbert norm~$w\mapsto\tnorm{\abs{w}_\mssg}_{L^2(\mu)}$. The superscript `$\Der$' stands for \emph{derivation}, see~\cite[\S6.1]{LzDS19b}.
It is well-established in the optimal transport theory that the \emph{tangent space}~$T_\mu\msP_2$ to the `Riemannian manifold'~$\msP_2$ at~$\mu$ is
\begin{align*}
T_\mu\msP_2\eqdef \cl_{T^\Der_\mu\msP_2} \set{\nabla f: f\in\mcC^\infty(\man)} \comma
\end{align*}
e.g.,~\cite[2.31 and \S7.2]{AmbGig11}, cf. also \cite{Gig11, Gig12, GanKimPac10}.
The inclusion~$T_\mu\msP_2\subset T^\Der_\mu\msP_2$ is generally a strict one. We shall make use of both definitions, the interplay of which was detailed in~\cite{LzDS19b}. 
The `directional derivative' of functions~$u\in \hTF{}{0}$ or~$\TF{}$ in the smooth `direction'~$w$ is given by
\begin{align}\label{eq:Intro:DirDer}
\grad_w u(\mu)\eqdef& \diff_t\restr_{t=0} u\ttonde{\fl^{w,t}_\pfwd \mu} \comma
\end{align}
Lemma~\ref{l:DirDer}. For~$u,v\in\hTF{}{0}$ and any smooth vector field~$w$, we show that there exists some small~$\eps=\eps_{u,v}>0$ such that we have the \emph{integration-by-parts formula}
\begin{equation*}
\int_\msP \grad_w u \cdot v \diff\DF_\mssm=-\int_\msP u\cdot \grad_w v  \diff\DF_\mssm - \int_\msP u \cdot v\cdot \mbfB_\eps[w] \diff\DF_\mssm \comma
\end{equation*}
Theorem~\ref{l:IbP}, where
\begin{align}\label{eq:Intro:Drift}
\mbfB_\eps[w](\eta)\eqdef \sum_{x: \eta_x> \eps} \div^\mssm_x w \fstop
\end{align}

Provided that~$w\mapsto \grad_w u(\mu)$ be a $T^\Der_\mu\msP_2$-continuous linear functional, a gradient
\begin{align}\label{eq:Intro:Grad}
\mu\longmapsto\grad u(\mu)\in T^\Der_\mu\msP_2
\end{align}
is induced by Riesz Representation Theorem. The latter integration-by-parts formula is then a main tool in establishing the following theorem.

\begin{thm}[See Thm.~\ref{t:DForm} and Cor.~{\ref{c:Regularity}}]
The quadratic form~$(\mcE,\TF{})$ defined by
\begin{align*}
\mcE(u,v)\eqdef \int_\msP \tscalar{\grad u(\mu)}{\grad v(\mu)}_{T_\eta\msP_2} \diff\DF_\mssm(\mu) \comma \qquad u,v\in\TF{}
\end{align*}
is closable. Its closure~$(\mcE,\dom{\mcE})$ is a regular strongly local recurrent Dirichlet form on $L^2(\msP_2,\DF_\mssm)$ with carr\'e du champ operator
\begin{align}\label{eq:Intro:CdCRed}
\boldGamma(u,v)(\mu)\eqdef \tscalar{\grad u(\mu)}{\grad v(\mu)}_{T_\mu\msP_2}\comma \qquad u,v\in \TF{}\fstop
\end{align}
\end{thm}

\subsection{Generator} In addition to the carr\'e du champ operator, the generator of~$(\mcE,\dom{\mcE})$ entails further geometrical information. Up to Friedrichs extension,
\begin{align*}
\mbfL u=\mbfL_1u+\mbfL_2 u\comma \qquad u\in\hTF{}{0}\comma
\end{align*}
where
\begin{equation}\label{eq:Intro:Generator}
\begin{aligned}
\mbfL_1 u(\eta)\eqdef& \tfrac{1}{2} \sum_{i,j}^k (\partial_{ij}^2 F)\ttonde{\hat f_1^\trid(\eta),\dotsc, \hat f_k^\trid(\eta)} \int_\man\diff \eta(x)\, \varrho_i(\eta_x) \varrho_j(\eta_x)\, \gscal{\nabla_x f_i}{\nabla_x f_j}\comma
\\
\mbfL_2u(\eta)\eqdef& \tfrac{1}{2} \sum_i^k (\partial_i F)\ttonde{\hat f_1^\trid(\eta),\dotsc, \hat f_k^\trid(\eta)} \sum_{x: \eta_x>0} \varrho_i(\eta_x) \, \Delta^\mssg_x f_i \fstop
\end{aligned}
\end{equation}

For functions in the core~$\TF{}$ of~$(\mcE,\dom{\mcE})$, the first operator takes the form
\begin{align*}
\mbfL_1 u(\eta)=& \tfrac{1}{2}\sum_{i,j}^k (\partial_{ij}^2 F)\ttonde{\eta f_1,\dotsc,\eta f_k}\cdot  \scalar{\nabla f_i}{\nabla f_j}_{T_\eta\msP_2}\comma
\end{align*}
the \emph{diffusion part} of the generator.
The first order operator~$\mbfL_2$ represents instead the \emph{drift part} of the generator, constraining the process~$\Eta_\bullet$ on~$\Ppaiso$. We note here that the expression of~$\mbfL_2$ in~\eqref{eq:Intro:Generator} does \emph{not} converge for functions in~$\TF{}$, that is, in the pointwise limit~$\varrho_i\rar \car_I$. This is consistent with the heuristic observation of N.~Gigli that the Laplacian of potential energies on~$\msP_2$ should \emph{not} exist, see~\cite[Rmk.~5.6]{Gig12}. On the other hand though, this does not prevent the closability of~$(\mcE,\TF{})$ above.

This seeming contradiction is resolved in the understanding that the operator~$\mbfL_2$ is in fact a ``\emph{boundary term}'', and, as such, it was not accounted for in~\cite[\emph{ibid.}]{Gig12}. Indeed ---~in the present framework~--- the set~$\Ppaiso$ whereon~$\DF_\mssm$ is concentrated ought to be thought of as part of the \emph{geodesic boundary} of~$\msP_2$. Here, we say that a point~$\mu_0\in\msP$ is a geodesic-boundary point if there exists some $W_2$-geodesic~$\seq{\mu_t}_t$ for which~$\mu_0$ is extremal, that is,~$\seq{\mu_t}_t$ may not be further prolonged \emph{through}~$\mu_0$. The fact that measures with atoms satisfy this property is a consequence of the same result for Dirac masses, originally proved by J.~Bertrand and~B.~R.~Kloeckner in~\cite[Lem.~2.2]{BerKlo15}, and of the known fact that transport optimality is inherited by restrictions,~\cite[Thm.~4.6]{Vil09}.

\subsection{Quasi-invariance, representations and Helmoltz decomposition}
If~$G$ is a group acting measurably on a probability space~$(\Omega,\mcF,\mbfP)$, write~$G\acts \Omega$, we say that~$\mbfP$ is \emph{quasi-invariant} with respect to the action of~$h\in G$, write~$h.\omega$, if
\begin{align*}
\mbfP^h\eqdef (h.)_\pfwd \mbfP= \mbfR[h]\cdot \mbfP
\end{align*}
for some $\mcF$-measurable Radon--Nikod\'ym derivative~$\mbfR[h]\colon \Omega\rar [0,\infty]$; \emph{invariant} if~$\mbfP^h=\mbfP$.

\smallskip

In the case when~$\man=\mbbS^1$ and~$G$ is the Virasoro group~$\Diff^\infty_+(\mbbS^1)$ of smooth orientation-preserving diffeomorphisms of~$\mbbS^1$, the quasi-invariance of the entropic measure~$\mbbP^\beta$ and of the Dirichlet--Ferguson measure~$\DF$ has been a key tool in establishing the closability of the form~\eqref{eq:Intro:DirForm} in~\cite{vReStu09,Sha11}.
Let us recall the definition of the actions in~\cite{vReStu09,Sha11,vReYorZam08}.
Following~\cite[\S2.2]{vReStu09}, set
\begin{align*}
\msG(\R)\eqdef\set{\text{right-cont., non-decreasing } g\colon \R\rar \R\comma : g(x+1)=g(x)+1} \fstop
\end{align*}
Let~$\pr^{\mbbS^1}\colon \R\rar \mbbS^1\cong \R/\!\Z$ be the quotient projection, and $\msG(\mbbS^1)\eqdef \pr^{\mbbS^1}(\msG(\R))$.
By equi-variance, $g\colon \mbbS^1\rar \mbbS^1$ for every~$g\in\msG(\mbbS^1)$ and the set~$\msG(\mbbS^1)$, endowed with the usual composition of functions, is a semi-group with identity~$\id_{\mbbS^1}$.
In particular, the group~$\Diff^\infty_+(\mbbS^1)$ injects into~$\msG(\mbbS^1)$, e.g.,~\cite{vReStu09} or~\cite{LzDS19b}. Again following~\cite{vReStu09}, let
\begin{align}\label{eq:Intro:G1}
\msG_1\eqdef \quotient{\msG(\mbbS^1)}{\mbbS^1}\comma
\end{align}
where~$g,h\in\msG(\mbbS^1)$ are identified if~$g(\emparg)=h(\emparg+a)$ for some~$a\in\mbbS^1$,
and define the maps
\begin{equation}\label{eq:Intro:ZetaChi}
\begin{aligned}
\zeta\colon \msG(\mbbS^1)&\longrar\msP(\mbbS^1)
\\
\zeta\colon g&\longmapsto \diff g
\end{aligned}\qquad \textrm{and} \qquad 
\begin{aligned}
\chi\colon \msG_1&\longrar\msP(\mbbS^1)
\\
\chi\colon g&\longmapsto g_\pfwd \n\mssm
\end{aligned} \comma
\end{equation}
where~$\diff g$ is the Lebesgue--Stieltjes measure induced by~$g$ and~$\n\mssm$ denotes here the normalized Lebesgue measure on~$\mbbS^1$. Both maps are invertible. Namely, the inverse~$\zeta^{-1}$ assigns to~$\mu$ its cumulative distribution function, while~$\chi^{-1}$ assigns to~$\mu$ its generalized inverse distribution function; see~\cite[Eqn.~(2.2)]{vReStu09}. In particular, up to passing from~$\mbbS^1$ to~$I$,
\begin{align}\label{eq:ChiZeta}
\chi^{-1}=\emparg^{-1}\circ\zeta^{-1}\comma
\end{align}
where~$\emparg^{-1}\colon g\mapsto g^{-1}$ is the right-inversion map defined by
\begin{align*}
g^{-1}(t_0)\eqdef \inf\set{t\in I : g(t)>t_0}\fstop
\end{align*}

For~$h\in \Diff^\infty_+(\mbbS^1)$ we consider the left and right action on~$\msG(\mbbS^1)$ defined by
\begin{subequations}
\begin{align}
\addtocounter{equation}{11}
\label{eq:LeftAct}\tag{$\theparentequation_\ell$}
\ell_h\colon g&\longmapsto h\circ g \comma
\\
\addtocounter{equation}{5}
\label{eq:RightAct}\tag{$\theparentequation_r$}
r_h\colon g&\longmapsto g\circ h \fstop
\end{align}
\end{subequations}

It is then the content of~\cite[Thm.~4.1]{vReStu09} and~\cite[Thm.~4.1]{vReYorZam08} that the measure on~$\msG(\mbbS^1)$ defined as~$\mbbQ^\beta\eqdef \zeta^{-1}_\pfwd \DF_{\beta\n\mssm}$
is quasi-invariant with respect to the \emph{left} action~\eqref{eq:LeftAct}. This fact has two consequences.
On the one hand, the measure~$\DF_\mssm$ is quasi-invariant with respect to the ``left'' action~$L_h\eqdef \zeta\circ \ell_h \circ \zeta^{-1}$ of~$\Diff^\infty_+(\mbbS^1)$ on~$\msP$ corresponding to~\eqref{eq:LeftAct} on~$\msG(\mbbS^1)$ via~$\zeta$; see~\cite[Thm.~3.4]{Sha11}.
On the other hand, the entropic measure~$\mbbP^\beta=\chi_\pfwd \mbbQ^\beta$ is quasi-invariant with respect to the ``right'' ---~because of~\eqref{eq:ChiZeta}~--- action~$R_h\eqdef \chi\circ r_h\circ \chi^{-1}$ of~$\Diff^\infty_+(\mbbS^1)$ on~$\msP$ corresponding to~\eqref{eq:RightAct} on~$\msG(\mbbS^1)$ via~$\chi$; see~\cite[Cor.~4.2]{vReStu09}.

The action~\eqref{eq:LeftAct} is meaningful only for one-dimensional base spaces, where the representation of~$\mu$ via its cumulative distribution function is available. As a consequence, it is not possible to generalize the results of~\cite{Sha11} to base spaces of arbitrary dimension.
Analogously, since the $R_h$-quasi-invariance of the entropic measure~$\mbbP^\beta$ is a consequence of the $\ell_h$-quasi-invariance of~$\mbbQ^\beta$, it is bound to hold only in the case of one-dimensional base spaces.

Notwithstanding this fact, let us note that
\begin{align*}
g_{h_\pfwd \mu}\eqdef\zeta^{-1}(h_\pfwd \mu)=g_\mu\circ h \defeq r_h(g_\mu)\comma
\end{align*}
thus, the action~$K_h\eqdef \zeta\circ r_h\circ \zeta^{-1}$ of~$\Diff^\infty_+(\mbbS^1)$ on~$\msP(\mbbS^1)$ is meaningful in the general case, as we detail now.
Indeed, let~$\mfG\eqdef \Diff^\infty_+(\man)$ be the Lie group of orientation-preserving smooth diffeomorphisms of~$\man$. The natural action of~$\mfG$ on~$\man$ lifts to an action of~$\mfG$ on~$\msP$, given by
\begin{equation}\label{eq:Intro:GAct}
\begin{aligned}
.\colon \mfG \times \msP&\longrar \msP
\\
(\uppsi\;, \; \mu)\,&\longmapsto \uppsi_\pfwd \mu
\end{aligned} \fstop
\end{equation}

The quasi-invariance of~$\dvol_{\msP_2}$ with respect to the action~$\mfG\acts \msP$ is a natural question within representation theory, e.g.,~\cite{KonLytVer15, AirMal06, VerGelGra75}, where it corresponds to the action above defining a quasi-regular representation of the infinite-dimensional Lie group~$\mfG$ on~$L^2(\msP)$.
In turn, this relates to the closability of the gradient~\eqref{eq:Intro:Grad} on~$\msP_2$. Indeed, the Lie algebra of~$\mfG$ is the algebra~$\Vect^\infty\eqdef \Gamma^\infty(T\man)$ of smooth vector fields on~$\man$ and its exponential curves based at~$\id_\mfG=\id_\man$ are precisely the shifts~$\fl^{w,t}$ defining the directional derivative~\eqref{eq:Intro:DirDer}.

It turns out that the Dirichlet--Ferguson measure~$\DF_\mssm$ is \emph{not} quasi-invariant with respect to the action of~$\mfG$\,: Were this the case, then the Gamma measure~$\GP_\mssm=\DF_\mssm\otimes \Gam[1,\beta]$ too would be quasi-invariant with respect to the analogous action~$\mfG \acts \Mbp(\man)=\msP(\man)\times\R_+$. However, this does not hold; see the introduction to~\S2.4 in~\cite{KonLytVer15}.

In order to address this issue, we recall the following definition from~\cite[Dfn.~9]{KonLytVer15}. As in the beginning of this section, let~$G$ be a group acting on a probability space~$(\Omega,\mcF,\mbfP)$.

\begin{defs}[Partial quasi-invariance]\label{d:PQI} $\mbfP$ is \emph{partially quasi-invariant} with respect to $G\acts \Omega$ if there exists a filtration $\mcF_\bullet\eqdef\seq{\mcF_n}_n$ so that:
\begin{enumerate}[label=\bfseries({\itshape\roman*})]
\item\label{i:d:PQI1} $\mcF=\mcF_\infty$, the $\sigma$-algebra generated by~$\mcF_\bullet$;
\item\label{i:d:PQI2} for each $h\in G$ and~$n\in\N$ there exists~$n'\in\N$ so that~$h.\mcF_n=\mcF_{n'}$;
\item\label{i:d:PQI3} for each $h\in G$ and~$n\in \N$ there exists an $\mcF_n$-measurable~$\mbfR_n[h]\colon \Omega\rar [0,\infty]$ so that
\begin{align*}
\int_\Omega u(\omega) \diff \mbfP^h(\omega) = \int_\Omega u(\omega) \mbfR_n[h](\omega) \diff \mbfP(\omega)
\end{align*}
for each $\mcF_n$-measurable semi-bounded $u\colon \Omega \rar [-\infty,\infty]$.
\end{enumerate}

If~$\mbfP$ is quasi-invariant with respect to~$G\acts \Omega$, then it is partially quasi-invariant for the choice~$\mcF_n=\mcF$. Finally,~$\mbfR_n[h]$ is $\mbfP$-a.e.~uniquely defined; see~\cite[Rmk.~10]{KonLytVer15}.
\end{defs}

Now, let~$\Bo{\bullet}(\msP)\eqdef \tseq{\Bo{\eps}(\msP)}_{\eps\in I}$ be the filtration of~$\sigma$-algebras on~$\msP_2$ generated by the functions
\begin{align}\label{eq:Intro:RND}
\mbfR_\eps[\uppsi]\colon \eta \longmapsto \prod_{x: \eta_x>\eps} \frac{\diff \uppsi_\pfwd \mssm}{\diff \mssm}(x) \comma \qquad \uppsi\in \Diff^\infty_+(\man) \fstop
\end{align}

Then,~$\Bo{1}(\msP)$ is the trivial $\sigma$-algebra and the restriction~$\Bo{0}(\msP)_{\Ppa}$ of~$\Bo{0}(\msP)$ to~$\Ppa$ coincides with the Borel $\sigma$-algebra~$\Bo{}(\msP)_{\Ppa}$, Lemma~\ref{l:SAlg}. We shall prove the following

\begin{thm}[See Prop.~\ref{p:PQI} and~Cor.~\ref{c:Martingale}]
Let~$\uppsi\in \Diff^\infty_+(\man)$. Then,
\begin{enumerate*}[label=$\boldsymbol{(}$\bfseries{\itshape\roman*}$\boldsymbol{)}$]
\item $\DF_\mssm$ is partially quasi-invariant w.r.t.\ the action of~$\uppsi$ on the filtration~$\seq{\Bo{1/n}(\msP)}_n$;

\item $\DF_\mssm$ is quasi-invariant w.r.t. the action of~$\uppsi$ if and only if~$\uppsi_\pfwd \mssm=\mssm$, in which case it is in fact invariant;

\item if~$\fl^{w,t}$ is the flow of a smooth vector field~$w$, then~$\mbfB_\bullet$, defined in~\eqref{eq:Intro:Drift}, satisfies
\end{enumerate*}
\begin{align*}
\mbfB_\bullet[w]= \diff_t\restr_{t=0}\mbfR_\bullet[\fl^{w,t}] \comma
\end{align*}
and is a centered square-integrable $\DF_\mssm$-martingale adapted to~$\Bo{\bullet}(\msP)$.
\end{thm}

By the theorem, the algebra~$\Vect$ is decomposed, as a vector space, into a direct sum~$\Vect^\inv\oplus \Vect^\pqi$, where~$\Vect^\inv$, resp.~$\Vect^\pqi$, denotes the space of vectors such that~$\DF_\mssm$ is invariant, resp.\ partially quasi-invariant not quasi-invariant, with respect to the action of~$\fl^{w,t}$.
Now, it is readily checked that, if~$\uppsi=\fl^{w,1}$ for some~$w\in \Vect$, then~$\uppsi$ is $\mssm$-measure-preserving, i.e.~$\uppsi_\pfwd \mssm=\mssm$, if and only if~$w$ is divergence-free; see e.g.,~\cite[Rmk.~1.29]{AmbGig11}.
Thus,~$\Vect^\inv=\Vect^\div$, the space of divergence-free vector fields, and~$\Vect^\pqi=\Vect^\nabla$ the space of gradient-type vector fields. This is but an instance of the classical \emph{Helmholtz decomposition}, and extends for every~$\eta$ to an orthogonal decomposition of the tangent space~$T^\Der_\eta\msP_2$ into the subspaces~$T_\eta\msP_2=\cl_{T^\Der_\eta \msP_2}(\Vect^\nabla)$ and~$\cl_{T^\Der_\eta\msP_2}(\Vect^\div)$; also cf.~\cite[Prop.~1.28]{AmbGig11}.

\subsection{Properties of the process}
 By the standard theory of Dirichlet forms there is a Markov process~$\Eta_\bullet$ with state space~$\msP$ properly associated to~$(\mcE,\dom{\mcE})$ in the sense of~\cite[Dfn.~IV.2.5(i)]{MaRoe92}. In order to show, as anticipated, that~$\Eta_\bullet=\boldPhi\circ\widehat\BM_\bullet$, we shall construct finite-dimensional approximations of~$\Eta_\bullet$ and~$\widehat\BM_\bullet$ and prove their coincidence up to a suitable restriction of the map~$\boldPhi$. Namely, we construct
\begin{itemize}
\item a sequence of Dirichlet forms~$(\mcE^n,\dom{\mcE^n})$ defined as a martingale-type approximation of~$(\mcE,\dom{\mcE})$ w.r.t.~the filtration $\tseq{\Bo{1/n}(\msP)}_n$ given by~\eqref{eq:Intro:RND};

\item a sequence of Dirichlet form~$(\hat\mssE^n,\dom{\hat\mssE^n})$ associated to the processes~$\widehat\mssW^n_\bullet$ obtained by truncation of~$\widehat\BM_\bullet$ onto the first $n$ components of the product space~$\mbfM$ and onto the first $n$ elements of~$\mbfs\in\To$; see Prop.~\ref{p:ProdForm}.
\end{itemize}

We show their coincidence and their generalized Mosco convergence to~$(\mcE,\dom{\mcE})$ in the sense of Kuwae--Shioya~\cite{KuwShi03}; see~Prop.~\ref{p:MoscoEn}.
This approximation allows to identify, up to quasi-homeomorphism of Dirichlet forms~\cite{CheMaRoe94}, the Dirichlet form~$\mcE$ with the Dirichlet form~$\widehat\dirE$ associated to~$\widehat\BM_\bullet$, hence to prove~$\Eta_\bullet$'s sample-continuity properties and to classify its invariant sets and invariant measures, Thm.~\ref{t:PropertiesEta}.

Profiting the essential self-adjointness of the generator~$\mbfL$ on~$\hTF{}{0}$, Prop.~\ref{p:ESA}, we are able to show the $\DF_\mssm$-a.e.\ differentiability of~$W_2$-Lipschitz functions, Prop.~\ref{p:RademacherE}, and to provide a one-sided Varadhan-type estimate of the short-time asymptotics for the heat kernel of~$\mcE$, Cor.~\ref{c:Varadhan}.

\section{Comparison with other measure-valued processes}
Before exploring the construction and the properties of~$\mcE$ and its relation to~$\Eta_\bullet$, it is worth to compare its carr\'e du champ operator~\eqref{eq:Intro:CdCRed} and generator with those of other processes on~$\msP$. 

\subsection{A comparison with the Fleming--Viot process} 
Following L.~Overbeck, M.~R\"ockner~and B.~Schmuland, let
\begin{align}\label{eq:Intro:DeltaORS95}
\ttonde{\tfrac{\partial}{\partial \delta_x}u}(\mu)(x)\eqdef \diff_t\restr_{t=0} u(\mu+t\delta_x)\fstop
\end{align}
It was shown in~\cite{OveRoeSch95}, that the Dirichlet form~$(\mcE^{\mathsc{fv}},\dom{\mcE^{\mathsc{fv}}})$ with carr\'e du champ operator
\begin{align*}
\boldGamma^{\mathsc{fv}}(u)(\mu)\eqdef \Var_\mu\ttonde{\tfrac{\partial}{\partial \delta_\emparg} u(\mu)}\comma \qquad u\in \hTFB{}\comma
\end{align*}
and invariant measure~$\DF_\mssm$ is associated with the \emph{Fleming--Viot process}~\cite{FleVio79} with parent independent mutation.
In~\cite{Sha11}, J.~Shao observed that the increment in~\eqref{eq:Intro:DeltaORS95} is not internal to~$\msP$. To overcome the issue, he considered the map~$S_f$~\cite[Eqn.~(2.7)]{Sha11}, originally introduced by K.~Handa in~\cite{Han02},
\begin{align}\label{eq:ShiftHanda}
S_f(\mu)\eqdef \frac{e^f\cdot \mu}{\mu(e^f)}\comma \qquad f\in \mcC(\man)\comma
\end{align}
and termed `exponential map', see Remark~\ref{r:Shifts} below.

For~$\mu\in\msP$ we recall the notation~\eqref{eq:Intro:CoCo} and set
\begin{align*}
\ttonde{\widetilde{\tfrac{\partial}{\partial\delta_x}}}u(\mu)\eqdef \diff_t\restr_{t=0} u(\mu^x_t)\comma
\end{align*}
also cf.~\cite[Eqn.~(1.1)]{Sch97}.
Then,
\begin{align*}
\diff_t\restr_{t=0} S_{tf}(\mu)= \scalar{\ttonde{\widetilde{\tfrac{\partial}{\partial\delta_\emparg}}u}(\mu)}{f}_{L^2(\mu)} \comma \qquad u\in \TF{} \comma
\end{align*}
and
\begin{align*}
\boldGamma^{\mathsc{fv}}(u)(\mu)= \norm{\widetilde{\tfrac{\partial}{\partial\delta_\emparg}} u}_{L^2(\mu)}^2 \fstop
\end{align*}

As noted by M.~D\"oring and W.~Stannat in~\cite[Rmk.~1.5]{DoeSta09}, the carr\'e du champ operator~\eqref{eq:Intro:CdCRed} is strictly stronger than~$\boldGamma^{\mathsc{fv}}$ and one has in fact
\begin{align*}
\boldGamma(u)(\mu)=\norm{\nabla_\emparg \tfrac{\partial }{\partial\delta_\emparg} u}_{L^2(\mu)}^2=\norm{\nabla_\emparg \widetilde{\tfrac{\partial }{\partial\delta_\emparg}} u}_{L^2(\mu)}^2 \fstop
\end{align*}

In the case~$\man=\mbbS^1$, it holds that~$T_\mu\msP_2=T^\Der_\mu\msP_2=L^2(\mbbS^1,\mu;\R)$, and
\begin{align}\label{eq:Intro:CdCWD}
\boldGamma=\boldGamma^\vReStu\comma
\end{align}
the carr\'e du champ~\cite[Dfn.~7.24]{vReStu09} of the Wasserstein diffusion~\cite{vReStu09}. 
Letting $(\mcE^\vReStu,\dom{\mcE^\vReStu})$ be the Dirichlet form~\cite[Thm.~7.25]{vReStu09} of the Wasserstein diffusion,~\eqref{eq:Intro:CdCWD} is interpreted as follows. By definition~$u\in \TF{}\subset \dom{\mcE^\vReStu}$; then, for each~$u\in \TF{}$ there exist a continuous $\mbbP^\beta$-representative~$\widetilde{\boldGamma^\vReStu(u)}$ of~$\boldGamma^{\vReStu}(u)$ and a continuous $\DF_{\mssm}$-representative~$\widetilde{\boldGamma(u)}$ of~$\boldGamma(u)$ so that~$\widetilde{\boldGamma^\vReStu(u)}=\widetilde{\boldGamma(u)}$ everywhere on~$\msP(\mbbS^1)$.

\subsection{A comparison with the Wasserstein Diffusion}
For the purpose of further comparison, let us recall the form of the generator of the Wasserstein diffusion on~$\man=I$, or~$\mbbS^1$. We have
\begin{align*}
\mbfL^{\vReStu}=&\mbfL^{\vReStu}_1+\mbfL^{\vReStu}_2+\mbfL^{\vReStu}_3 && \text{over $\man = I$}\comma
\\
\mbfL^{\vReStu}=&\mbfL^{\vReStu}_1+\mbfL^{\vReStu}_2 && \text{over $\man = \mbbS^1$}\comma
\end{align*}
where~$\mbfL^\vReStu_1=\mbfL_1$ is the diffusion part of~$\mbfL^{\mathsc{wd}}$ in the decomposition~\cite[Thm.~7.25]{vReStu09}, and the equality is interpreted as in~\eqref{eq:Intro:CdCWD}. The drift part~$\mbfL^{\vReStu}_2+\mbfL^{\vReStu}_3$ is discussed below.
As noted in~\cite[Rmk.~7.18]{vReStu09},~$\mbfL_1^{\mathsc{wd}}$ \guillemotleft describes the [Wasserstein] diffusion [...] in all directions of the respective tangent spaces\guillemotright.
Thus, the process~$\Mu^\vReStu_\bullet$ associated with~$(\mcE^{\mathsc{wd}},\dom{\mcE^{\mathsc{wd}}})$ \guillemotleft experiences [...] the full tangential noise\guillemotright.
In the present case, the same statement may be formulated rigorously, in terms of Hino's \emph{index}~\cite{Hin09} of the form, Prop.~\ref{p:Hino}.

For simplicity, we only recall the form of~$\mbfL^{\vReStu}_2$,~$\mbfL^{\vReStu}_3$ on functions of the form~$f^\trid$, by means of the martingale problem corresponding to~$\mu^{\vReStu}_\bullet$. Denote by~$\scalar{\emparg}{\emparg}$ the standard pairing of distributions.
\begin{prop}[{\cite[Cor.~7.20]{vReStu09}}]
For each~$f\in\mcC^2(I)$ with $f'(0)=f'(1)=0$,
\begin{align}\label{eq:vReStuWD1}
M^f_t\eqdef \scalar{\mu^{\vReStu}_t}{f}-\int_0^t \scalar{\mbfL^{\vReStu}_2 \mu^{\vReStu}_s}{f} \diff s -\int_0^t \scalar{\mbfL^{\vReStu}_3 \mu^{\vReStu}_s}{f} \diff s
\end{align}
is a continuous martingale with quadratic variation process
\begin{align*}
\quadre{M^f}_t=\int_0^t \scalar{\mu^{\vReStu}_s}{(f')^2} \diff s \fstop
\end{align*}
\end{prop}

Here, for~$\mu$ singular continuous w.r.t.~$\Leb^1$, the operator~$\mbfL^\vReStu_2 \mu$ acts as the distribution
\begin{align*}
(\mbfL^\vReStu_2 \mu)f\eqdef\,& \sum_{J\in \mathrm{gaps}(\mu)}\quadre{\frac{f''(J_+)+f''(J_-)}{2}-\frac{f'(J_+)-f'(J_-)}{J_+-J_-}}\comma 
\end{align*}
where~$\mathrm{gaps}(\mu)$ is the set of maximal intervals~$J\eqdef (J_-,J_+)$ with~$\mu J=0$. The operator~$\mbfL^\vReStu_3$ acts as the distribution
\begin{align*}
\mbfL^\vReStu_3 \mu = \tfrac{\beta}{2}\mu''\comma
\end{align*}
a scalar multiple of the second distributional derivative.

When~$\man=\mbbS^1$, the operator~$\mbfL_2^{\mathsc{wd}}$ may be given the same interpretation of~$\mbfL_2$, analogously to the case of~$\mbfL_1$ and~$\mbfL_1^{\mathsc{wd}}$,~\cite[Rmk.~7.18, Thm.~7.25]{vReStu09}.
Finally, we note that the operator~$\mbfL_3^{\mathsc{wd}}$ in~\cite[Rmk.~7.18]{vReStu09} has no counterpart in our case ---~which should rather be compared with~\cite[Thm.~7.25]{vReStu09}~---, since~$\mbfL_3^{\mathsc{wd}}$ is an artifact of the boundary of~$I$.

\subsection{A comparison with the modified massive Arratia flow}
In~\cite{Kon17} V.~V.~Ko\-na\-rov\-skyi introduced the \emph{modified massive Arratia flow}, a random element~$y(\emparg,\emparg)$ in the Skorokhod space~$D\ttonde{I;\mcC([0,T])}$. The corresponding measure-valued process~$\mu^{\KonvRe}_t\eqdef y(\emparg,t)_\pfwd \Leb^1$ is studied by Konarovskyi and M.-K.\ von Renesse in~\cite{KonvRe18}.
Its generator has the form
\begin{align*}
\mbfL^{\KonvRe}=\mbfL^{\KonvRe}_1+\mbfL^{\KonvRe}_2\comma
\end{align*}
where~$\mbfL^{\KonvRe}_1=\mbfL_1$ is the diffusion part, and the equality is interpreted as in~\eqref{eq:Intro:CdCWD}. The drift part~$\mbfL^{\KonvRe}_2$ is discussed below.
The process is a solution to the following martingale problem.

\begin{prop}[{\cite[Prop.~1.2]{KonvRe18}}] For each~$f\in\Cb^2(\R)$,
\begin{align}\label{eq:KonvReMMAF}
M^f_t\eqdef \scalar{\mu^{\KonvRe}_t}{f}-\int_0^t \scalar{\mbfL^{\KonvRe}_2 \mu^{\KonvRe}_s}{f} \diff s
\end{align}
is a continuous local martingale with quadratic variation process
\begin{align*}
\quadre{M^f}_t=\int_0^t \scalar{\mu^{\KonvRe}_s}{(f')^2} \diff s \comma
\end{align*}
where, for~$\mu\in\msP(\R)$ with finitely many atoms, the operator~$\mbfL_2^{\KonvRe}\mu$ acts as the distribution
\begin{align}\label{eq:IntroL2}
\mbfL^\KonvRe_2 \mu= \tfrac{1}{2}\sum_{x\,:\, \mu\!\set{x}>0} \delta_x''\fstop
\end{align}
\end{prop}

For the relation of~\eqref{eq:vReStuWD1} and~\eqref{eq:KonvReMMAF} with the Dean--Kawasaki dynamics for supercooled liquid models, cf., e.g.,~\cite[\S1.1]{KonvRe17}. For the relation with the Dawson--Watanabe process, the Dobrushin--Doob process, and the empirical-measure process of a single Brownian motion, see~\cite[Rmk.~7.21]{vReStu09}.

\smallskip

In the case when~$d\geq 2$, the process~$\Eta_\bullet$ may be regarded as a counterpart on multidimensional base spaces to the modified massive Arratia flow, as a consequence of the following multi-dimensional analogue of~\eqref{eq:KonvReMMAF}.
\begin{prop}[see Theorem~\ref{t:PropertiesEta}\iref{i:t:PropertiesEta:6}]\label{p:IntroSPDE}
Assume~$d\geq 2$. For each~$\hat f$ as in~\eqref{eq:Intro:fTrid},
\begin{align*}
M^{\hat f}_t\eqdef \hat f^\trid(\Eta_t)-\int_0^t \tscalar{\mbfL_2 \Eta_s}{\hat f^\trid} \diff s
\end{align*}
is a continuous martingale with quadratic variation process
\begin{align*}
\quadre{M^{\hat f}}_t=\int_0^t \ttonde{\ttabs{\nabla \hat f}^2}^\trid \Eta_s \diff s \comma
\end{align*}
where, for~$\mu\in\msP(\man)$, the operator~$\mbfL_2\mu$ acts as the distribution
\begin{align*}
\mbfL_2\mu=\tfrac{1}{2}\sum_{x\,:\, \mu\!\set{x}>0} \Delta^\mssg\, \delta_x
\end{align*}
and~$\Delta^\mssg$ is the distributional Laplace--Beltrami operator on~$\man$.
\end{prop}

\section{Product spaces and spaces of measures}\label{s:Prelim}
Let~$S$ be any Hausdorff topological space with topology~$\tau$ and Borel $\sigma$-algebra~$\Bo{\T}$.
We denote by~$\mcC(S)$, resp.\ $\Cc(S)$, $\mcC_0(S)$, $\Cb(S)$, $\mcB_b(S)$, the space of (real-valued) continuous, resp.\ continuous compactly supported, continuous vanishing at infinity, continuous bounded, bounded (Borel) measurable, functions on~$(S,\tau)$. Whenever~$U\in\tau$, we always regard the spaces~$\Cc(U)$ and~$\mcC_0(U)$ as embedded into~$\mcC(S)$ by taking the trivial extension of~$f\in\Cc(U)$ identically vanishing on~$U^\complement\eqdef S\setminus U$.

For~$n\in\N$ and~$\mbfh\colon S\rar \R^k$,~$\mbfh\eqdef\seq{h_1,\dotsc, h_k}$, we set~$\displaystyle{\norm{\mbfh}_\infty\eqdef \sup_{s\in S} \max_{i\leq k}\abs{h_i(s)}}$.

\subsection{Riemannian manifolds}
The main object of our analysis are Riemannian manifolds satisfying Assumption~\ref{ass:RM}. We refer the reader to the monograph~\cite{Gri09} for a detailed account of (stochastic) analysis on manifolds. We state here without proof the main results we shall assume in the following.

For~$w\in\Vect^\infty$, the algebra of smooth vector fields on~$\man$, we denote by~$\fl^{w,t}$ its flow, satisfying
\begin{align}\label{eq:FlowDef}
\begin{aligned}
\diff_t \fl^{w,t}(x)=&w(\fl^{w,t}(x))\\
\fl^{w,0}(x)=&x
\end{aligned}\comma
\qquad x\in \man\comma t\in \R \fstop
\end{align}

As a consequence of the compactness of~$\man$:
\begin{enumerate*}[label=$\boldsymbol{(}$\bfseries{\itshape\alph*}$\boldsymbol{)}$]
\item for every~$w\in \Vect^\infty$ the flow~$\fl^{w,t}$ is well-defined and a smooth diffeomorphism for every~$t\in \R$, with inverse $(\fl^{w,t})^{-1}$ identical to~$\fl^{w,-t}$;
\item the manifold~$\man$ is geodesically complete, that is, every geodesic curve is infinitely prolongable to a locally length-minimizing curve;
\item the Laplace--Beltrami operator~$\Delta^\mssg$ is a densely defined linear operator on~$L^2(\man,\mssm)$, essentially self-adjoint on~$\mcC^\infty(\man)$ and with discrete spectrum;

\item the heat semigroup~$\mssH_t\eqdef e^{-t\Delta^\mssg}$ on $L^2(\man,\n\mssm)$ has (absolutely continuous) kernel with density $y\mapsto\mssh_t(\emparg,y)$;

\item the manifold~$\man$ is stochastically complete, that is
\end{enumerate*}
\begin{align*}
(\mssH_t\car_\man)(x)= \int_\man \mssh_t(x,y) \diff\n\mssm(y)=1\comma \qquad x\in \man\comma t>0 \fstop
\end{align*}

For~$f\in\mcC^1(\man)$ and~$w\in\Vect^\infty$ let
\begin{enumerate*}[label=$\boldsymbol{(}$\bfseries{\itshape\alph*}$\boldsymbol{)}$]
\item $\nabla^\mssg_w f=(\diff f)w$ be the directional derivative of~$f$ in the direction~$w$; 
\item $\nabla^\mssg f$ be the gradient of~$f$; 
\item $\div^\mssm w$ be the divergence of~$w$ induced by the volume measure~$\mssm$, satisfying the integration-by-parts formula
\end{enumerate*}
\begin{align*}
\int_\man (\nabla^\mssg_w f_1) \cdot f_2 \diff\mssm = -\int_\man f_1\cdot (\nabla^\mssg_w f_2) \diff\mssm -\int_\man f_1\cdot f_2\cdot \div^\mssm w \diff\mssm\fstop
\end{align*}

Whenever no confusion may arise, we drop the superscript~$\mssg$ from the notation. We denote variables by a superscript: e.g.~$\nabla^z\restr_{z=x} f$ denotes the gradient of~$f$ in the variable~$z$ computed at the point~$x\in \man$.

\paragraph{Canonical Dirichlet forms}
We endow~$(\man,\mssg)$ with its \emph{canonical Dirichlet form} $(\mssE^\mssg,\dom{\mssE^\mssg})$, defined as the closure of the pre-Dirichlet form
\begin{align}\label{eq:DirFormX}
\mssE^\mssg(f_1,f_2)\eqdef \int_\man \Gamma^\mssg(f_1, f_2) \diff\n\mssm \comma\qquad f_i\in\mcC^\infty(\man) \comma
\end{align}
where~$\Gamma^\mssg$ is the carr\'e du champ operator~$\Gamma^\mssg(f_1,f_2)\eqdef  \tfrac{1}{2}\gscal{\nabla^\mssg f_1}{\nabla^\mssg f_2}$.
We stress that the reference measure here is the \emph{normalized} volume~$\n\mssm$ (as opposed to the volume measure~$\mssm$) and that we adhere to the stochastic convention, taking~$\tfrac{1}{2}\Delta^\mssg$ as generator of~$\mssE^\mssg$ (as opposed to~$\Delta^\mssg$). If not otherwise stated, by a \emph{Brownian motion on~$\man$} we shall mean the diffusion process associated to~$\mssE^\mssg$; due to the normalization in the measure, this differs from the usual Brownian motion by a linear deterministic time change.

\paragraph{Conformal rescaling} We will make extensive use of conformal rescaling for metric objects on~$\man$, some of which are listed below. Let~$a>0$. Then,
\begin{equation}\label{eq:Rescal}
\begin{aligned}
a\nabla^{a\mssg}=&\nabla^\mssg\comma & a\Delta^{a\mssg}=&\Delta^\mssg\comma &\Gamma^\mssg(\emparg)\eqdef& \abs{\nabla^\mssg\emparg}_{\mssg}^2=a\Gamma^{a\mssg}(\emparg) \comma
\\
\mssm^{a\mssg}=&a^{d/2} \mssm^\mssg \comma & \mssh^{a\mssg}_t=&\mssh^\mssg_{t/a}\comma &\mssH^{a\mssg}_t\eqdef& e^{-t\Delta^{a\mssg}}=\mssH^\mssg_{t/a} \fstop
\end{aligned}
\end{equation}

\subsection{Product manifolds}\label{ss:ProdMan}
Everywhere in the following let~$\mbfM$, resp.~$\mbfI$, be the infinite-product manifold of~$\man$, resp.~$I$, endowed with the respective product topologies, and~$\Tso\subset \To\subset \boldDelta\subset \mbfI$ as in~\eqref{eq:DeltaTau} be endowed with the trace topology.
Define the topology~$\T_\mrmp$ on~$\hM\eqdef\To\times\mbfM$ as the product topology of the spaces~$\mbfM$ and~$\To$ and denote by the same symbol the trace topology on any subset of~$\hM$. For the sake of notational simplicity we let further~$\andihM\eqdef \Tso\times\andi{\mbfM}$.
We always endow~$\hM$ with the fully supported measure~$\hvolm_\beta\eqdef \PD_\beta\otimes\n\volm$, concentrated on~$\andihM$.

\begin{rem}\label{r:Notation} We adopt the following conventions:
\begin{enumerate*}[label=$\boldsymbol{(}$\bfseries{\itshape\alph*}$\boldsymbol{)}$]
\item {\bfseries bold}face symbols always denote infinite-product objects or objects on spaces of measures;
\item for an object~$\square$ on~$\man$, resp.~$\mbfM$, the symbol~$\widehat\square$ denotes some extension of~$\square$ to~$\hman\eqdef \man\times [0,1]$, resp.~$\hM$. 
\item for an object~$\square$ with infinitely many coordinates, the symbol~$\andi{\square}$ indicates that such coordinates are always assumed to be pairwise different.
\end{enumerate*}
\end{rem}

\begin{defs}\label{d:NProd}
For~$n\in \N$ and~$\mbfs\in\Tso$, we denote by~$\man^{n,\mbfs}$ the product manifold~$\man^{\times n}$ endowed with Riemannian metric~$\mssg^{n,\mbfs}\eqdef \bigoplus_\ell^n s_\ell\, \mssg$, normalized volume measure~$\n\mssm^n\eqdef \n\mssm^{\otimes n}$, canonical form~$(\mssE^{n,\mbfs},\dom{\mssE^{n,\mbfs}})$, heat semigroup~$\mssH^{n,\mbfs}_t$ with kernel~$\mssh^{n,\mbfs}_t$, and Brownian motion
\begin{align}\label{eq:BMn}
\mssW^{n,\mbfs}_\bullet\eqdef \seq{\Omega,\mcF,\seq{\mcF_t}_{t\geq 0}, \seq{\mssW^{n,\mbfs}_t}_{t\geq 0}, \seq{P^{n,\mbfs}_\mbfx}_{\mbfx\in \man^{n,\mbfs}}} \comma
\end{align}
where~$\mcF_\bullet$ is the natural filtration in the sense of~\cite[Dfn.~IV.1.8]{MaRoe92}.
For~$\mbfx_0\eqdef \ttseq{x_0^\ell}_{\ell\leq n}\in \man^{\times n}$ and~$t>0$ one has~$\mssW^{n,\mbfs;\mbfx_0}_t=\tseq{x^1_{t/s_1},\dotsc, x^n_{t/s_n}}$, where~$\ttseq{x^\ell_\bullet}_{\ell\leq n}$ are independent Brownian motions on~$\man$ respectively starting at~$x^\ell_0$.

Denote further by~$\mbfM^\mbfs$ the infinite-product manifold~$\mbfM$ endowed with the symmetric tensor field~$\mssg^\mbfs\eqdef \bigoplus_\ell^\infty s_\ell\,\mssg$ and normalized volume measure~$\n\volm\eqdef \n\mssm^{\otimes\infty}$. Each of the above objects is well-defined since~$\mbfs\in \Tso$, as opposed to~$\To$.
\end{defs}

We summarize several results about the canonical form and heat semigroup (kernel) on~$\mbfM^\mbfs$. For the sake of notational simplicity,~$A\eqdef\set{\ell_1,\dotsc, \ell_{\abs{A}}}$ shall always denote a \emph{finite} subset of~$\N$, and we set~$\mbfx_A\eqdef \seq{x_{\ell_1},\dotsc, x_{\ell_{\abs{A}}}}$. 
For~$m\in \N$, define the algebra of cyl\-in\-der functions~\cite[Eqn.~(3)]{AlbDalKon00}
\begin{align}\label{eq:Cyl}
\Cyl{m}\eqdef\set{\begin{matrix} u\colon \mbfM\longrar \R : u(\mbfx)= F(\mbfx_A)\comma F\in \mcC^m(\man^{\abs{A}}) \end{matrix}}\fstop
\end{align}

\begin{thm}[{\cite{AlbDalKon00, BenSaC97}}]\label{t:ADKBSC}
Fix~$\mbfs\in\Tso$. Then,
\begin{enumerate*}[label=$\boldsymbol{(}$\bfseries{\itshape\roman*}$\boldsymbol{)}$]
\item\label{i:t:ADKBSC:1} $\mbfM^\mbfs$ is a Banach manifold modelled on the space~$\ell^\infty(\N;\R^d)$ with norm $\norm{\mbfa}_\mbfs\eqdef \sup_\ell s_\ell \abs{a_\ell}_{\R^d}$.
\item\label{i:t:ADKBSC:2} The form~$(\dirE^\mbfs,\Cyl{1})$ given by
\end{enumerate*}
\begin{align}\label{eq:ADKBSC0}
\dirE^\mbfs(u)\eqdef \tfrac{1}{2}\int_\mbfM \sum_{\ell\in A} s_\ell^{-1}\tabs{\nabla^{\mssg,x_\ell}F(\mbfx_A)}_{\mssg_{x_\ell}}^2 \diff\n\volm(\mbfx)\comma \qquad u(\mbfx)= F(\mbfx_A)\in\Cyl{1}\comma
\end{align}
is closable, cf.~\cite[Eqn.~(25)]{AlbDalKon00}.
\begin{enumerate*}[resume*]
\item\label{i:t:ADKBSC:3} Its closure is a regular strongly local Dirichlet form on~$L^2(\mbfM,\n\volm)$ with standard core~$\Cyl{\infty}$;
\item\label{i:t:ADKBSC:4} generator~$(\boldDelta^\mbfs,\dom{\boldDelta^\mbfs})$, essentially self-adjoint on~$\Cyl{2}$, given by
\end{enumerate*}
\begin{align*}
\boldDelta^\mbfs u (\mbfx)=\tfrac{1}{2}\sum_{\ell\in A} s_\ell^{-1} \Delta^{\mssg,x_\ell} F(\mbfx_A)\comma \qquad u(\mbfx)= F(\mbfx_A)\in\Cyl{2}\semicolon
\end{align*}
and 
\begin{enumerate*}[resume*]
\item\label{i:t:ADKBSC:5} heat kernel~$\heat^\mbfs_\bullet$, defined as in~\eqref{eq:Intro:heat}, absolutely continuous w.r.t.~$\n\volm$ with density in~$\Cb(\mbfM)$;
\item\label{i:t:ADKBSC:6} associated Brownian motion~$\BM^\mbfs_\bullet$, defined as in~\eqref{eq:Intro:BM}, satisfying
\end{enumerate*}
\begin{align*}
\Proj{n}{}\circ \BM^{\mbfs;\mbfx_0}_t=\mssW^{\mbfs,n;\Proj{n}{}(\mbfx_0)}_t\comma \qquad \mbfx_0\in \mbfM\comma \quad t>0\comma
\end{align*}
where~$\mssW^{n,\mbfs;\Proj{n}{}(\mbfx_0)}_t$ is as in Definition~\ref{d:NProd}.

\begin{proof} Throughout the proof we shall refer to results in~\cite{AlbDalKon00,AlbDalKon97} concerned with the infinite-product manifold~$\mbfM=\mbfM^{\uno}$, rather than with~$\mbfM^\mbfs$. However, as noted in~\cite[Rmk.~2.1]{AlbDalKon00}, this construction is possible and nearly identical for arbitrary~$\mbfs\in\To$. Again throughout the proof, we refer to the form~\eqref{eq:ADKBSC0} as coinciding with the one in~\cite[Eqn.~(25)]{AlbDalKon00}.

Assertion~\iref{i:t:ADKBSC:1} is claimed in~\cite[p.~284]{AlbDalKon00}.
\iref{i:t:ADKBSC:2} The closability of the form is claimed in~\cite[p.~289]{AlbDalKon00}; it is a consequence of the integration-by-parts formula~\cite[Eqn.~(44)]{AlbDalKon97} with~$\Lambda_k\equiv \zero$ for all~$k$, in the notation of~\cite{AlbDalKon97}. 
\iref{i:t:ADKBSC:3} Since~$\man$ is compact,~$\mbfM$ is compact as well by Tikhonov's Theorem, and the density of~$\Cyl{1}$ in~$\Cb(\mbfM)$ follows by Stone--Weierstra\ss\ Theorem. The density of~$\Cyl{1}$ in~$\dom{\dirE^\mbfs}$ is by definition, thus~$\Cyl{1}$ is a core for~$(\dirE^\mbfs,\dom{\dirE^\mbfs})$ and the form is regular; the standardness is immediate.
By~\cite[Cor.~I.5.1.4, Rmk.~I.5.1.5]{BouHir91}, it is sufficient to check strong locality on the core~$\Cyl{1}$; by finiteness of the set~$A$ in~\eqref{eq:Cyl} this is in turn a standard finite-dimensional fact.
\iref{i:t:ADKBSC:4} is claimed in~\cite[Thm.~4]{AlbDalKon97} and~\cite[Thm.~4.1]{AlbDalKon00}.

Provided we can identify the semigroup~$\TT^\mbfs_\bullet$ of~$(\dirE^\mbfs,\dom{\dirE^\mbfs})$ with~$\Heat^\mbfs_\bullet$ as in~\eqref{eq:Intro:ProdSemigroup},~\iref{i:t:ADKBSC:5} is the content of~\cite[Thm.~1.1]{BenSaC97} since
\begin{align}\label{eq:NonObvious}
\sum_n^\infty e^{-2\,\lambda_1\, t/ s_n}<\infty\comma \qquad \mbfs\eqdef\seq{s_n}_n^\infty\in\Tso\comma t>0\comma
\end{align}
where~$\lambda_1$ denotes the spectral gap of the Laplace--Beltrami operator of~$(\man,\mssg)$. In order to prove~\eqref{eq:NonObvious} it is sufficient to show that~$\nlimsup e^{-2\lambda_1 t/(n\,s_n)}<1$, by the root test. In turn, this is equivalent to~$\nliminf n\, s_n<\infty$. In fact, since~$\mbfs\in\Tso$, there exists~$\nlim n\,s_n=0$ by the Abel--Olivier--Pringsheim criterion.
The rest of the proof is devoted to the identification of the semigroup~$\TT^\mbfs_\bullet$ with~$\Heat^\mbfs_\bullet$. Since $\heat^\mbfs_t(\mbfx,\diff\mbfy)\ll\n\volm$, by~\cite[Lem.~6]{Kak48} we have in fact
\begin{align*}
\heat^\mbfs_t(\mbfx,\diff\mbfy)= \tonde{\prod_\ell^\infty \mssh_{t/s_\ell}(x_\ell,y_\ell)}\diff\n\volm(\mbfy)\comma \qquad \mbfx,\mbfy\in \mbfM \fstop
\end{align*}
In particular, the product of the densities converges. For~$u\in \Cyl{2}$ one has
\begin{align*}
\lim_{t\downarrow 0} \tfrac{1}{t} \ttonde{\Heat^\mbfs_t u-u}(\mbfx)=&\lim_{t\downarrow 0} \frac{1}{t}\tonde{\int_\mbfM \tonde{\prod_\ell^\infty \mssh_{t/s_\ell}(x_\ell,y_\ell)} F(\mbfy_A)\diff\n\volm(\mbfy) - F(\mbfx_A)}
\\
=&\lim_{t\downarrow 0} \frac{1}{t}\tonde{\int_{\man^{\times \abs{A}}} \tonde{\prod_{\ell\in A} \mssh_{t/s_\ell}(x_\ell,y_\ell)} F(\mbfy)\diff\n\mssm^{\abs{A}}(\mbfy)-F(\mbfx)} \comma
\end{align*}
also cf.~\cite[\S4.1]{Ber76}.
A standard finite-dimensional computation now shows that the generator~$(\Gen^\mbfs,\dom{\Gen^\mbfs})$ of~$\Heat^\mbfs_\bullet$ satisfies~$\Gen^\mbfs=\boldDelta^\mbfs$ on~$\Cyl{2}$. This concludes the proof of~\iref{i:t:ADKBSC:5} by essential self-adjointness,~\iref{i:t:ADKBSC:4}, of~$\boldDelta^\mbfs$ on~$\Cyl{2}$.~\iref{i:t:ADKBSC:6} is a direct consequence of~\iref{i:t:ADKBSC:5}.
\end{proof}
\end{thm}

\begin{rem} If~$\man=\mbbS^1$, Theorem~\ref{t:ADKBSC}\iref{i:t:ADKBSC:4} and~\iref{i:t:ADKBSC:5} are~\cite[\S4.1, Thm.s~4.3,~4.6]{Ber76}.
\end{rem}

\subsection{Spaces of measures} Let~$\msP\eqdef\msP(\man)$, resp.~$\Mp_1\eqdef \Mp_1(\man)$, denote the space of Borel probability, resp. subprobability, measures on~$\man$. On~$\msP$ we consider
\begin{enumerate*}[label=$\boldsymbol{(}$\bfseries{\itshape\alph*}$\boldsymbol{)}$]
\item the \emph{narrow} (or \emph{weak}) \emph{topology}~$\T_\mrmn$, induced by duality with~$\mcC_b(\man)$;

\item the \emph{weak atomic topology} $\T_\mrma$~\cite[\S2]{EthKur94}, see Proposition~\ref{p:EthierKurtz} below.
\end{enumerate*}

For $i=1,2$ let~$\mu_i\in\msP$. We denote by~$\Cpl(\mu_1,\mu_2)$ the set of couplings of~$\mu_1,\mu_2$, i.e.\ the set of probability measures~$\pi$ on~$\man^{\times 2}$ such that~$\pr^i_\pfwd \pi=\mu_i$, where~$\pr^i\colon \man^{\times 2}\rar \man$ is the projection on the $i^\text{th}$ component of the product. For~$p\in [1,\infty)$, the $L^p$-Kantorovich--Rubinshtein distance~$W_p$ on~$\msP$ is defined as
\begin{align}\label{eq:WD}
W_p(\mu_1,\mu_2)\eqdef \inf_{\pi\in\Cpl(\mu_1,\mu_2)} \tonde{\int_{\man^{\times 2}}\mssd_\mssg(x,y)^p \diff\pi(x,y)}^{1/p} \fstop
\end{align}

Since~$\man$ is compact, the narrow topology coincides with both the \emph{vague topology}, induced by duality with~$\mcC_c(\man)$, and with the topology induced by~$W_p$ for every~$p\in[1,\infty)$; see e.g.,~\cite[Cor.~6.13]{Vil09}.
For the reader's convenience, we collect here the main properties of the weak atomic topology which we shall dwell upon in the following, taken, almost \emph{verbatim}, from~\cite{EthKur94}.

\begin{rem}[{\cite[p.~5]{EthKur94}}]\label{r:SAlgP} The Borel $\sigma$-algebra~$\Bo{\mrmn}(\msP)\eqdef \Bo{\T_\mrmn}(\msP)$ generated by the narrow topology coincides with the Borel $\sigma$-algebra~$\Bo{\T_\mrma}(\msP)$ generated by the weak atomic topology.
\end{rem}

\begin{prop}[Ethier--Kurtz~{\cite{EthKur94}}]\label{p:EthierKurtz} The following hold:
\begin{enumerate*}[label=$\boldsymbol{(}$\bfseries{\itshape\roman*}$\boldsymbol{)}$]
\item\label{i:p:EthierKurtz0} $\T_\mrma$ is strictly finer than~$\T_\mrmn$;
\item\label{i:p:EthierKurtz1} $(\msP,\T_\mrma)$ is a Polish space;
\item\label{i:p:EthierKurtz2} suppose~$\T_\mrmn$-$\nlim\mu_n=\mu_\infty$. For~$N\in \overline{\N}_0$ let~$\seq{s_{N,i}\delta_{x_{N,i}}}_i^{m_N}$ be the set of atoms of~$\mu_N$, ordered in such a way that $s_{N,i-1}\geq s_{N,i}$ for all~$i\leq m_N\in \overline{\N}_0$. Then,
\end{enumerate*}
\begin{align*}
\text{$\T_\mrma$-$\nlim\mu_n=\mu_\infty$ if and only if~$s_{n,i}\longrar s_{\infty,i}$ for all~$i\leq m_\infty$} \semicolon
\end{align*} 
\begin{enumerate*}[resume*]
\item\label{i:p:EthierKurtz3} suppose~$\T_\mrmn$-$\nlim\mu_n=\mu_\infty$ and that~$\mu_\infty$ is purely atomic. Then,
\end{enumerate*}
\begin{align*}
\text{$\T_\mrma$-$\nlim\mu_n=\mu_\infty$ if and only if~$\nlim \sum_i\abs{s_{n,i}-s_{\infty,i}}=0$} \fstop
\end{align*}
Conversely, if~$\T_\mrma$-$\nlim\mu_n=\mu_\infty$ and~$s_{\infty,k}>s_{\infty,k-1}$ for some~$k\leq m_\infty$, then the set of locations~$\set{x_{n,1},\dotsc,x_{n,k}}$ converges to~$\set{x_{\infty,1},\dotsc, x_{\infty, k}}$. In particular, if~$s_{\infty,i}>s_{\infty,i+1}$ for all~$i$, then~$x_{n,i}\longrar x_{\infty,i}$ for all~$i$;
\begin{enumerate*}[resume*]
\item\label{i:p:EthierKurtz4} $(\msP,\T_\mrma)$ is \emph{not} compact, even if~$\man$ is.
\end{enumerate*}
\begin{proof}
For a metric metrizing~$(\msP,\T_\mrma)$ see~\cite[Eqn.~(2.2)]{EthKur94}. For separability and completeness see~\cite[Lem.~2.3]{EthKur94}. The inclusion in~\iref{i:p:EthierKurtz0} follows by comparison of~\cite[Eqn.~(2.2)]{EthKur94} with the Prokhorov metric; it is strict by~\cite[Example~2.7]{EthKur94}. 
For~\iref{i:p:EthierKurtz2}--\iref{i:p:EthierKurtz3} see~\cite[Lem.~2.5(b)]{EthKur94}.
Assume~$\man$ compact. We sketch a proof of~\iref{i:p:EthierKurtz4}. By e.g.~\cite[Rmk.~6.19]{Vil09}, $(\msP,\T_\mrmn)$ is compact. Argue by contradiction that~$(\msP,\T_\mrma)$ is compact. A continuous injection from a compact Hausdorff space is a homeomorphism onto its image,~\cite[Thm.~3.1.13]{Eng89}. Applying this to~$\id\colon (\msP,\T_\mrma)\rar (\msP,\T_\mrmn)$ contradicts~\iref{i:p:EthierKurtz0}.
\end{proof}
\end{prop}

\paragraph{The Dirichlet--Ferguson measure}

Everywhere in the following let~$\beta\in(0,\infty)$ be defined by~$\mssm=\beta\n\mssm$.
Let~$I\eqdef [0,1]$, resp.~$I^\circ\eqdef (0,1)$, be the closed, resp. open, unit interval and set~$\hman\eqdef \man\times I$ and~$\hman^\circ\eqdef \man\times I^\circ$, always endowed with the product topology, Borel $\sigma$-algebra~$\mcB(\hman)$ and with the measure~$\widehat\mssm_\beta\eqdef  \n\mssm\otimes \Beta_\beta$.
The next result may be regarded as a consequence of the Mecke identity~\eqref{eq:Mecke}.

\begin{cor}[cf.~{\cite[Prop.~1]{Fer73}}]\label{c:MeckeDF}
Let~$\eta$ be a~$\DF_\mssm$-distributed $\msP$-valued random field. Then, for every semi-bounded measurable~$f\colon \man\rar \R$,
\begin{align}\label{eq:diffuse}
\mbfE[\eta f]=\n\mssm f \fstop
\end{align}
In particular, for every measurable~$A$ it holds with $\DF_\mssm$-probability~$1$ that~$\eta A>0$ iff~$\mssm A>0$.
\end{cor}

\begin{defs} Denote by~$\Ppa$ the set of \emph{purely atomic} probability measures on~$\man$; by~$\Ppaiso$ the set of measures~$\eta\in\Ppa$ with infinitely many atoms and such that~$\eta_{x_1}\neq \eta_{x_2}$ whenever~$\eta_{x_1}>0$ and~$x_1\neq x_2$, cf.~Remark~\ref{r:Notation}; by~$\Pfs$ the set of \emph{fully supported} Borel probability measures on~$\man$. Finally set~$\Ppaisofs\eqdef\Pfs\cap \Ppaiso$.
\end{defs}

\begin{prop}\label{p:SuppDF} The following assertions hold:
\begin{enumerate*}[label=$\boldsymbol{(}$\bfseries{\itshape\roman*}$\boldsymbol{)}$]
\item\label{i:p:SuppDF1} $\Ppaisofs,\Ppaiso\in\Bo{\mrmn}(\msP)$ and $\Ppaisofs$ has full $\DF_\mssm$-measure;
\item\label{i:p:SuppDF2} let~$\boldPhi$ be defined as in~\eqref{eq:Intro:Phi}; then, its \emph{non-relabeled} restriction
\end{enumerate*}
\begin{align*}
\boldPhi\colon \ttonde{\andihM,\T_\mrmp,\hvolm_\beta}\longrar \tonde{\Ppaiso,\T_\mrma,\DF_{\beta\n\mssm}}
\end{align*}
is a homeomorphism and an isomorphism of measure spaces;
\begin{enumerate*}[resume*]
\item\label{i:p:SuppDF3} the set
\end{enumerate*}
\begin{align}\label{eq:p:SuppDF0}
\msN\eqdef \set{(\eta,x)\in \msP\times \man : \eta_x>0}\times I\subset \msP\times \hman
\end{align}
is $\Bo{\mrmn}(\msP)\otimes \Bo{}(\hman)$-measurable and $\DF_\mssm\otimes\widehat\mssm_\beta$-negligible.

\begin{proof}
Assertion~\iref{i:p:SuppDF1} is known, see e.g.~\cite{EthKur94} or~\cite[\S4]{Fer73}. 
In order to prove~\iref{i:p:SuppDF2} note that each of the considered topologies is metrizable, including~$\T_\mrma$, by Prop.~\ref{p:EthierKurtz}\iref{i:p:EthierKurtz1}, thus it suffices to show the continuity of~$\boldPhi$, resp.~$\boldPhi^{-1}$, along sequences.
To this end, for~$N\in \overline{\N}$ let
\begin{align*}
\mbfx_N\eqdef \seq{x_{N,i}}_i^\infty\in \andi\mbfM \comma \qquad \mbfs_N\eqdef& \seq{s_{N,i}}_i^\infty\in \Tso\comma \qquad \mu_N\eqdef \boldPhi(\mbfs_N,\mbfx_N)\in\Ppa\fstop
\end{align*}
The latter association is unique, since~$\boldPhi$ is bijective. 

Assume first~$\T_\mrmp$-$\nlim (\mbfx_n,\mbfs_n)=(\mbfx_\infty,\mbfs_\infty)$. In particular~$\ell^1$-$\nlim\mbfs_n=\mbfs_\infty$. For every~$f\in\mcC_b(\man)$,
\begin{align*}
\abs{\int_\man f \diff\mu_n -\int_\man f \diff\mu_\infty}=&\abs{\sum_i^\infty \ttonde{s_{n,i}f(x_{n,i})-s_{\infty,i}f(x_{\infty,i})}}
\\
\leq&\sum_i^\infty s_{n,i} \abs{f(x_{n,i})-f(x_{\infty,i})}+\sum_i^\infty \abs{s_{n,i}-s_{\infty,i}}\abs{f(x_{\infty,i})}
\\
\leq&\sum_i^\infty s_{n,i} \abs{f(x_{n,i})-f(x_{\infty,i})}+\norm{f}_{\infty} \norm{\mbfs_n-\mbfs_\infty}_{\ell^1} \fstop
\end{align*}

The first term vanishes as~$n\rar\infty$ by Dominated Convergence Theorem with varying dominating functions~$\norm{f}_{\infty}\mbfs_n\in \ell^1(\N)$; the second term vanishes by assumption. By arbitrariness of~$f$, $\T_\mrmn$-$\nlim \mu_n=\mu_\infty$, whence~$\T_\mrma$-$\nlim \mu_n=\mu_\infty$ by Proposition~\ref{p:EthierKurtz}\iref{i:p:EthierKurtz3}. This shows the continuity of~$\boldPhi$.
The continuity of~$\boldPhi^{-1}$ is precisely the converse statement in Proposition~\ref{p:EthierKurtz}\iref{i:p:EthierKurtz3}.
It follows that~$\boldPhi$ is bi-measurable. The measure isomorphism property is also known, \emph{Sethuraman stick-breaking representation}, cf.~\cite{Set94}.
\iref{i:p:SuppDF3} The set~$\msN$ is measurable by Lemma~\ref{l:JointMeas}. Moreover, its sections $\msN_x\eqdef\set{\eta\in\msP: (\eta, x, r)\in \msN}$ are $\DF_\mssm$-negligible for every~$x\in \man$ by Corollary~\ref{c:MeckeDF}, hence~$\msN$ is $\DF_\mssm\otimes\widehat\mssm_\beta$-negligible itself.
\end{proof}
\end{prop}

\begin{rem}
It is not possible to extend the homeomorphism in Proposition~\ref{p:SuppDF}\iref{i:p:SuppDF2}, in the sense that the spaces~$(\Ppa,\T_\mrma)$ and~$\To\times\andi{\mbfM}$ are \emph{not} homeomorphic. Clearly, the same holds for~$(\Ppa,\T_\mrma)$ and~$\Tso\times\mbfM$, for which~$\boldPhi$ is not even bijective.
\end{rem}

\begin{rem}\label{r:LocCompact}
It is noted in~\cite[Prop.~3.1]{BanRad97} that~$\Ppa$ is an~$F_{\sigma\delta}$-set in~$(\msP,\T_\mrmn)$, and thus so is~$\Ppafs$. The same holds in~$\T_\mrma$. Neither subspace is locally compact in~$\T_\mrmn$, nor in~$\T_\mrma$.
\end{rem}

\section{The Dirichlet form}
In this section, we assume~$d\geq 1$ whenever not stated otherwise.

\subsection{Cylinder functions}
We introduce some spaces of suitably differentiable functions.

\begin{defs}[Cylinder functions]\label{d:CylFunc}
Let~$k,\ell,m,n\in \N_0$ and~$\mcC^m_b(\R^k)$ be the space of real-valued bounded $m$-differentiable functions on~$\R^k$ with bounded derivatives of any order up to~$m$.

For~$\hat f\in \mcB_b(\hman)$ and~$\eta\in\msP$ set
\begin{align}\label{eq:Testing}
\hat f^\trid(\eta)\eqdef \sum_{x\in \eta} \eta_x\, \hat f(x,\eta_x)=\int_\man \diff\eta(x)\, \hat f(x,\eta_x) \fstop
\end{align}

For~$\hat f_i\in \mcB_b(\hman)$ for~$i\leq k$, set~$\hat\mbff\eqdef\tseq{\hat f_1,\dotsc, \hat f_k}$ and~$\hat\mbff^\trid(\eta)\eqdef \tseq{\hat f_1^\trid(\eta),\dotsc, \hat f_k^\trid(\eta)}$.

For~$\eps\in I$ set further~$\hman_\eps\eqdef \man\times (\eps,1]$. We always regard~$\Cc^m(\hman_\eps)$ as the subspace of~$\mcC^m(\hman)$ obtained by extension by~$0$. Consistently with this identification, we put~$\mcC^0(\hman_1)\eqdef \set{0}$ by convention.
Note that~$\Cz^0(\hman_0)\subsetneq \mcC^0(\hman)$.
We define the following families of \emph{cylinder functions}
\begin{equation}\label{eq:TestF}
\begin{aligned}
\hTF{m}{}\eqdef&\set{\begin{matrix} u\colon \msP\rar \R : u= F\circ\hat\mbff^\trid, \\
F\in \mcC^m_b(\R^k), \hat \mbff\in\mcC^m(\hman)^{\otimes k} \end{matrix}}\comma
\\
\hTF{m}{-}\eqdef&\set{\begin{matrix} u\colon \msP\rar \R : u= F\circ\hat\mbff^\trid, F\in \mcC^m_b(\R^k),\\ \hat f_i=\car_\man\otimes \varrho_i, \varrho_i\in \mcC^m(I)\comma i\leq k \end{matrix}}\comma
\\
\hTF{m}{\eps}\eqdef&\set{\begin{matrix} u\colon \msP\rar \R : u= F\circ\hat\mbff^\trid, \\
F\in \mcC^m_b(\R^k), \hat \mbff\in\Cc^m(\hman_\eps)^{\otimes k} \end{matrix}}\comma
\\
\TF{m}\eqdef&\set{\begin{matrix} u\colon \msP\rar \R : u= F\circ\hat\mbff^\trid, F\in \mcC^m_b(\R^k),\\ \hat f_i=f_i\otimes \car_I, f_i\in \mcC^m(\man) \comma i\leq k \end{matrix}}\comma
\\
\hTF{m}{-,\eps}\eqdef&\hTF{m}{-}\cap \hTF{m}{\eps} \fstop
\end{aligned}
\end{equation}

For~$u\in\hTF{m}{}$ we define the \emph{vanishing threshold~$\eps_u$ of~$u$} by
\begin{align*}
\eps_u\eqdef \sup\tset{\eps\in I : u\in \hTF{m}{\eps}}\fstop
\end{align*}

Note that~$\eps_u>0$ for all~$u\in\hTF{0}{0}$.
Finally, for~$\eps\in I$, define the family of $\sigma$-algebras
\begin{align}\label{eq:SAlgebras}
\Bo{\eps}(\msP)\eqdef \sigma_0\ttonde{\hTF{\infty}{\eps}} \fstop
\end{align}
\end{defs}

\begin{rem}[Representation of cylinder functions]\label{r:ReprCyl}
The representation of~$u$ by~$F$ and~$\hat\mbff$ is never unique.
Indeed, assume~$u\in \hTF{m}{}$ may be written as~$u=F\circ\hat\mbff^\trid$ for appropriate~$F$ and~$\hat\mbff$. By compactness of~$\msP$ and~$\hman$ and by our definition of the test functions~$\hat f^\trid$, if~$G\colon \R^k\rar\R$ satisfies~$G\equiv F$ identically on~$\prod_i^k \im \hat f_i$, then~$u=G\circ \hat\mbff^\trid$.
As a consequence:
\begin{enumerate*}[label=$\boldsymbol{(}$\bfseries{\itshape\alph*}$\boldsymbol{)}$]
\item the families in~\eqref{eq:TestF} remain unchanged if we replace~$\mcC_b^m(\R^k)$ with~$\Cc^m(\R^k)$ or~$\mcC^m(\R^k)$;

\item in particular, if~$\hat f\in \mcC^m(\hman)$, then the induced test function~$\hat f^\trid$ belongs to~$\hTF{m}{}$ (and analogously for the other families of functions in~\eqref{eq:TestF});

\item if additionally~$F$ is constant in the direction~$\mbfe_j$ on~$\im \hat f_j$ for some~$j\leq k$, then $u=G\circ \hat\mbfg^\trid$, where $\hat\mbfg\eqdef \tseq{\hat f_1,\dotsc, \hat f_{j-1}, \hat f_{j+1},\dotsc, \hat f_k}$ and~$G\in\mcC^m_b(\R^{k-1})$ is such that, for some, hence any,~$\bar t\in\im\hat f_j$,
\end{enumerate*}
\begin{align*}
G(\mbfs)=F(s_1,\dotsc, s_j, \bar t, s_{j+1},\dotsc,s_{k-1})\comma \qquad \mbfs\eqdef (s_1,\dotsc, s_{k-1})\in \prod_{i \, : \, i\neq j}^k \im \hat f_i \semicolon
\end{align*}
\begin{enumerate*}[resume*]
\item if~$u\in \hTF{m}{}$ there exists a minimal~$k$ such that~$u=F\circ \hat\mbff^\trid$ for~$F\in\mcC^m(\R^k)$ and appropriate~$\hat \mbff^\trid$. If this is the case, we say that~$u$ is written in minimal form.
\end{enumerate*}
In the following, we shall always assume every cylinder function to be written in minimal form.
\end{rem}

\begin{rem}[Measurability, continuity of cylinder functions]\label{r:MeasContCyl}
\begin{enumerate*}[label=$\boldsymbol{(}$\bfseries{\itshape\alph*}$\boldsymbol{)}$]
\item Every~$u\in\hTF{0}{}$ is measurable as a consequence of Lemma~\ref{l:Meas};
\item every non-constant~$u\in\hTF{0}{0}$ is $\T_\mrmn$-dis\-con\-tin\-u\-ous at~$\DF_\mssm$-a.e.~$\mu$, even for~$u\in \hTF{\infty}{\eps}$, for every~$\eps\in I$; 
\item\label{i:r:MeasContCyl:3} every function~$u\in\hTF{0}{0}$ is $\T_\mrma$-continuous as a consequence of~\cite[Rmk.~2.6]{EthKur94};
\item\label{i:r:MeasContCyl:4} every function in~$\TF{0}$ is $\T_\mrmn$-continuous by definition of~$\T_\mrmn$;
\item\label{i:r:MeasContCyl:5} $\TF{m}$, $\hTF{m}{-}$, $\hTF{m}{\eps}$ and $\hTF{m}{}$ are algebras with respect to the pointwise multiplication of real-valued functions on~$\msP$ and are closed with respect to pre-composition with $m$-differentiable functions, i.e.\ if e.g.~$u\in\TF{m}$ and~$\psi\in\mcC^m(\R;\R)$, then~$\psi\circ u\in \TF{m}$;

\item\label{i:r:MeasContCyl:6} the sequences~$m\mapsto \TF{m}$, $\hTF{m}{-}$, $\hTF{m}{\eps}$, $\hTF{m}{}$ are decreasing;

\item\label{i:r:MeasContCyl:7} $\TF{m}$, $\hTF{m}{-}$, $\hTF{m}{\eps}\subsetneq \hTF{m}{}$ (\emph{strict} inclusion) for every~$m\in \overline{\N}_0$ and~$\eps\in I$;

\item\label{i:r:MeasContCyl:8} $\eps\mapsto \hTF{m}{\eps}$ is decreasing and left-continuous, in the sense that~$\hTF{m}{\eps}=\cup_{\delta>\eps} \hTF{m}{\delta}$ for every~$m\in \overline{\N}_0$ and~$\eps\in I$;

\item\label{i:r:MeasContCyl:9} $\hTF{0}{0}\cap \TF{0}=\R$ (constant functions);

\item\label{i:r:MeasContCyl:10} $\Bo{\eps}(\msP)$ does not separate points in~$\msP$ for any~$\eps\geq 0$.
\end{enumerate*}
\end{rem}

\begin{lem}\label{l:SAlg} It holds that 
\begin{enumerate*}[label=$\boldsymbol{(}$\bfseries{\itshape\roman*}$\boldsymbol{)}$]
\item\label{i:l:SAlg1} $\Bo{1}(\msP)=\set{\emp,\msP}$;
\item\label{i:l:SAlg2} $\Bo{\mrmn}(\Ppa)=\Bo{0}(\msP)_{\Ppa}$; 
and
\item\label{i:l:SAlg3} $\cl_{L^2(\msP,\Bo{\mrmn},\DF_\mssm)}(\hTF{\infty}{\eps})=\cl_{L^2(\msP,\Bo{\eps},\DF_\mssm)}(\hTF{\infty}{\eps})=L^2(\msP,\Bo{\eps},\DF_\mssm)$.
\end{enumerate*}

\begin{proof}
\iref{i:l:SAlg1} is immediate, since~$\hTF{\infty}{1}=\R$. As for~\iref{i:l:SAlg2}, note that the family of pointwise limits of sequences in~$\hTF{\infty}{0}$ contains the algebra
\begin{align*}
\hTFB{0}_0\eqdef&\set{\begin{matrix} u\colon \msP\rar \R : u\eqdef F\circ\hat\mbff^\trid, F\in \mcC^0_b(\R^k),\\ \hat f_i=f_i\otimes \car_{(0,1]}, f_i\in\mcC^0_b(\man) \comma i\leq k \end{matrix}}\comma
\end{align*}
and, for~$u=F\circ \ttonde{\seq{f_i\otimes \car_{(0,1]}}_{i\leq k}}^\trid\in \hTFB{0}_0$ let~$\tilde u\eqdef F\circ\mbff^\trid\in \TF{0}$. Clearly~$u(\eta)=\tilde{u}(\eta)$ for every~$\eta\in \Ppa$, hence
\begin{align*}
\Bo{0}(\msP)_{\Ppa}\eqdef \sigma_0(\hTF{\infty}{0})_{\Ppa}\supset \sigma_0(\hTFB{0}_0)_{\Ppa}=\sigma_0(\TF{0})_{\Ppa} \fstop
\end{align*}
Since~$\T_\mrmn$ on~$\msP$ is generated by the linear functionals~$f^\trid\in \TF{0}$ varying~$f\in \mcC_b(\man)$, one has~$\sigma_0(\TF{0})_{\Ppa}=\Bo{\mrmn}(\Ppa)$.
Thus,~$\Bo{0}(\msP)_{\Ppa}\supset \Bo{\mrmn}(\Ppa)$.
On the other hand,
\begin{align*}
\Bo{\mrmn}(\msP)=\Bo{\mrma}(\msP)\supset\Bo{0}(\msP)
\end{align*}
by Rmk.s~\ref{r:SAlgP} and~\ref{r:MeasContCyl}\iref{i:r:MeasContCyl:3} respectively, hence~$\Bo{\mrmn}(\Ppa)\supset \Bo{0}(\msP)_{\Ppa}$ and the conclusion follows.
The first equality in~\iref{i:l:SAlg3} is immediate, since~$\hTF{\infty}{\eps}\subset L^2(\msP,\Bo{\eps},\DF_\mssm)$ by boundedness of functions in~$\hTF{\infty}{\eps}$ and finiteness of~$\DF_\mssm$.
The second equality is not entirely straightforward, cf.~Rmk.~\ref{r:MeasContCyl}\iref{i:r:MeasContCyl:10}.
It is however a consequence of Proposition~\ref{p:Dynkin} which we postpone to the Appendix.
\end{proof}
\end{lem}

As a consequence of the proof of Lemma~\ref{l:SAlg}\iref{i:l:SAlg3}, we have that~$\Bo{\eps}(\msP)_{\Ppa}=\sigma_0(\hTF{m}{\eps})_{\Ppa}$ and that we may replace~$\hTF{\infty}{\eps}$ with~$\hTF{m}{\eps}$ in the statement of Lemma~\ref{l:SAlg}\iref{i:l:SAlg3} for any~$m\in \N_0$.

\subsubsection{Directional derivatives of cylinder functions}
In the following, if~$\upphi\colon \man\rar \man$ is measurable, set~$\Phi\eqdef \upphi_\pfwd\colon\msP\rar \msP$. In particular,
\begin{align*}
\Psi^{w,t} \eqdef \fl^{w,t}_\pfwd \fstop
\end{align*}

\begin{defs} For~$w\in \Vect^\infty$ and~$u\in \hTF{1}{}$ we define the derivative of~$u$ in the direction of $w$
\begin{align}\label{eq:DirDer}
\grad_w u(\eta)\eqdef \diff_t\restr_{t=0} u\ttonde{\Fl^{w,t} \eta}
\end{align}
whenever it exists.
\end{defs}

\begin{rem}[Geometries of~$\msP$]\label{r:Shifts}
The \emph{shift} $S_{tf}$ in~\eqref{eq:ShiftHanda} (considered in~\cite{Han02, Sch02, Sha11}; see~\cite[p.~546]{Han02} for the terminology) is \emph{not} the `exponential map' of~$\msP_2$, i.e.\ in the sense of the $L^2$-Wasserstein geometry of~$\msP$. Rather, it is associated to $\msP_1$, where the convex combination~$\mu\mapsto \mu^x_t$ \emph{is} a geodesic curve.
In fact, the map~$\Fl^{w,t}$ ---~one might suggestively write~$(e^{tw})_\pfwd$~--- is also not the exponential map~$\bexp$ of~$\msP_2$, studied in~\cite{Gig11}. However,~$\Fl^{w,t}$ is tangent to~$\bexp t\,\cdot \boldgamma$ for some appropriately chosen `tangent plan'~$\boldgamma \in \msP_2(T\man)$ depending on~$w$, as shown in~\cite[Lem.~4.3]{LzDS19b}.
\end{rem}

\begin{lem}[Directional derivative]\label{l:DirDer}
Let~$u\in \hTF{1}{}$,~$w\in \Vect^\infty$ and~$\eta\in\msP$. Then, there exists
\begin{align}\label{eq:l:DirDer0}
\grad_w u(\eta) = \sum_i^k (\partial_i F)(\hat \mbff^\trid\eta)\cdot \int_\man \tgscal{\nabla \hat f_i(x,\eta_x)}{w(x)} \diff \eta(x) \fstop
\end{align}

Furthermore $\grad_w\colon \mfZ\longrar \mfZ$ for $\mfZ=\hTF{\infty}{},\hTF{\infty}{\eps},\TF{\infty}$ while~$\grad_w\colon \hTF{1}{-}\longrar \set{0}$, and
\begin{align*}
\tnorm{\grad_w u}_{L^2(\DF_\mssm)}\leq \sqrt k \norm{\nabla F}_\infty \max_i \tnorm{\nabla \hat f_i}_{\Vect^0} \norm{w}_{\Vect_{\n\mssm}} \fstop
\end{align*}

\begin{proof} We show that the curve~$t\mapsto u\ttonde{\Fl^{w,t} \eta}$ is differentiable for every~$t$. Indeed
\begin{align*}
\diff_t u\tonde{\Fl^{w,t}\eta}
=&\sum_i^k (\partial_i F)(\hat \mbff^\trid\eta) \cdot \diff_t \int_\man \hat f_i\ttonde{y, \tonde{\Fl^{w,t}\eta}_y} \diff\!\tonde{\Fl^{w,t}\eta}(y)\\
%
%
=&\sum_i^k (\partial_i F)(\hat \mbff^\trid\eta) \cdot \diff_t \int_\man \hat f_i\tonde{\fl^{w,t}(y), \eta_y}  \diff \eta(y) \fstop
\intertext{\indent Since~$\hat f\in \mcC^1(\hman)$, differentiation under integral sign yields}
\diff_t u\tonde{\Fl^{w,t}\eta}=&\sum_i^k (\partial_i F)(\hat \mbff^\trid\eta) \cdot \int_\man \gscal{\nabla \hat f_i\tonde{\fl^{w,t}(y), \eta_y} }{ \dot \fl^{w,t}(y)} \diff \eta(y)
%
\end{align*}

Computing at~$t=0$ yields~\eqref{eq:l:DirDer0}.
For the second claim, note that, by smoothness of~$w$, $\tgscal{\nabla\hat f_i(\emparg,\emparg)}{w(\emparg)}\in \mcC^\infty(\hman)$ as soon as~$\hat f_i$ is. One can estimate
\begin{align*}
\tnorm{\grad_wu}_{L^2(\DF_\mssm)}^2 \leq& \norm{\nabla F}_\infty^2\cdot k \max_i \int_\msP  \abs{\int_\man \gscal{\nabla \hat f_i(x, \eta_x)}{w(x)} \diff\eta (x)}^2 \diff\DF_\mssm(\eta)\\
\leq&k\norm{\nabla F}_\infty^2 \max_i \tnorm{\nabla \hat f_i}_{\Vect^0}^2 \int_\msP\norm{w}^2_{\Vect_\eta} \diff\DF_\mssm(\eta) \\
=&k \norm{\nabla F}_\infty^2 \max_i \tnorm{\nabla \hat f_i}_{\Vect^0}^2 \norm{w}_{\Vect_{\n\mssm}}^2
\end{align*}
by~\eqref{eq:diffuse}, which concludes the proof.
\end{proof}
\end{lem}

\subsubsection{Integration-by-parts formula}
We discuss integration by parts for cylinder functions.

\begin{lem}[Local derivative and Laplacian]\label{l:LocDerLap}
Let $u\eqdef F\circ \hat\mbff^\trid\in\hTF{0}{0}$. 
Then, the function $U\colon(\eta,z,r)\mapsto u(\eta^z_r)$ is~$\Bo{\mrmn}(\msP)\otimes \Bo{}(\hman)$-measurable.
Furthermore,
\begin{enumerate}[label=$\boldsymbol{(}$\bfseries{\itshape\roman*}$\boldsymbol{)}$]
\item~if~$u$ is in~$\hTF{1}{0}$, then $\forallae{\DF_\mssm\otimes \widehat\mssm_\beta} (\eta,x,r)$ the map~$z\mapsto U(\eta,z,r)$ is differentiable in a neighborhood of~$z=x$ and
\begin{align}\label{eq:l:IbP1}
\nabla^z_w\restr_{z=x} U(\eta,z,r)=&\, r \sum_i^k (\partial_i F)\ttonde{\hat \mbff^\trid (\eta^x_r)}\cdot \nabla_w \hat f_i(x,r)\fstop
\end{align}

\item~if~$u$ is in~$\hTF{2}{0}$, then $\forallae{\DF_\mssm\otimes \widehat\mssm_\beta} (\eta,x,r)$ the map~$z\mapsto U(\eta,z,r)$ is twice differentiable in a neighborhood of~$z=x$ and
\begin{equation}\label{eq:l:IbP2}
\begin{aligned}
\Delta^z\restr_{z=x} U(\eta,z,r)=&\, r^2\sum_{i,j}^k (\partial_{ij}^2 F)\ttonde{\hat\mbff^\trid(\eta^x_r)} \cdot \tgscal{\nabla\hat f_i(x,r)}{\nabla\hat f_j(x,r)}\\
&\quad +r\sum_i^k (\partial_i F)\ttonde{\hat\mbff^\trid(\eta^x_r)} \cdot\Delta \hat f_i (x,r) \fstop
\end{aligned}
\end{equation}
\end{enumerate}

Furthermore, the right-hand sides of~\eqref{eq:l:IbP1} and~\eqref{eq:l:IbP2} are~$\Bo{\mrmn}(\msP)\otimes \Bo{}(\hman)$-mea\-sur\-able.

\begin{proof}
By continuity of the Dirac embedding~$x\mapsto \delta_x$ and Lem.~\ref{l:Meas} the function~$(\eta,x,r)\mapsto \eta^x_r$ is continuous. As a consequence,~$U$ is measurable by Remark~\ref{r:MeasContCyl}.
Let~$\msN$ be as in~\eqref{eq:p:SuppDF0}.
For~$(\eta,z,r)\not\in\msN$ and every~$\hat f\in\Cz^0(\hman)$ one has
\begin{align}\label{eq:l:LocDerLap1}
\hat f^\trid(\eta^z_r)=\sum_{y\in \eta^z_r} (\eta^z_r)_y \hat f\ttonde{y,(\eta^z_r)_y}= r\hat f(z,r)+ \sum_{y\in \eta} (1-r) \eta_y \hat f\ttonde{y,(1-r)\eta_y} \fstop
\end{align}

Thus,
\begin{align*}
\nabla^z_w\restr_{z=x} F(\hat\mbff^\trid\eta_r^z)=& \sum_i^k (\partial_i F)(\hat \mbff^\trid\eta^x_r) \cdot \nabla^z_w\restr_{z=x} \hat f_i^\trid (\eta_r^z)\\
=&\sum_i^k (\partial_i F)(\hat\mbff^\trid\eta^x_r) \cdot \nabla^z_w\restr_{z=x}\sum_{y\in \eta^z_r} (\eta^z_r)_y\, \hat f_i\ttonde{y,(\eta^z_r)_y}\\
=&\sum_i^k (\partial_i F)(\hat \mbff^\trid\eta^x_r)\cdot 
\\
&\cdot\tonde{\nabla^z_w\restr_{z=x} (\eta^z_r)_z\, \hat f_i\ttonde{z,(\eta^z_r)_z}+\sum_{y\in \eta} \nabla^z_w\restr_{z=x} (\eta^z_r)_y\, \hat f_i\ttonde{y,(\eta^z_r)_y} }\comma
\end{align*}
where the gradient may be exchanged with the sum, since the latter is always over a finite number of points by the choice of $\hat f_i$. In light of~\eqref{eq:l:LocDerLap1},
\begin{align*}
\nabla^z_w\restr_{z=x}& F\ttonde{\hat\mbff^\trid(\eta_r^z)}=
\\
=&\sum_i^k (\partial_i F)\cdot
\\
&\qquad \cdot\ttonde{\hat \mbff^\trid(\eta^x_r)}\tonde{\nabla^z_w\restr_{z=x} r \hat f_i(z,r)+\sum_{y\in \eta} \nabla^z_w\restr_{z=x} (1-r)\eta_y\, \hat f_i\ttonde{y,(1-r)\eta_y} }
\\
=&\sum_i^k (\partial_i F)\ttonde{\hat \mbff^\trid(\eta^x_r)}\cdot r \nabla_w \hat f_i(x,r)\fstop
\end{align*}

By~\eqref{eq:l:IbP1} and arbitrariness of~$w$ one has
\begin{align*}
\nabla^z\restr_{z=x} u(\eta^z_r)=& \ r\sum_i^k (\partial_i F)\ttonde{\hat\mbff^\trid(\eta^x_r)} \cdot \nabla \hat f_i(x,r)\comma
\intertext{hence, if~$u$ is sufficiently regular,}
\Delta^z\restr_{z=x} u(\eta^z_r)=&\ (\div^{\mssm,z}\circ\nabla^z)\restr_{z=x}u(\eta^z_r)\\
=&\ r\sum_i^k \tgscal{\nabla^z\restr_{z=x} (\partial_i F)\ttonde{\hat \mbff^\trid(\eta^z_r)}}{\nabla\hat f_i(x,r)}
\\
&\ +r\sum_i^k (\partial_i F)\ttonde{\hat\mbff^\trid(\eta^x_r)} \cdot \Delta^z\restr_{z=x}\hat f_i(z,r)\\
%
%
=&r^2 \sum_{i,j}^k (\partial_{ji}^2 F)\ttonde{\hat\mbff^\trid(\eta^x_r)} \cdot \tgscal{\nabla \hat f_j(x,r)}{\nabla\hat f_i(x,r)}
\\
&\ +r\sum_i^k (\partial_i F)\ttonde{\hat\mbff^\trid(\eta^x_r)} \cdot \Delta\hat f_i (x,r) \fstop
\end{align*}

This shows~\eqref{eq:l:IbP1} and~\eqref{eq:l:IbP2} outside the $\DF_\mssm\otimes \widehat\mssm_\beta$-negligible set~$\msN$.
\end{proof}
\end{lem}

\begin{thm}[Integration by parts]\label{l:IbP}
Let $w\in\Vect^\infty$ and $u\eqdef F\circ \hat\mbff^\trid$, $v\eqdef G\circ\hat \mbfg^\trid$ be cylinder functions in $\hTF{1}{0}$. Set~$\eps\eqdef \eps_u\wedge\eps_v>0$. Then, the following integration-by-parts formula holds:
\begin{equation}\label{eq:l:IbP00}
\begin{aligned}
\int_\msP \grad_w u\cdot v \diff\DF_\mssm =&-\int_\msP u\cdot \grad_w v \diff\DF_\mssm - \int_\msP u\cdot v \cdot \mbfB_\eps[w] \diff\DF_\mssm\comma
\end{aligned}
\end{equation}
where
\begin{align}\label{eq:l:IbP0}
\mbfB_\eps[w](\eta)\eqdef \sum_{x: \eta_x> \eps} \div^\mssm_x w \fstop
\end{align}

\begin{proof}
We can compute
\begin{align*}
\int_\msP& \grad_w u \cdot v \diff\DF_\mssm=
\\
=& \int_\msP v(\eta) \cdot \sum_i^k (\partial_i F)\ttonde{\hat\mbff^\trid(\eta)} \int_\man \tgscal{\nabla \hat f_i(x, \eta_x)}{w(x)} \diff \eta (x) \diff \DF_\mssm(\eta)
\\
=&\int_\msP \int_\man v(\eta) \cdot \sum_i^k (\partial_i F)\ttonde{\hat\mbff^\trid(\eta)} \cdot \tgscal{\nabla \hat f_i (x, \eta_x)}{w(x)} \diff\eta(x) \diff \DF_\mssm(\eta)\comma\\
\intertext{whence, by the Mecke identity~\eqref{eq:Mecke} and by~\eqref{eq:l:IbP1},}
=&\int_\msP \int_{\hman} v(\eta_r^x) \cdot \sum_i^k (\partial_i F)\ttonde{\hat\mbff^\trid(\eta_r^x)} \, \tgscal{\nabla \hat f_i(x, r)}{w(x)} \diff\widehat\mssm_\beta(x,r) \diff \DF_\mssm(\eta)\comma\\
=&\int_\msP \int_0^1 \car_{(\eps,1]}(r) \int_\man v(\eta^x_r) \cdot \tfrac{1}{r}\nabla_w^z\restr_{z=x} u(\eta^z_r) \diff \mssm(x) \frac{\diff \Beta_\beta(r)}{\beta} \diff \DF_\mssm(\eta) \fstop
\end{align*}

Since~$\bd \man=\emp$, by integration by parts on~$\man$
\begin{align}
\nonumber
\int_\msP& \grad_w u \cdot v \diff \DF_\mssm=
\\
\nonumber
=& -\int_\msP \int_{\hman} \nabla_w^z\restr_{z=x} v(\eta^z_r) \cdot u(\eta^x_r) \diff\widehat\mssm_\beta(x,r) \diff \DF_\mssm(\eta)
\\
\nonumber
&\quad-\int_\msP \int_{\hman} \frac{\car_{(\eps,1]}(r)}{r} (uv)(\eta^x_r)\cdot \div^\mssm_x w \, \diff\widehat\mssm_\beta(x,r) \diff \DF_\mssm(\eta)
\\
\label{eq:l:IbP:1}
=&- \int_\msP \int_{\hman} u(\eta_r^x) \cdot \sum_j^h (\partial_j G)(\hat\mbfg^\trid\eta_r^x) \,\gscal{\nabla \hat g_j(x, r)}{w(x)} \diff\widehat\mssm_\beta(x,r) \diff \DF_\mssm(\eta)
\\
\label{eq:l:IbP:2}
&\quad-\int_\msP \int_{\hman} \frac{\car_{(\eps,1]}(r)}{r} (uv)(\eta^x_r)\cdot \div^\mssm_x w \, \diff\widehat\mssm_\beta(x,r)\diff \DF_\mssm(\eta)\fstop
\end{align}

Applying the Mecke identity~\eqref{eq:Mecke} to~\eqref{eq:l:IbP:1} and~\eqref{eq:l:IbP:2},
\begin{align*}
\int_\msP& \int_{\hman} u(\eta_r^x) \cdot \sum_j^h (\partial_j G)\ttonde{\hat\mbfg^\trid(\eta_r^x)} \,\gscal{\nabla \hat g_j(x, r)}{w(x)} \diff\widehat\mssm_\beta(x,r) \diff \DF_\mssm(\eta)=
\\
&=\int_\msP u(\eta)\cdot \grad_w v(\eta) \diff\DF_\mssm(\eta) \comma
\\
\int_\msP&\int_{\hman} \frac{\car_{(\eps,1]}(r)}{r}  (uv)(\eta^x_r)\cdot \div^\mssm_x w\, \diff\widehat\mssm_\beta(x,r) \diff \DF_\mssm(\eta)=
\\
&=\int_\msP \int_\man \frac{\car_{(\eps,1]}(\eta_x)}{\eta_x}  (uv)(\eta)\cdot \div^\mssm_x w\, \diff\eta(x)\diff \DF_\mssm(\eta)
\\
&=\int_\msP (uv)(\eta)\cdot \sum_{ \eta_x> \eps} \div^\mssm_x w\, \diff \DF_\mssm(\eta)\fstop \qedhere
\end{align*}
\end{proof}
\end{thm}

\subsection{Gradient and Dirichlet form on~\texorpdfstring{$\msP$}{P}}
At each point~$\mu$ in~$\msP$, the directional derivative~$\grad_w u$ of any~$u\in \dom{\grad_w}$ defines a linear form~$w\mapsto \grad_w u(\mu)$ on~$\Vect^\infty$. Let~$\norm{\emparg}_\mu$ be a pre-Hilbert norm on~$\Vect^\infty$ such that this linear form is continuous, and let~$T_\mu\msP$ denote the completion of~$\Vect^\infty$ with respect to the said norm. By Riesz Representation Theorem there exists a unique element~$\grad u(\mu)$ in~$T_\mu\msP$ such that~$ \grad_w u(\mu)=\tscalar{\grad u(\mu)}{w}_\mu$, where~$\scalar{\emparg}{\emparg}_\mu$ denotes the scalar product of the Hilbert space~$T_\mu\msP$.
Different choices of~$\norm{\emparg}_\mu$, hence of~$T_\mu\msP$, yield different gradient maps~$\grad u$, namely, as suggested by Lemma~\ref{l:DirDer}, the closures of the operator
\begin{align}\label{eq:DefGrad}
\grad u(\mu)(x)\eqdef& \sum_i^k (\partial_i F)\ttonde{\hat\mbff^\trid(\mu)} \cdot \nabla \hat f_i(x,\mu_x)\comma\qquad u\eqdef F\circ \hat\mbff^\trid\in \hTF{1}{} \fstop
\end{align}

\begin{rem}[Measurability of gradients]\label{r:MeasContGrad}
The function $x\mapsto \grad u(\mu)(x)$ is measurable for every~$u\in\hTF{1}{0}$ and~$\mu$ in~$\msP$ by measurability of~$x\mapsto \mu_x$, whereas it is generally discontinuous at~$\DF_\mssm$-a.e.~$\mu$, even for~$u\in \hTF{\infty}{}$. 
\end{rem}

\subsubsection{The Dirichlet form~\texorpdfstring{$\mcE$}{hE}}

Throughout this section fix~$T_\mu\msP=T^\Der_\mu\msP_2\eqdef \Vect_\mu$. The $\norm{\emparg}_{\Vect_\mu}$-continuity of~$w\mapsto \grad u(\mu)$, granting that our choice is admissible, readily follows from~\eqref{eq:l:DirDer0}. We refer the reader to~\cite[\S5.1]{LzDS19b} for the geometrical meaning of this choice.

For~$u\eqdef F\circ\hat\mbff^\trid\in\hTF{0}{}$ where~$\hat\mbff\eqdef\tseq{\hat f_1,\dotsc, \hat f_k}$, denote by~$\tilde{u}$ the extension of~$u$ to~$\Mb(\man)$ defined by extending~$\hat f_i^\trid\colon \msP\rar \R$ to~$\hat f_i^\trid\colon \Mb(\man)\rar\R$ in the obvious way for all~$i$'s. For the purpose of clarity, in the \emph{statement} of the following theorem we distinguish~$u$ from~$\tilde u$. Everywhere else, with slight abuse of notation, we will denote both~$u$ and~$\tilde u$ simply by~$u$.
\begin{thm}\label{t:DForm}
Assume~$d\geq 1$. For~$u,v$ in~$\hTF{2}{0}$ set
\begin{align}
\label{eq:FormMain}
\mcE(u,v)\eqdef&\ \tfrac{1}{2} \int_\msP \tscalar{\grad u(\eta)}{\grad v(\eta)}_{\Vect_\eta} \diff \DF_{\beta\n\mssm}(\eta) \comma \\
\label{eq:hL}
\mbfL u(\eta)\eqdef&\ \tfrac{1}{2} \int_\man \frac{\Delta^z\restr_{z=x} \tilde u(\eta+\eta_x\delta_z-\eta_x\delta_x)}{(\eta_x)^2} \diff\eta(x) \comma && \eta\in\Ppa\comma
\\
\label{eq:hGamma}
\boldGamma(u,v)(\eta)\eqdef&\ \tfrac{1}{2} \tscalar{\grad u(\eta)}{\grad v(\eta)}_{\Vect_\eta} \comma && \eta\in\msP \fstop
\end{align}

Then,~$(\mbfL,\hTF{2}{0})$ is a symmetric operator on~$L^2(\msP,\DF_\mssm)$ satisfying
\begin{align*}
\mcE(u,v)=\tscalar{u}{-\mbfL v}_{L^2(\DF_\mssm)} \comma \qquad u,v\in\hTF{2}{0} \fstop
\end{align*}
The bilinear form~$(\mcE,\hTF{2}{0})$ is a closable symmetric form on~$L^2(\msP,\DF_\mssm)$.
Its closure~$(\mcE,\dom{\mcE})$ is a strongly local recurrent (in particular: conservative) Dirichlet form with generator the Friedrichs extension~$(\mbfL_{\mathsc{f}},\dom{\mbfL_{\mathsc{f}}})$ of~$(\mbfL,\hTF{2}{0})$.
Moreover,~$(\mcE,\dom{\mcE})$ has carr\'e du champ operator~$(\boldGamma,\dom{\boldGamma})$ where~$\dom{\boldGamma}\eqdef\dom{\mcE}\cap L^\infty(\msP,\DF_\mssm)$, that is, for all~$u,v,z\in\dom{\boldGamma}$,
\begin{align}\label{eq:CdC}
2\int_\msP z\, \boldGamma(u,v) \diff\DF_\mssm= \mcE(u,vz)+\mcE(uz,v)-\mcE(uv,z) \fstop
\end{align}
\begin{proof} 
By definition of~$\grad u$,
\begin{align}
\label{eq:DoubleGrad}
2\, &\mcE(u,v)=
\\
\nonumber
=&\int_\msP \int_\man \sum_{i,j}^{k,h} \gscal{(\partial_i F)\ttonde{\hat\mbff^\trid(\eta)} \nabla \hat f_i(x,\eta_x)}{(\partial_j G)\ttonde{\hat\mbfg^\trid(\eta)} \nabla \hat g_j(x,\eta_x)} \diff\eta(x)\diff\DF_\mssm(\eta)\comma
\intertext{whence, by the Mecke identity~\eqref{eq:Mecke} and by integration by parts on~$\man$,}
\nonumber
=& \int_\msP\int_I \int_{\man} \gscal{\tfrac{1}{r}\nabla^z\restr_{z=x} u(\eta^z_r) }{\tfrac{1}{r} \nabla^z\restr_{z=x} v(\eta^z_r) } \diff\n\mssm(x)\diff\Beta_\beta(r)\diff\DF_\mssm(\eta)
\\
\nonumber
=&- \int_\msP\int_I \int_{\man} u(\eta^x_r) \cdot \Delta^z\restr_{z=x} v(\eta^z_r) \diff\mssm(x) \frac{\diff \Beta_\beta(r)}{\beta r^2}\diff\DF_\mssm(\eta)
\\
%
\nonumber
=&- \int_\msP \int_{\hman} u(\eta^x_r) \cdot \frac{\Delta^z\restr_{z=x}v(\eta^x_r+r\delta_z-r\delta_x)}{r^2} \diff\widehat\mssm_\beta(x,r) \diff\DF_\mssm(\eta) \comma
\intertext{thus, again by the Mecke identity~\eqref{eq:Mecke},}
\nonumber
=&-2\int_\msP u(\eta) \cdot \mbfL v(\eta) \diff\DF_\mssm(\eta) \fstop
\end{align}

Let~$\mcH\eqdef \cl_{L^2(\msP,\DF_\mssm)} \hTF{2}{0}$ as a Hilbert subspace of~$L^2(\msP,\DF_\mssm)$. Clearly, $\hTF{2}{}\subset \mcH$, hence in particular~$\TF{2}\subset \mcH$ and the family~$\TF{2}$ is a unital nowhere-vanishing algebra of $\T_\mrmn$-continuous functions (cf.~Rem.~\ref{r:MeasContCyl}) separating points in~$\msP$, thus it is uniformly dense in~$\mcC(\msP,\T_\mrmn)$ by compactness of~$\msP$ and the Stone--Weierstra{\ss} Theorem.
Since~$(\msP,\mcB_\mrmn(\msP),\DF_\mssm)$ is a compact Polish probability space,~$\cl_{L^2(\msP,\DF_\mssm)}\mcC(\msP)=L^2(\msP,\DF_\mssm)$, thus finally~$\mcH=L^2(\msP,\DF_\mssm)$.
It is straightforward that~\eqref{eq:hL} defines a linear operator~$\mbfL\colon\hTF{2}{0}\rar L^2(\msP,\DF_\mssm)$. The symmetry (and coercivity) of the bilinear form~$(\mcE,\hTF{2}{0})$ is obvious. Its closability on~$L^2(\msP,\DF_\mssm)$ and the existence of the Friedrichs extension~$(\mbfL,\dom {\mbfL})$ follow from~\cite[Thm.~X.23]{ReeSim75}. The Markov and strong local properties are also straightforward since~$\hTF{2}{0}$ is closed w.r.t. post-composition with smooth real-valued functions.

By the Leibniz rule for~$\grad$,~\eqref{eq:CdC} holds for all~$u,v,z\in \hTF{1}{0}$. Arbitrary~$u,v,z\in\dom{\boldGamma}$ may be respectively approximated both in~$\mcE^{1/2}_1$ and~$\DF_\mssm$-a.e. by uniformly bounded sequences~$u_n,v_n,z_n\in\hTF{1}{0}$. Thus~$\nlim u_nv_n=uv$,~$\nlim u_nz_n= uz$ and~$\nlim v_nz_n= vz$ in~$\mcE^{1/2}_1$ and
\begin{align*}
\nlim& \int_\msP \tabs{z\, \boldGamma(u,v)-z_n \boldGamma(u_n,v_n)}\diff\DF_\mssm
\\
&\quad \leq \nlim \int_\msP \abs{z-z_n} \boldGamma(u,v)\diff\DF_\mssm+\nlim \int_\msP \abs{z_n}\tabs{\boldGamma(u,v)-\boldGamma(u_n,v_n)}\diff\DF_\mssm=0 \comma
\end{align*}
whence~\eqref{eq:CdC} carries over from~$\hTF{1}{0}$ to~$\dom{\mcE}\cap L^\infty(\msP,\DF_\mssm)$.
%
Since~$\car\in \dom\mcE$ and $\mcE(\car)=0$, the form is recurrent, thus conservative~\cite[Thm.~1.6.3, Lem.~1.6.5]{FukOshTak11}.
\end{proof}
\end{thm}

\begin{rem} Note that~$(\mcE,\dom{\mcE})$,~$(\mbfL_{\mathsc{f}}, \dom{\mbfL_{\mathsc{f}}})$ and~$(\boldGamma,\dom{\boldGamma})$ all depend on~$\beta$. 
We assume~$\beta>0$ to be fixed and drop it from the notation.
\end{rem}

\begin{rem}
For~$u=F\circ\hat\mbff^\trid\in\hTF{2}{0}$ with vanishing threshold~$\eps_u$, set, consistently with~\eqref{eq:l:IbP0},
\begin{align*}
\mbfB[\nabla \hat f_i](\eta)=\sum_{x\in\eta} \Delta\hat f_i(x,\eta_x)=\sum_{x: \eta_x>\eps_u} \Delta \hat f_i(x,\eta_x)\comma\qquad i\leq k\comma
\end{align*}
and
\begin{align}\label{eq:L1L2}
\mbfL_1 u\eqdef&  \tfrac{1}{2} \sum_{i,p}^k (\partial^2_{ip} F\circ \hat\mbff^\trid)\cdot \boldGamma(\hat f_i^\trid, \hat f_p^\trid)\comma
&
\mbfL_2 u\eqdef& \tfrac{1}{2} \sum_i^k(\partial_i F\circ \hat\mbff^\trid) \cdot \mbfB[\nabla \hat f_i]\fstop
\end{align}

Also note that
\begin{align}\label{eq:GammaGammaStar}
\boldGamma(\hat f_i^\trid,\hat f_p^\trid)=\Gamma(\hat f_i,\hat f_p)^\trid \comma\qquad i,p\leq k\comma
\end{align}
where~$\Gamma(\hat f_i,\hat f_p)(x)\eqdef \tfrac{1}{2}\tgscal{\nabla_x \hat f_i}{\nabla_x \hat f_p}$ is but the carr\'e du champ operator of the Laplace--Beltrami operator on~$(\man,\mssg)$.
Then, consistently with~\eqref{eq:l:IbP0} and~\eqref{eq:l:IbP00},
\begin{align}\label{eq:LapBis}
\mbfL=& \ \mbfL_1+\mbfL_2 \comma
\end{align}
which makes apparent that~$\mbfL u$ is defined everywhere on~$\msP$ and is identically vanishing on the subspace of diffuse measures. for every~$u\in\hTF{2}{0}$.
\end{rem}

\subsubsection{The \texorpdfstring{$\T_\mrmn$}{tn}-regularity of~\texorpdfstring{$\mcE$}{E}}

In view of Remark~\ref{r:MeasContCyl},~$\dom{\mcE}$ might appear unsuitable for the form to be regular, since we defined the latter on a core of \emph{non-continuous} functions. The goal of this section is to show that, in fact,~$\dom{\mcE}$ contains sufficiently many continuous functions.

\begin{defs}[Sobolev functions of mixed regularity and Sobolev cylinder functions]
Denote by~$I_\beta$ the measure space~$(I, \Beta_\beta)$ and consider the space
\begin{align*}
W_{\widehat\mssm_\beta}\eqdef L^2\ttonde{I_\beta;W^{1,2}_{\n\mssm}(\man)}\cong L^2(I_\beta)\, \widehat\otimes\, W^{1,2}_{\n\mssm}(\man)\comma
\end{align*}
where~$\widehat\otimes$ denotes the tensor product of \emph{Hilbert} spaces. It coincides with the completion of~$\mcC^\infty(\hman)$ with respect to the norm defined by
\begin{align*}
\tnorm{\hat f}_{W_{\widehat\mssm_\beta}}^2\eqdef \int_{\hman} \ttonde{\tabs{\nabla\hat f}^2_\mssg+\tabs{\hat f}^2} \diff\widehat\mssm_\beta \fstop
\end{align*}

To fix notation, let~$\hat f\colon (x,s)\mapsto \hat f_s(x)\eqdef\hat f(x,s)\in W_{\widehat\mssm_\beta}$.
We denote further by~$D$ the distributional gradient on~$\man$ and by~$D^{1,0}\eqdef (D\otimes \id_{I^\circ})$ the distributional differential operator given, locally on a chart of~$\hman^\circ$, by differentiation along coordinate directions on~$\man$.

By~\cite[Prop.~3.105]{AmbFusPal00}, for $\hat f\in W_{\widehat\mssm_\beta}$ and for a.e.~$s\in I^\circ$ one has that, locally on~$\hman^\circ$, differentiation in the $\man$-directions commutes with restriction in the $I^\circ$-direction, i.e.~$D\hat f_s=(D^{1,0}\hat f)_s$. Thus, for any such~$\hat f$, the notation~$D\hat f$ is unambiguous.
For~$\hat f\in W_{\widehat\mssm_\beta}$, we denote by~$[\hat f,D\hat f]$ any of its Borel representatives. We write~$[\hat f, D\hat f]_1$ when referring only to the representative of~$\hat f$, and~$[\hat f, D\hat f]_2$ when referring only to the representative of~$D\hat f$.
Finally set
\begin{align}\label{eq:SobolevCylFunc}
\hTFW{2,2}{b}\eqdef&\set{\begin{matrix} u\colon \msP\rar \R : u= F\circ\tseq{[\hat f_1,D\hat f_1]_1^\trid,\dotsc, [\hat f_k,D\hat f_k]_1^\trid}, \\ 
F\in \mcC^2_b(\R^k)\comma \hat f_i\in W_{\widehat\mssm_\beta} \comma [\hat f_i,D\hat f_i]_1\in\mcB_b(\hman) \comma i\leq k \end{matrix}}\fstop
\end{align}
\end{defs}

\begin{rem} The specification of representatives for both~$\hat f$ and~$D\hat f$ in the definition of~$\hTFW{2,2}{b}$ is instrumental to the statement of Lemma~\ref{l:Domains} below. It is then the content of the Lemma that such a specification is in fact immaterial.
\end{rem}

\begin{lem}\label{l:Domains}
Let~$(\mcE,\dom{\mcE})$ be defined as in Theorem~\ref{t:DForm}. Then,
\begin{enumerate*}[label=$\boldsymbol{(}$\bfseries{\itshape\roman*}$\boldsymbol{)}$]
\item\label{i:l:Domains1} $\hTF{2}{}\subset \dom{\mcE}$;
\end{enumerate*}
and
\begin{enumerate*}[resume*]
\item\label{i:l:Domains2} $\hTFW{2,2}{b}\subset \dom{\mcE}$ and~$u\in \dom{\mcE}$ of the form~\eqref{eq:SobolevCylFunc} does not depend on the choice of the representatives for~$\hat f_i$.
\end{enumerate*}
Furthermore, for any such~$u$, for~$\DF_\mssm$-a.e.~$\eta\in\msP$,
\begin{align*}
\boldGamma(u)(\eta)=&\sum_{i,p} (\partial_i F)\ttonde{[\hat\mbff,D\hat\mbff]_1^\trid(\eta)} \cdot (\partial_p F)\ttonde{[\hat\mbff,D\hat\mbff]_1^\trid(\eta)} \cdot
\\
&\cdot\int_\man \gscal{[\hat f_i,D\hat f_i]_2}{[\hat f_p,D\hat f_p]_2}(x,\eta_x)\diff\eta(x)\comma
\end{align*}
with usual meaning of the notation~$\hat\mbff$, does not depend on the choice of representatives for~$\hat f_i$.
\begin{proof}
\iref{i:l:Domains1} Let~$u=F\circ\hat\mbff\in \hTF{2}{}$ and~$\varrho_n\in \mcC^\infty(I)$ be such that~$\varrho_n\uparrow_n \car_I$ pointwise and~$\supp\varrho_n\subset [1/n,1]$. For~$i\leq k$ set~$\hat f_{n,i}\eqdef \varrho_n\cdot\hat f_i$ and note that~$\hat f_{n,i}\in \mcC^2(\hman_{1/n})$, hence~$u_n\eqdef F\circ \hat\mbff_n\in \hTF{2}{1/n}$ for every~$n\in \N$. It is straightforward that
\begin{align*}
\max_{i\leq k} \nlim \tnorm{\hat f_{n,i}-\hat f_i}_{W_{\widehat\mssm_\beta}}=0\comma
\end{align*}
hence there exists~$C_u>0$ such that~$\max_{i\leq k} \sup_n \tnorm{\hat f_{n,i}}_{W_{\widehat\mssm_\beta}}\leq C_u$. Thus,
\begin{align}
\nonumber
2\,\mcE&(u_n-u_m)=
\\
=&\ \int_\msP\int_\man \abs{\sum_i^k \tonde{(\partial_i F)\ttonde{\hat\mbff_n^\trid(\eta)}\cdot\nabla\hat f_{n,i}-(\partial_i F)\ttonde{\hat\mbff_m^\trid(\eta)\cdot \nabla\hat f_{m,i}}}}^2_\mssg\!\!\!(x,\eta_x) \,\diff\eta(x) \diff\DF_\mssm(\eta)
\\
\nonumber
\leq& \ 2\int_\msP\int_\man \abs{\sum_i^k (\partial_i F)\ttonde{\hat\mbff_n^\trid(\eta)}\cdot\ttonde{\nabla\hat f_{n,i}-\nabla\hat f_{m,i}}(x,\eta_x)}_\mssg^2 \diff\eta(x) \diff\DF_\mssm(\eta)
\\
\nonumber
&\ +2\int_\msP\int_\man \abs{\sum_i^k \tonde{(\partial_i F)\ttonde{\hat\mbff_n^\trid(\eta)}-(\partial_i F)\ttonde{\hat\mbff_m^\trid(\eta)}}\cdot \nabla\hat f_{m,i}(x,\eta_x)}_\mssg^2 \diff\eta(x) \diff\DF_\mssm(\eta)
\\
\nonumber
\leq&\ 2^k k \cdot \Lip(F)^2 \cdot \sum_i^k \int_\msP\int_\man \abs{\nabla\hat f_{n,i}-\nabla\hat f_{m,i}}_\mssg^2(x,\eta_x) \diff\eta(x) \diff\DF_\mssm(\eta)
\\
\nonumber
&\ +2^k k \cdot \max_i\Lip(\partial_i F)^2 \cdot C_u^2 \cdot \sum_i^k \int_\msP\abs{\int_\man \ttonde{\hat f_{n,i}-\hat f_{m,i}}(x,\eta_x)}^2 \diff\eta(x)\diff\DF_\mssm(\eta)
\\
\label{eq:Estimate}
\leq&\ K_u \cdot \sum_i^k \mcE_1\ttonde{\hat f_{n,i}^\trid-\hat f_{m,i}^\trid} \comma
\end{align}
for some appropriate constant~$K_u$, independent of~$i,n,m$. Now, by Jensen inequality
\begin{align*}
2\,\mcE_1\ttonde{\hat f_{n,i}^\trid-\hat f_{m,i}^\trid}=&\int_\msP \Gamma(\hat f_{n,i}-\hat f_{m,i})^\trid+\tabs{(\hat f_{n,i}-\hat f_{m,i})^\trid}^2\diff\DF_\mssm(\eta)
\\
\leq& \int_\msP \Gamma\ttonde{\hat f_{n,i}-\hat f_{m,i}}^\trid+\ttonde{\tabs{\hat f_{n,i}-\hat f_{m,i}}^2}^\trid \diff\DF_\mssm(\eta)\comma
\end{align*}
thus, by the Mecke identity~\eqref{eq:Mecke}
\begin{align}
\nonumber
2\,\mcE_1(\hat f_{n,i}^\trid-\hat f_{m,i}^\trid)\leq&\ \int_\msP \int_{\hman} \tonde{\Gamma(\hat f_{n,i}-\hat f_{m,i})+\tabs{\hat f_{n,i}-\hat f_{m,i}}^2} \diff\widehat\mssm_\beta \diff\DF_\mssm
\\
\nonumber
=&\ \int_{\hman}\tonde{\Gamma(\hat f_{n,i}-\hat f_{m,i})+\tabs{\hat f_{n,i}-\hat f_{m,i}}^2} \diff \widehat\mssm_\beta
\\
\label{eq:FundSeq}
=&\ \tnorm{\hat f_{n,i}-\hat f_{m,i}}_{W_{\widehat\mssm_\beta}}^2 \fstop
\end{align}

This shows that the sequence~$\tseq{\mcE(u_n)}_n$ is fundamental, therefore bounded. Analogously, one can show that the sequence~$u_n$ converges to~$u$ strongly in~$L^2(\msP,\DF_\mssm)$. Thus,~$u\in \dom{\mcE}$ by~\cite[Lem.~I.2.12]{MaRoe92}. Since~$\tseq{\mcE_1(u_n)}_n$ is fundamental, letting~$n\rar\infty$ in~\eqref{eq:FundSeq} and combining it with~\eqref{eq:Estimate} yields
\begin{align}\label{eq:l:Domains1}
\nlim\mcE_1(u_n- u)=0\comma \qquad u\in\hTF{2}{}\comma u_n\in \hTF{2}{1/n}\fstop
\end{align}


Note that the condition~$[\hat f_i,D\hat f_i]_1\in\mcB_b(\hman)$ grants that~$[\hat f_i,D\hat f_i]_1^\trid(\eta)$ is well-defined by~\eqref{eq:Testing}. The measurability of~$u$ follows as in Remark~\ref{r:MeasContCyl}. As a consequence,~\iref{i:l:Domains2} may be proven similarly to~\iref{i:l:Domains1} by $\norm{\emparg}_{W_{\widehat\mssm_\beta}}$-density of~$\mcC^2(\hman)$ in~$W_{\widehat\mssm_\beta}$.
\end{proof}
\end{lem}

\begin{rem} In view of Lemma~\ref{l:Domains}\iref{i:l:Domains2}, everywhere in the following we write~$\hat f^\trid$ in place of~$[\hat f,D\hat f]_1^\trid\in \hTFW{2,2}{b}\setminus \hTF{2}{}$. Analogously, for any such~$\hat f$ we will write~$\boldGamma(\hat f^\trid)=\ttonde{\tabs{D\hat f}_\mssg^2}^\trid$ omitting any explicit indication of the representative of~$\hat f\in W_{\widehat\mssm_\beta}$. We note that, with a little more effort, one could show that~$\hat f^\trid$ is a well-defined element of~$L^2(\msP,\DF_\mssm)$ for any~$\hat f\in L^2(\hman,\widehat\mssm_\beta)$ and independent of the choice of representatives for~$\hat f$. In a similar way,  one can show that~$\hTFW{2,2}{}\subset \dom{\mcE}$, with obvious meaning of the notation~$\hTFW{2,2}{}$, as opposed to~$\hTFW{2,2}{b}$.
\end{rem}

\begin{cor}\label{c:Domains}
Let~$u\in\dom{\mcE}$. Then, there exists~$\seq{u_n}_n$ such that
\begin{enumerate*}[label=$\boldsymbol{(}$\bfseries{\itshape\roman*}$\boldsymbol{)}$]
\item\label{i:c:Domains1}$u_n\in\hTF{2}{1/n}$ for all~$n\in \N$;
\item\label{i:c:Domains2}$\mcE_1^{1/2}$-$\nlim u_n=u$;
\item\label{i:c:Domains3}$\DF_\mssm$-$\nlim u_n=u$;
\item\label{i:c:Domains4}$\DF_\mssm$-$\nlim \tnorm{\grad u_n-\grad u}_{\Vect_\emparg}=0$.
\end{enumerate*}

\begin{proof} \iref{i:c:Domains3} and~\iref{i:c:Domains4} are a standard consequence of~\iref{i:c:Domains2} up to passing to a suitable subsequence. Thus, it suffices to show~\iref{i:c:Domains1} and~\iref{i:c:Domains2}, which in turn follow by Lemma~\ref{l:Domains}\iref{i:l:Domains1} and an $\eps/3$-argument. 
\end{proof}
\end{cor}

\begin{cor}\label{c:Regularity}
Assume~$d\geq 2$. The Dirichlet form~$(\mcE,\dom \mcE)$ on~$L^2(\msP,\DF_\mssm)$ is a $\T_\mrmn$-regular strongly local recurrent (hence: conservative) Dirichlet form with standard core
\begin{align*}
\TFL{2}{1}\eqdef&\set{\begin{matrix} u\colon \msP\rar \R: u\eqdef F\circ\hat\mbff^\trid, F\in \mcC^2_b(\R^k),\\ \hat f_i=f_i\otimes \car_I\comma f_i\in \Lip(\man) \comma i\leq k \end{matrix}}\fstop
\end{align*}

\begin{proof}
The family~$\TF{\infty}\subset \hTFW{2,2}{b} \subset \dom{\mcE}$ is uniformly dense in~$\mcC(\msP)$ as in the proof of Thm.~\ref{t:DForm}. A proof that~$\TF{\infty}$ is also dense in~$\dom{\mcE}$ if~$d\geq 2$ is postponed to Lemma~\ref{l:TFDensity}.

One has~$\TFL{2}{1}\subset \hTFW{2,2}{b} \subset \dom{\mcE}$ and in fact~$\TFL{2}{1}\subset\dom{\mcE}$ similarly to the proof of Lemma~\ref{l:Domains}\iref{i:l:Domains2}. In particular, for any~$f\in\Lip(\man)$, the section~$D f$ is defined on~$A^\complement$, where the singular set~$A$ of~$f$ satisfies~$A\in\Bo{\mssg}$ and $\mssm A=0$ by the classical Rademacher Theorem. For~$u=F\circ\mbff^\trid\in \TFL{2}{1}$ with~$\mbff\eqdef\seq{f_1,\dotsc, f_k}$, let~$A_i$ denote the singular set of~$f_i$ and set~$A\eqdef \cup_{i\leq k} A_i$. Then,
\begin{equation*}
\begin{aligned}
\boldGamma(u)(\eta)=& \tfrac{1}{2}\int_\man \sum_{i,p}^k (\partial_i F)(\mbff^\trid\eta) \cdot  (\partial_p F)(\mbff^\trid\eta)\, \cdot
\\ 
&\qquad \cdot \gscal{D f_i}{D f_p}\diff \eta(x)
\end{aligned}
\qquad\forallae{\DF_\mssm}\eta\comma \qquad u\in\TFL{2}{1}\comma
\end{equation*}
is well-defined, since the set of measures~$\eta\in\Ppaisofs$ charging~$A$ is $\DF_\mssm$-negligible by Corollary~\ref{c:MeckeDF}. The fact that~$\TFL{2}{1}$ is a standard core is a straightforward consequence of the definition of~$\TFL{2}{1}$ and of the classical chain rule. 
\end{proof}
\end{cor}

\subsubsection{Partial quasi-invariance of~\texorpdfstring{$\DF_\mssm$}{DFm}}
The following result is heuristically clear from the analogous result~\cite[Thm.~13]{KonLytVer15} for the Gamma measure. However, it seems to me that it cannot be rigorously deduced from it. Thus, we provide an independent proof.

\begin{prop}\label{p:PQI}
The measure~$\DF_\mssm$ is partially quasi-invariant with respect to the action~$\mfG\acts\msP$ as in~\eqref{eq:Intro:GAct} on the filtration~$\mcB_\bullet(\Ppa)\eqdef \seq{\mcB_{1/n}(\Ppa)}_n$ as in~\eqref{eq:SAlgebras}, with Radon--Nikod\'ym derivatives~$\mbfR_{1/n}[\uppsi]$ as in~\eqref{eq:Intro:RND}.

\begin{proof} It suffices to establish~\iref{i:d:PQI3} in Definition~\ref{d:PQI}. Indeed~\iref{i:d:PQI1} was noted in Definition~\ref{d:CylFunc} and~\iref{i:d:PQI2} is straightforward with~$n'=n$. In order to check~\iref{i:d:PQI3} it suffices to restrict ourselves to functions~$u\in \hTF{\infty}{1/n}$, since they generate~$\mcB_{1/n}(\Ppa)$ by definition. Now, for~$\uppsi\in \Diff^\infty_+(\man)$,
\begin{align*}
\int_\msP u(\eta) \diff(\uppsi.)_\pfwd\DF_\mssm(\eta) = \int_{\Ppaiso} u(\uppsi_\pfwd\eta) \diff \ttonde{\boldPhi_\pfwd\hvolm_\beta}(\eta)
\end{align*}
by Proposition~\ref{p:SuppDF}\iref{i:p:SuppDF2}. Let~$\eta=\sum_i^\infty s_i\delta_{x_i}$. Since~$s_1>s_2>\dotsc$, one has~$s_{n+1}< 1/n$, hence, for every~$i\leq k$, every~$x\in \man$ and every~$n'>n$ it holds that~$\hat f_j(x,s_{n'})=0$ by definition of~$\hat f_i$.
Thus,
\begin{align*}
\int_\msP& u(\eta) \diff(\uppsi.)_\pfwd\DF_\mssm(\eta) =
\\
=&\int_{\andi{\mbfM}}\int_{\Tso} F\tonde{\sum_i^\infty s_i \hat f_1(\uppsi(x_i),s_i),\dotsc,\sum_i^\infty s_i \hat f_k(\uppsi(x_i),s_i)} \diff \PD_\beta(\mbfs)\diff \n\volm(\mbfx)
\\
=&\int_{\man^{\times n}} \int_{\Tso} F\tonde{\sum_i^n s_i \hat f_1(\uppsi(x_i),s_i),\dotsc,\sum_i^n s_i \hat f_k(\uppsi(x_i),s_i)} \diff\PD_{\beta}(\mbfs)\diff\n\mssm^n \ttonde{\mbfx^{(n)}}
\\
=&\int_{\man^{\times n}} \int_{\Tso} F\tonde{\sum_i^n s_i \hat f_1(x_i,s_i),\dotsc,\sum_i^n s_i \hat f_k(x_i,s_i)} \diff\PD_{\beta}(\mbfs)\diff(\uppsi_\pfwd\n\mssm)^{\otimes n} \ttonde{\mbfx^{(n)}}
\\
=&\int_{\man^{\times n}} \prod_i^n J^\mssm_{\uppsi}(x_i) \int_{\Tso} F\tonde{\sum_i^n s_i \hat f_1(x_i,s_i),\dotsc,\sum_i^n s_i \hat f_k(x_i,s_i)} \diff\PD_{\beta}(\mbfs) \diff\n\mssm^n \ttonde{\mbfx^{(n)}}
\\
=&\int_{\andi{\mbfM}} \int_{\Tso} \prod_{i : s_i> 1/n} \!\!\!\!\! J^\mssm_{\uppsi}(x_i) \,\, F\tonde{\sum_i^n s_i \hat f_1(x_i,s_i),\dotsc,\sum_i^n s_i \hat f_k(x_i,s_i)} \diff\PD_\beta(\mbfs)\diff\n\volm(\mbfx) 
\\
=&\int_{\msP} \mbfR_{1/n}[\uppsi](\eta)\cdot u(\eta) \diff\DF_{\mssm}(\eta) \fstop \qedhere
\end{align*}
\end{proof}
\end{prop}

Together with the definition of partial quasi-invariance, Proposition~\ref{p:PQI} suggests that some of the quantities we defined in terms of the $\sigma$-algebras~$\mcB_\eps(\msP)$ ought to be martingales with respect to the filtration~$\seq{\mcB_\eps(\msP)}_{\eps\in I}$. Indeed, this turns out to be the case. The following is a corollary of~Theorem~\ref{l:IbP}.

\begin{cor}\label{c:Martingale}
Let~$w\in\Vect^\infty$ and~$\mbfB_\eps[w]$ be defined as in~\eqref{eq:l:IbP0}. Then, with the same notation of Proposition~\ref{p:PQI},
\begin{enumerate*}[label=$\boldsymbol{(}$\bfseries{\itshape\roman*}$\boldsymbol{)}$]
\item\label{i:c:Martingale1} the stochastic process~$\mbfB_\bullet[w]\eqdef\tseq{\mbfB_{\eps}[w]}_{\eps\in I}$ is a centered square-integrable martingale on $\ttonde{\Ppa, \mcB_\mrmn(\Ppa), \DF_\mssm}$ with respect to the filtration
\end{enumerate*}
\begin{align*}
\mcB_\bullet(\Ppa)\eqdef\seq{\mcB_{\eps}(\Ppa)}_{\eps\in I}\comma
\end{align*}
cf.~\eqref{eq:SAlgebras};
\begin{enumerate*}[resume,label=$\boldsymbol{(}$\bfseries{\itshape\roman*}$\boldsymbol{)}$]
\item\label{i:c:Martingale1.5} it holds that
\end{enumerate*}
\begin{align*}
\mbfB_\bullet[w]=\diff_t\restr_{t=0} \mbfR_\bullet[\psi^{w,t}] \semicolon
\end{align*}
\begin{enumerate*}[resume*]
\item\label{i:c:Martingale2} the quadratic form
\end{enumerate*}
\begin{align*}
\mcA_0^w\colon \hTF{1}{0}\times \hTF{1}{0} \ni (u,v)\longmapsto \mbfE_{\DF_\mssm}\tquadre{u \cdot v \cdot \mbfB_\eps[w]}\comma \qquad \eps\eqdef \eps_u\wedge \eps_v\comma
\end{align*}
is $\mcE^{1/2}_1$-bounded, and uniquely extends to an $\mcE^{1/2}_1$-bounded quadratic form~$\mcA^w$ on $\dom{\mcE}$;

\begin{proof}
\iref{i:c:Martingale1} Let~$\delta>\eps>0$ and~$\varrho_{\eps,\delta}\in \Cc^\infty([\eps,1];\R)$ be such that~$\varrho_{\eps,\delta}(r)= 1/r$ for every~$r\geq \delta$. Set~$\hat f_{w,\eps,\delta}\eqdef  \div^\mssm w\, \otimes \varrho_{\eps,\delta}\in \hTF{\infty}{\eps}$ and note that~$\mbfB_\eps[w]=\lim_{\delta\downarrow \eps}\hat f_{w,\eps,\delta}^\trid$ pointwise on~$\Ppa$.
Thus,~$\mbfB_\eps[w]$ is~$\mcB_\eps(\Ppa)$-measurable. This shows that the process~$\mbfB_\bullet[w]$ is adapted to~$\mcB_\bullet(\Ppa)$.
Moreover,
\begin{align}\label{eq:c:Mg0}
\tabs{\mbfB_\eps[w](\eta)}\leq \sum_{x: \eta_x> \eps} \norm{\div^\mssm w}_{\mcC^0}\leq \floor{\eps^{-1}}\norm{\div_\mssm w}_{\mcC^0} \fstop
\end{align}

Choosing~$v=\car$ in~\eqref{eq:l:IbP00} yields
\begin{align*}
\int_\msP \grad_w u(\eta) \diff\DF_\mssm(\eta)= \int_\msP u(\eta)\, \mbfB_{\eps}[w](\eta) \diff\DF_\mssm(\eta) \fstop
\end{align*}

Since~$u$ is~$\mcB_{\eps}(\Ppa)$-measurable, it is also~$\mcB_{\eps}(\Ppa)$-measurable for all~$\delta\leq \eps$, hence
\begin{align*}
\int_\msP u(\eta)\, \mbfB_{\delta}[w](\eta) \diff\DF_\mssm(\eta) = \int_\msP u(\eta)\, \mbfB_{\eps}[w](\eta) \diff\DF_\mssm(\eta)\comma \qquad \delta\leq \eps \fstop
\end{align*}

By arbitrariness of~$u\in\hTF{1}{\eps}$ and Lemma~\ref{l:SAlg}\iref{i:l:SAlg3},
\begin{align}\label{eq:c:Mg1}
\mbfE_{\DF_\mssm}\tquadre{\mbfB_{\delta}[w] \mid \mcB_{\eps}(\Ppa)}=\mbfB_{\eps}[w] \fstop
\end{align}

Then,~$\mbfB_\bullet[w]$ is a martingale by~\eqref{eq:c:Mg0} and \eqref{eq:c:Mg1}.
Finally, choosing~$u=v=\car$ in~\eqref{eq:l:IbP00} yields~$\mbfE_{\DF_\mssm}\tquadre{\mbfB_{\eps}[w]}=0$.
\iref{i:c:Martingale1.5} For every~$\eps>0$ one has
\begin{align*}
\diff_t\restr_{t=0}\mbfR_\eps[\fl^{w,t}](\eta)=&\diff_t\restr_{t=0} \exp\quadre{\int_\man \car_{(\eps,1]}(\eta_x) \ln \frac{\diff \fl^{w,t}_\pfwd \mssm}{\diff\mssm}(x) \diff\eta (x)}
\\
=&\diff_t\restr_{t=0} \int_\man \car_{(\eps,1]}(\eta_x) \ln \frac{\diff \fl^{w,t}_\pfwd \mssm}{\diff\mssm}(x) \diff\eta (x)
\\
=&\int_\man \car_{(\eps,1]}(\eta_x) \diff_t\restr_{t=0} \ln \frac{\diff \fl^{w,t}_\pfwd \mssm}{\diff\mssm}(x) \diff\eta (x)
\end{align*}
by Dominated Convergence Theorem. Finally, since~$\fl^{w,t}$ is orientation-preserving, $\frac{\diff \fl^{w,t}_\pfwd \mssm}{\diff\mssm}=\det \diff \fl^{w,t}$, whence
\begin{align*}
\diff_t\restr_{t=0}&\ \mbfR_\eps[\fl^{w,t}](\eta)=
\\
=&\ \int_\man \car_{(\eps,1]}(\eta_x) \diff_t\restr_{t=0} \ln \det \diff \fl^{w,t}(x) \diff\eta (x)
\\
=&\ \int_\man \car_{(\eps,1]}(\eta_x)\, \tr\quadre{\diff_t\restr_{t=0} \diff \fl^{w,t}(x)} \diff\eta (x)
\\
=&\ \int_\man \car_{(\eps,1]}(\eta_x)\, \div^\mssm_x w \diff\eta (x)
\\
=&\ \mbfB_\eps[w](\eta) \fstop
\end{align*}

\iref{i:c:Martingale2} By~\eqref{eq:l:IbP00}, Cauchy--Schwarz inequality,~\eqref{eq:l:DirDer0} and~\eqref{eq:DefGrad}
\begin{align*}
\abs{\int_\msP u \cdot v \cdot \mbfB_\eps[w] \diff\DF_\mssm}\leq&\ \int_\msP \tabs{\grad_w u\cdot v} \diff\DF_\mssm + \int_\msP \tabs{u \cdot \grad_w v} \diff\DF_\mssm
\\
\leq&\ \tnorm{\grad_w u}_{L^2(\DF_\mssm)} \norm{v}_{L^2(\DF_\mssm)} + \norm{u}_{L^2(\DF_\mssm)} \tnorm{\grad_w v}_{L^2(\DF_\mssm)}
\\
\leq&\ \norm{w}_{\Vect_{\n\mssm}} \norm{u}_{\mcE^{1/2}_1}\norm{v}_{\mcE^{1/2}_1} \fstop
\end{align*}

The existence and uniqueness of~$\mcA^w$ are a standard consequence.
\end{proof}
\end{cor}

Next, we show that~$(\mcE,\dom{\mcE}$ describes a truly infinite-dimensional diffusion. We refer to~\cite[Dfn.~2.9]{Hin09} for the concept of \emph{index} of a Dirichlet form.

\begin{prop}\label{p:Hino}
The form~$(\mcE,\dom{\mcE}$ has pointwise index~$p(\eta)=\infty$ $\DF_\mssm$-a.e.. Moreover, the index is `full' in the following sense:
For~$\DF_\mssm$-a.e.~$\eta\in\msP$ there exists an orthonormal basis~$\seq{e_i}_i$ of~$T^\Der_\eta\msP$ and a function~$u=u_i\in\hTF{2}{0}\subset\dom{\mcE}$ with~$\grad u=e_i$ for any choice of~$i$.
\begin{proof}
Since~$\DF_\mssm(\Ppaisofs)=1$, we can restrict our attention to~$\eta=\sum_i^\infty s_i\delta_{x_i}\in\Ppaisofs$. For all such~$\eta$ one has~$T^\Der_\eta\msP\cong \bigoplus_i^\perp (T_{x_i}\man,s_i\mssg)$, where~$\oplus^\perp$ denotes the orthogonal direct sum. 
For the rest of the proof we tacitly assume this identification.

As a basis for~$T^\Der_\eta\msP$ we fix~$\seq{e_{i,\ell}}_{i\in \N, \ell\leq d}$, where~$\seq{e_{i,\ell}}_{\ell\leq d}$ is a $\mssg$-orthonormal basis for~$T_{x_i}\man$ for every~$i$.
In order to show the second assertion, let~$f=f_{i,\ell}\in \mcC^\infty(\man)$ be so that~$e_{i,\ell}(f)_{x_i}=1$ and~$e_{i,\ell'}(f)_{x_i}=0$ for every~$\ell'\neq \ell$ and~$\varrho=\varrho_i\in\mcC^\infty(I)$ be so that~$\varrho(s_i)=1$ and~$\varrho(s_{i'})=0$ for every~$i\neq i'$.
The existence of~$f$ is standard, while the existence of~$\varrho$ follows from the fact that~$\eta\in \Ppaiso$, hence~$(\Proj{\To}{}\circ \boldPhi^{-1})(\eta)\in\Tso$.
Letting~$u=u_{i,\ell}\eqdef (f\otimes\varrho)^\trid$, one has~$\grad u(\eta)=e_{i,\ell}$.

By definition of~$\seq{e_{i,\ell}}_{i,\ell}$, one has~$\boldGamma(u_{i,\ell},u_{i',\ell'})(\eta)=\delta_{ii'}\gscal{e_{i,\ell}}{e_{i,\ell'}}=0$ for every~$(i,\ell)\neq(i',\ell')$. As a consequence, setting
\begin{align*}
\mbfA_{ii'}(\eta)\eqdef \tquadre{\boldGamma(u_{i,\ell},u_{i',\ell'})(\eta)}_{\ell\leq d}^{\ell'\leq d}\in \R^{d\times d} \comma \qquad \mbfA(\eta)\eqdef \tquadre{\mbfA_{ii'}(\eta)}_{i\leq n}^{i'\leq n}\in \R^{d^2\times n^2}\comma
\end{align*}
one has~$\mbfA(\eta)=\id_{\R^{nd}}$ for every~$\eta\in\Ppaisofs$.
Thus~$p(\eta)\geq \rk(\mbfA(\eta))=nd$ for every~$n$, which shows the first assertion.
\end{proof}
\end{prop}

\begin{rem}[A comparison with the Cheeger energy] A known object in metric measure space's analysis is the \emph{Cheeger energy} of a (complete and separable) metric measure space~$(Y,\mssd,\mssn)$
\begin{align*}
\Ch_{\mssd,\mssn}(f)\eqdef& \inf \set{\tfrac{1}{2}\nliminf \int_Y \slo{f_n}^2 \diff\mssn : f_n\in\Lip(Y,\mssd)\comma f_n\rar f \textrm{ in } L^2(Y,\mssn)} \comma
\\
\slo{f}(y)\eqdef& \limsup_{z\rar y} \frac{\abs{f(y)-f(z)}}{\mssd(y,z)}\comma \qquad f\in \Lip(Y,\mssd) \fstop
\end{align*}

A comparison of the Cheeger energy~$\Ch_{W_2,\DF_\mssm}$ of~$(\msP,W_2,\DF_\mssm)$ with the form $(\mcE,\dom{\mcE})$ constructed in Theorem~\ref{t:DForm} is here beyond our scope.
However, let us note that, at a merely heuristic level, we do not expect~$\Ch_{W_2,\DF_\mssm}$ to be a quadratic form.
Indeed, $\DF_\mssm$-a.e.~$\eta\in\msP$ is \emph{not} a regular measure in the sense of optimal transport (e.g.,~\cite{Gig11}), hence the tangent space at~$\eta$ accessed by Lipschitz functions is the full `abstract tangent space'~$\textrm{AbstrTan}_\eta$ \cite[Dfn.~3.7]{Gig11}.
By the results in~\cite[\S6]{Gig11},~$\Vect^\nabla_\eta$ embeds canonically, non-surjectively into~$\textrm{AbstrTan}_\eta$ and the latter is $\DF_\mssm$-a.e.\ \emph{not} a Hilbert space (rather, it is merely a Banach space).
Additionally, it is not clear to me whether~$\Ch_{W_2,\DF_\mssm}$ is non-trivial, that is, not identically vanishing.
\end{rem}

\section{The associated process}
In the case~$d\geq 2$, by e.g.,~\cite[Thm.~IV.5.1]{MaRoe92}, the form $(\mcE,\dom{\mcE})$ is properly associated with a $\DF_\mssm$-symmetric recurrent diffusion process
\begin{align}\label{eq:Eta}
\Eta_\bullet\eqdef \seq{\Omega, \mcF,\seq{\mcF_t}_{t\geq 0},\seq{\Eta_t}_{t\geq 0}, \set{P_\eta}_{\eta\in \msP}}
\end{align}
which we now characterize.

\subsection{Finite-dimensional approximations}\label{ss:FiniteDimApprox}
Everywhere in this section, assume~$d\geq 1$ whenever not explicit stated otherwise.
In order to identify the process~$\Eta_\bullet$ in~\eqref{eq:Eta}, we firstly discuss sequences of finite-dimensional approximations for the corresponding Dirichlet form~\eqref{eq:FormMain}.

\paragraph{Proofs} Firstly, for~$\eps\in (0,1)$ we consider the form $(\mcE,\hTF{2}{\eps})$. Note that
\begin{enumerate*}[label=$\boldsymbol{(}$\bfseries{\itshape\alph*}$\boldsymbol{)}$]
\item $(\mcE,\hTF{2}{\eps})$ is a closable, strongly local Dirichlet form, with closure~$(\mcE^{\eps},\dom{\mcE^{\eps}})$, closability and strong locality being inherited by~$(\mcE,\hTF{2}{0})$;
\item$(\mcE^{\eps},\dom{\mcE^{\eps}})$ is not densely defined on~$L^2(\msP,\DF_\mssm)$, yet it is densely defined on $L^2(\msP,\Bo{\eps}(\msP),\DF_\mssm)$ by Lemma~\ref{l:SAlg}\iref{i:l:SAlg3}.
\end{enumerate*}

In the same spirit as in Corollary~\ref{c:Martingale}, it is shown below that~$(\mcE^\eps,\dom{\mcE^\eps})$ is a $\DF_\mssm$-martingale-approximation of~$(\mcE,\dom{\mcE})$ as~$\eps$ tends to~$0$.
Furthermore, for each fixed~$n$, the form~$(\mcE^{1/n},\dom{\mcE^{1/n}})$ is a `finite-dimensional' approximation of~$(\mcE,\dom{\mcE})$, in following sense:
There exists a sequence~$\seq{Q_n}_n$ of compact manifolds with boundary, of finite increasing dimension, and a sequence of regular strongly local Dirichlet forms~$(\hat\mssE^n,\dom{\hat\mssE^n})$ on~$L^2(Q_n)$, so that the Dirichlet spaces corresponding to~$(\mcE^{1/n},\dom{\mcE^{1/n}})$ and~$(\hat\mssE^n,\dom{\hat\mssE^n})$ are isomorphic.

\begin{figure}[htb!]
{\normalsize
\begin{align*}
\xymatrix@R=40pt@C=80pt@M=7pt{
{\begin{matrix}(\mcE,\dom{\mcE}) \textrm{~on~} \\ L^2(\msP,\Bo{0}(\msP),\DF_\mssm)\end{matrix}}&
\\
{\begin{matrix}(\mcE^{1/n},\dom{\mcE^{1/n}}) \textrm{~on~}\\ L^2(\msP,\Bo{1/n}(\msP),\DF_\mssm) \end{matrix}}\ar[u]^{n\rar\infty} \ar@{<=>}[r]_{\qquad\textrm{$H$-isomorphism}}& {\begin{matrix}(\hat\mssE^n,\dom{\hat\mssE^n})\\ \textrm{~on~} L^2(Q_n)\end{matrix}} \ar@{.>}[ul]_{n\rar\infty}
}
\end{align*}
}
\caption{Finite-dimensional approximations of~$(\mcE,\dom{\mcE})$}
\label{fig:1}
\end{figure}

On the one hand, profiting the fact that~$(\mcE^{1/n},\dom{\mcE^{1/n}})$ and~$(\mcE,\dom{\mcE})$ are defined on the same $L^2$-space, we show the martingale convergence in the solid arrow in Fig.~\ref{fig:1}.
On the other hand, by standard techniques of stochastic analysis on manifolds, we identify the Markov diffusion process associated to~$(\hat\mssE^n,\dom{\hat\mssE^n})$.

The identification of the process associated to~$(\mcE,\dom{\mcE})$ follows, provided we show the convergence in the dotted arrow in Fig.~\ref{fig:1} in a suitable sense. Since we have a sequence of varying Hilbert spaces~$L^2(Q_n)$, the right notion of convergence ---~of such Hilbert spaces and of the Dirichlet forms defined on them~--- is the generalized (Hilbert, resp.\ Mosco) convergence introduced by K.~Kuwae and T.~Shioya in~\cite{KuwShi03}, which we briefly recall in the Appendix~\S\ref{ss:KuwaeShioya}. As it turns out, this notion of convergence is invariant under isomorphism of Hilbert spaces. Thus, the convergence in the dotted arrow in Fig.~\ref{fig:1} is a consequence of the martingale convergence in the solid arrow, which concludes the argument.

\paragraph{Projection manifolds} Let us now construct the manifolds~$Q_n$.

\begin{defs}[Simplices and projections]\label{d:Projections} For~$n\in \overline\N$,~$\eps\in I$ and~$\beta>0$ set
\begin{align*}
\Sigma^n\eqdef& \Proj{n}{}(\To)
\comma & \hman^n\eqdef& \Proj{n}{}(\widehat\mbfM)\comma
\\
\Sigma^n_\eps\eqdef& \set{\mbfs\in\Sigma^n : s_1\geq \eps} \comma & \hman^n_\eps\eqdef& \Sigma^n_\eps\times \man^{\times n}\comma
\end{align*}
each endowed with the usual topology and $\sigma$-algebra. We endow~$\Sigma^n$ (resp.~$\Sigma^n_\eps$) with (the restriction of) the probability measure
\begin{align}\label{eq:PDn}
\pi^n_\beta\eqdef \Proj{n}{\pfwd}\PD_\beta\comma
\end{align}
and~$\hman^n$ (resp.~$\hman^n_\eps$) with (the restriction of) the measure~$\widehat\mssm^n_\beta\eqdef \pi_\beta^n \otimes \n\mssm^n$. Finally, we set
\begin{equation}\label{eq:PhiN}
\begin{aligned}
\Phi_n\colon \hman^n & \longrar \Mp_1
\\
(\mbfs,\mbfx)&\longmapsto \sum_\ell s_\ell\delta_{x_\ell}
\end{aligned}
\fstop
\end{equation}
\end{defs}

Everywhere in the following, for fixed~$n\in \N$ and~$\eps>0$ let
\begin{equation}\label{eq:d:Products}
\begin{aligned}
H^n_\eps\eqdef& \cl_{L^2(\widehat\mssm^n_\beta)} \mcC^\infty\ttonde{\hman^n_\eps}=L^2\ttonde{\hman^n_\eps,\widehat\mssm^n_\beta} \comma & H_n\eqdef& H^n_{1/n}\comma
\\
H_{1/n}\eqdef& L^2(\msP,\mcB_{1/n}(\msP),\DF_\mssm)\comma & H\eqdef& L^2(\msP,\DF_\mssm)\fstop
\end{aligned}
\end{equation}

\paragraph{Superpositions}
Analogously to the case of~$\hM$ discussed in~\S\ref{ss:ProdMan}, we discuss here how to construct a process on~$\hman^n_\eps$. We start by constructing the corresponding Dirichlet form which turns out to be a \emph{superposition} ---~in the sense of~\cite[\S{V.3.1}]{BouHir91}~--- of a family of Dirichlet forms on~$\man^n$ and indexed by~$\Sigma^n_\eps$. For a systematic treatment of superpositions of Dirichlet forms in the language of Direct integrals we refer to~\cite{LzDS20}, the results of which we will use throughout this section.

\smallskip

For fixed~$n\in \N$ and~$\eps>0$ let
\begin{align*}
\hat\mssE^{n,\eps}(h)=& \int_{\Sigma^n_\eps} \int_{\man^{n,\mbfs}} \tabs{\nabla^{\mssg^{n,\mbfs},\mbfz}\restr_{\mbfz=\mbfx} h(\mbfz,\mbfs)}_{(\mssg^{n,\mbfs})_\mbfx}^2 \diff\n\mssm^n(\mbfx,\mbfs)\diff \pi^n_\beta(\mbfs) \comma \qquad h\in \mcC^1(\hman^n_\eps) \fstop
\end{align*}

For every~$i\leq n$ let further~$x^i_\bullet\eqdef \seq{x^i_t}_{t> 0}$ be independent Brownian motions on~$(\man,\mssg)$ starting at~$x^i_0$ and defined on a common probability space~$(\Omega, \mcF,\mbfP)$. We define a stochastic process~$\widehat\mssW^{n,\eps}$ on~$(\Omega, \mcF,\mbfP)$ with state space~$\hman^n_\eps$ by
\begin{align*}
\widehat\mssW^{n,\eps,\mbfx_0,\mbfs}_t(\omega)\eqdef \seq{x^1_{t/s_1}(\omega),\dotsc,x^n_{t/s_n}(\omega)}\comma\qquad \mbfs\eqdef\seq{s_1,\dotsc,s_n}\in\Sigma^n_\eps\comma \quad \omega\in\Omega \comma
\end{align*}
where~$\widehat\mssW^{n,\eps,\mbfx_0,\mbfs}_\bullet(\omega)$ is any stochastic path of~$\widehat\mssW^{n,\eps}_\bullet$ starting at~$(\mbfx_0,\mbfs)\in \hman^n_\eps$.

\begin{prop}\label{p:ProdForm}
The following hold true:
\begin{enumerate*}[label=$\boldsymbol{(}$\bfseries{\itshape\roman*}$\boldsymbol{)}$]
\item\label{i:p:ProdForm:1} the form~$\tseq{\hat\mssE^{n,\eps},\mcC^1(\hman^n_\eps)}$ is closable. Its closure~$\tseq{\hat\mssE^{n,\eps},\dom{\hat\mssE^{n,\eps}}}$ is a regular strongly local Dirichlet form on~$H^n_\eps$ with special core~$\mcC^1(\hman^n_\eps)$, and \item\label{i:p:ProdForm:3} semigroup
\end{enumerate*}
\begin{align*}
(\mssH^{n,\eps}_t h)(\mbfx,\mbfs)=&\ \ttonde{(\mssH^\mbfs_t\otimes \id_{H^n_\eps})\, h}(\mbfx,\mbfs)
\\
=&\ c_{n,\eps}\int_{\man^{\times n}} \prod_i^n \mssh_{t/s_i}(x_i,y_i)\,  h(y_1,\dotsc, y_n,s_1,\dotsc, s_n) \diff\n\mssm^n(\mbfy) \comma
\end{align*}
where~$c_{n,\eps}\eqdef \pi^n_\beta\Sigma^n_\eps\uparrow_{\eps\downarrow 0} 1$;
\begin{enumerate*}[resume*]
\item\label{i:p:ProdForm:4} $(\hat\mssE^{n,\eps},\dom{\hat\mssE^{n,\eps}})$ is properly associated to the process~$\widehat\mssW^{n,\eps}_\bullet$;
\item\label{i:p:ProdForm:5}
$\andi{\man^{\times n}}\times \Sigma^n_\eps$ is $\widehat\mssW^{n,\eps}_\bullet$-coexceptional.
\end{enumerate*}
\begin{proof}
\iref{i:p:ProdForm:1}--\iref{i:p:ProdForm:4} follow by~\cite[Ex.~3.12]{LzDS20}.
\iref{i:p:ProdForm:5} follows by showing that~$\andi{\man^{\times n}}$ is~$\mssW^{n,\mbfs}_\bullet$-coexceptional for $\pi^n_\beta$-a.e.~$\mbfs\in \Sigma^n_\eps$, Lemma~\ref{l:Coexcept}, together with the capacity estimate~\cite[Ex.~3.12(iii)]{LzDS20}.
\end{proof}
\end{prop}

\begin{lem}\label{l:NestD}
Let~$D_n\eqdef \Phi_n^* \hTF{1}{1/n}$ and~$\seq{T^n_m}_m$ be a family of closed sets in~$\Sigma^n_{1/n}\cap\Proj{n}{}(\Tso)$ such that~$T^n_m\uparrow T^n$ with~$T^n$ of full~$\pi^n_\beta$-measure in~$\Sigma^n_{1/n}$.
Further let~$\seq{F^n_m}_{m}$ be a nest for~$\mssE^{n,\mbfs}$ for every~$\mbfs\in T^n$ and such that~$F^n_m\subset \inter F^n_{m+1}$. Set~$\hat F^n_m\eqdef T^n\times F^n_m\subset \hman^n_{1/n}$ and let
\begin{align*}
D_{n,m}=(D_n)_{\hat F^n_m}\eqdef \set{u\in D_n : u\equiv 0\;\; \pi^n_\beta\text{-a.e. on } (\hat F^n_m)^\complement}\fstop
\end{align*}

Then,~$\cup_m D_{n,m}$ is both dense in~$H_n$ and dense in~$\dom{\hat\mssE^n}$.

\begin{proof}
It suffices to show the second density statement. 
Set
\begin{align*}
(D_{n,m})_{\mbfs}\eqdef \set{h\ttonde{\emparg,\mbfs} : h\in D_{n,m}}\fstop
\end{align*}
In order to show that~$\cup_m D_{n,m}$ is dense in~$\dom{\hat\mssE^n}$, it suffices to show that~$\cup_m (D_{n,m})_\mbfs$ is dense in~$\dom{\mssE^{n,\mbfs}}$ for~$\pi^n_\beta$-a.e.~$\mbfs\in \Sigma^n_{1/n}$.
Since~$\seq{F^n_m}_{m}$ is a nest for~$(\mssE^{n,\mbfs},\dom{\mssE^{n,\mbfs}})$ for all~$\mbfs$ in the set of full $\pi^n_\beta$-measure~$T^n$, we have that~$\cup_m\dom{\mssE^{n,\mbfs}}_{F^n_m}$ is dense in~$\dom{\mssE^{n,\mbfs}}$ for all~$\mbfs\in T^n$. Thus, it suffices to show that
\begin{align}\label{eq:l:NestD1}
\cl_{(\mssE^{n,\mbfs})_1^{1/2}}\ttonde{(D_{n,m+1})_\mbfs}\supset \dom{\mssE^{n,\mbfs}}_{F^n_m}\comma\qquad m\in \N \fstop
\end{align}

To this end, we firstly show that~$(D_{n,m+1})_\mbfs\supset\mcA_{n,m+1}\eqdef\ttonde{\mcC^1(\man)^{\otimes n}}_{F^n_{m+1}}$. 
Indeed, since, in particular,~$\mbfs\in \Proj{n}{}(\Tso)$, for~$\ell\leq n$ there exists~$\varrho_\ell\in\Cc^\infty((1/n,1))$, depending on~$\mbfs$, so that~$\varrho_{\ell_1}(s_{\ell_2})=s_{\ell_1}^{-1}\delta_{\ell_1 \, \ell_2}\in (0,\infty)$, where~$\delta$ denotes Kronecker's delta. Thus, for any choice of~$\seq{f_\ell}_\ell^n\subset \mcC^1(\man)$, one has that
\begin{align*}
\bigotimes_\ell^n f_\ell=\tonde{\Phi^*\tonde{\prod_\ell^n (f_\ell \, \varrho_\ell)^\trid}}(\emparg,\mbfs)\comma\qquad \mbfs\in T^n \fstop
\end{align*}

Finally, it is clear that
\begin{align*}
\cl_{(\mssE^{n,\mbfs})_1^{1/2}}(\mcA_{n,m+1})=\cl_{(\mssE^{n,\mbfs})_1^{1/2}}\ttonde{\mcC^1(\man^{\times n})_{F^n_{m+1}}}\supset \dom{\mssE^{n,\mbfs}}_{F^n_m}\comma
\end{align*}
where the latter inclusion follows by a localization argument with smooth partitions of unity and regularization by convolution since~$F^n_m\subset\inter F^n_{m+1}$. We omit the details. This concludes the proof of~\eqref{eq:l:NestD1}.
\end{proof}
\end{lem}

We are now ready to show the isomorphism of forms in Figure~\ref{fig:1}. For notational simplicity
let~$(\hat\mssE^n,\dom{\hat\mssE^n})\eqdef (\hat\mssE^{n,1/n},\dom{\hat\mssE^{n,1/n}})$ be defined as in Proposition~\ref{p:ProdForm}.

\begin{lem}\label{l:Intertwining}
The forms~$(\hat\mssE^n,\dom{\hat\mssE^n})$ and~$(\mcE^{1/n},\dom{\mcE^{1/n}})$ are intertwined via~$\Phi_n^*$ as in~\eqref{eq:PhiN}.

\begin{proof} Let~$\eta\in \Ppa$ be of the form~$\eta=\sum_\ell^N s_\ell\delta_{x_\ell}$ for some~$N\in \overline\N$. For~$u=F\circ\hat\mbff\in \hTF{0}{\eps}$ define
\begin{align}
\label{eq:Shorthand}
\sum_\ell^N s_\ell\, \hat\mbff(x_\ell,s_\ell)\eqdef& \seq{\sum_\ell^N s_\ell \hat f_1(x_\ell,s_\ell),\cdots,\sum_\ell^N s_\ell \hat f_k(x_\ell,s_\ell)} \comma \quad N\in \overline\N
\\
\label{eq:t:Margin0}
U_n\colon u\mapsto& \tonde{\Phi_n^*u\colon(\mbfs,\mbfx)\mapsto F\tonde{\sum_\ell^n s_\ell\, \hat\mbff(x_\ell,s_\ell)} } \comma \quad n\in \N\comma (\mbfs,\mbfx)\in \Sigma^n \fstop
\end{align}

By definition of~$u$ one has
\begin{align}\label{eq:t:Margin00}
\hat f_i(x,s)=0\comma \qquad x\in \man\comma s\leq 1/n\comma i\leq k\comma
\end{align}
hence
\begin{align}\label{eq:IntertwiningN}
\boldPhi^*u = \Phi_n^* u \circ \Proj{}{n}\comma \qquad u\in\hTF{0}{1/n} \fstop
\end{align}

Now, by Proposition~\ref{p:SuppDF}\iref{i:p:SuppDF2} and the fact that~$\mcE^{1/n}=\mcE$ on~$\hTF{1}{1/n}$ one has for all~$u\in\hTF{1}{1/n}$
\begin{align}\label{eq:t:Margin1}
\mcE^{1/n}(u)=& \int_{\andihM} \sum_{i,p}^k (\partial_i F \cdot \partial_p F)\tonde{\sum_\ell^\infty s_\ell\, \hat\mbff(x_\ell,s_\ell)} \sum_\ell^\infty s_\ell\, \Gamma^\mssg(\hat f_i,\hat f_p)(x_\ell,s_\ell) \diff\hvolm_\beta(\mbfs,\mbfx) \fstop
\end{align}

If~$\ell>n$, then~$s_\ell\leq 1/n$ because~$s_1\geq s_2\geq \dotsc$. Thus, by~\eqref{eq:t:Margin1} and~\eqref{eq:Rescal},
\begin{align}
\nonumber
\mcE^{1/n}(u)=&\int_{\hman^n} \sum_{i,p}^k(\partial_i F\cdot \partial_p F)\tonde{\sum_\ell^n s_\ell\, \hat\mbff(x_\ell,s_\ell)} \sum_\ell^n s_\ell^2\, \Gamma^{s_\ell\mssg}(\hat f_i,\hat f_p)(x_\ell,s_\ell) \diff \widehat\mssm^n_\beta(\mbfs,\mbfx)
\\
\label{eq:t:Margin2}
=&\int_{\hman^n} \abs{\nabla^{\mssg^{n,\mbfs},\mbfz}\restr_{\mbfz=\mbfx} \, {F\tonde{\sum_\ell^n s_\ell\, \hat\mbff(z_\ell,s_\ell)}}}^2_{\mssg^{n,\mbfs}_\mbfx} \diff \widehat\mssm^n_\beta(\mbfs,\mbfx)
\\
\label{eq:t:Margin3}
=&\int_{\hman^n_{1/n}} \abs{\nabla^{\mssg^{n,\mbfs},\mbfz}\restr_{\mbfz=\mbfx} \, {F\tonde{\sum_\ell^n s_\ell\, \hat\mbff(z_\ell,s_\ell)}}}^2_{\mssg^{n,\mbfs}_\mbfx} \diff \widehat\mssm^n_\beta(\mbfs,\mbfx)
\\
\nonumber
=&\hat\mssE^n(\Phi_n^* u)\comma
\end{align}
where we may reduce the domain of integration in~\eqref{eq:t:Margin2} to the one in~\eqref{eq:t:Margin3} since the integrand vanishes identically on~$\hman^n\setminus \hman^n_{1/n}$ for all~$u\in \hTF{1}{1/n}$, again as a consequence of~\eqref{eq:t:Margin00}.
In particular,~$D_n\eqdef U_n\, \hTF{1}{\eps} \subset \dom{\hat\mssE^n}$. An analogous computation shows that~$\Phi_n^*u$ satisfies
\begin{equation*}
\norm{u}_{H_{1/n}}=\norm{\Phi_n^* u}_{H_n}\comma \qquad u\in \hTF{0}{1/n}\fstop
\end{equation*}

The family $\hTF{1}{1/n}$ is dense in~$H_{1/n}$ by Lemma~\ref{l:SAlg}\iref{i:l:SAlg3} and dense in~$\dom{\mcE^{1/n}}$ by definition of the latter. As a consequence, the operator~$U_n$ defined in~\eqref{eq:t:Margin0} uniquely extends to a non-relabeled isometric operator~$U_n\colon \dom{\mcE^{1/n}}_1\rar \dom{\hat\mssE^n}_1$, and, subsequently to an isometric operator~$U_n\colon H_{1/n}\rar H_n$.
It suffices to show the intertwining property on dense subsets. Thus, the conclusion follows by showing that~$D_n\eqdef U_n\, \hTF{1}{1/n}=\Phi_n^*\hTF{1}{1/n}$ is both dense in~$H_n$ and dense in~$\dom{\hat\mssE^n}$. This follows by Lemma~\ref{l:NestD} with~$F_m\eqdef \man$ for every~$m$.
\end{proof}
\end{lem}

\begin{prop}\label{p:MoscoEn}
Let~$\tseq{\hat\mssE^n,\dom{\hat\mssE^n}}\eqdef \tseq{\hat\mssE^{n,1/n},\dom{\hat\mssE^{n,1/n}}}$ be as in Pro\-position~\ref{p:ProdForm}. Then, the sequence~$\tseq{\hat\mssE^n, \dom{\hat\mssE^n}}_n$ Mosco converges to $(\mcE,\dom{\mcE})$ in the sense of Definition~\ref{d:Kuwae--Shioya}.

\begin{proof} Recall the notation in~\eqref{eq:d:Products}.
We claim that~$(\mcE^{1/n},\dom{\mcE^{1/n}})$ converges to~$(\mcE,\dom{\mcE})$ as~$n\rar \infty$ in the generalized Mosco sense.
Indeed let~$P_n\colon H\rar H$ be the projection operator~$P_n\eqdef \mbfE_{\DF_\mssm}[\,\emparg\!\mid \mcB_{1/n}(\Ppa)]$ given by the conditional expectation w.r.t.~$\mcB_{1/n}(\Ppa)$.
By definition,~$H_{1/n}=P_n(H)$. Since~$\mcB_0(\msP)_{\Ppa}=\mcB_\mrmn(\Ppa)$ by Lemma~\ref{l:SAlg}\iref{i:l:SAlg2}, the sequence~$\seq{P_n}_n$ converges strongly to~$\id_H$.
Regard~$(\mcE^{1/n},\dom{\mcE^{1/n}})$ as a (\emph{not} densely defined) quadratic form on~$H$. By Lemma~\ref{l:ConvHStandard} applied to the family~$\seq{P_n}_n$, it suffices to show the Mosco convergence of the form $(\mcE^{1/n},\dom{\mcE^{1/n}})$ to $(\mcE,\dom{\mcE})$ in the classical sense.
The strong $\Gamma$-$\limsup$ condition~\eqref{eq:GammaLimSup} is the content of Corollary~\ref{c:Domains}. 
Since~$\mcE(u)=\mcE^{1/n}(u)$ for every~$u\in\dom{\mcE^{1/n}}$, the weak $\Gamma$-$\liminf$ condition~\eqref{eq:GammaLimInf} is a consequence of the weak lower semi-continuity of~$\mcE$, Lemma~\ref{l:WLSC}.

Now, Lemma~\ref{l:UniEquiv} applies to~$Q=\mcE^{1/n}$ and~$Q^\sharp=\hat\mssE^n$ with~$U_n=\Phi_n^*$ as in (the proof of) Lemma~\ref{l:Intertwining}. Therefore,~$(\hat\mssE^n,\dom{\hat\mssE^n})$ Mosco converges to~$(\mcE,\dom{\mcE})$ as~$n\rar \infty$.\end{proof}
\end{prop}

\begin{prop}\label{p:ESA} The operator~$(-\mbfL,\hTF{\infty}{0})$ is essentially self-adjoint.
\begin{proof} We show that for every~$u_0\in \hTF{\infty}{0}$ and every~$T\in[0,\infty)$ there exists
\begin{enumerate*}[label=$\boldsymbol{(}$\bfseries{\itshape\alph*}$\boldsymbol{)}$]
\item\label{i:p:ESA:1} a sequence $\seq{u_{0,n}}_n\subset \hTF{\infty}{0}$ such that~$L^2(\msP,\DF_\mssm)$-$\nlim u_{0,n}=u_0$
and 
\item\label{i:p:ESA:2} strong solutions~$u_n$ of the Cauchy problems
\end{enumerate*}
\begin{equation}\label{eq:BackHeat}
\begin{cases}
(\diff_t u_n)(t)-(-\mbfL u_n)(t)=0
\\
u_n(T)=u_{0,n}
\end{cases}
\qquad u_n(t)\in\hTF{2}{0}\subset\dom{\mbfL}\comma \qquad t\in[0,T] \fstop
\end{equation}
Then, the assertion follows by~\cite[\S{II.5},~Thm.~1.10, p.~30]{BerKon95}. Condition~$(ii)$ there is trivially satisfied since we chose, in the notation of~\cite[\emph{ibid.}]{BerKon95},~$A_n=A$.

\paragraph{Solutions to the heat equation} For~$n\leq N\in \overline{\N}$ let~$\eta=\sum_i^N s_i\delta_{x_i}$ and~$u\in \hTF{\infty}{1/n}$, and recall the notation in~\eqref{eq:Shorthand}. Analogously to the proof of Lemma~\ref{l:Intertwining}, by~\eqref{eq:L1L2} and~\eqref{eq:GammaGammaStar} we have
\begin{align}
2\,(\mbfL u_0)(\eta)=&\ \sum_{i,p}^k (\partial_{ip}^2 F)\tonde{\sum_\ell^n s_\ell\hat\mbff(x_\ell,s_\ell)} \sum_{\ell}^n s_\ell\, \Gamma^\mssg(\hat f_i,\hat f_p)(x_\ell,s_\ell)
\\
\nonumber
&\ +\sum_i^k (\partial_i F) \tonde{\sum_\ell^n s_\ell \hat\mbff(x_\ell,s_\ell)} \sum_{\ell}^n \Delta^z\restr_{z=x_\ell} \hat f_i(z,s_\ell)
\\
\nonumber
=&\ \sum_{i,p}^k (\partial_{ip}^2 F)\tonde{\sum_\ell^n s_\ell \hat\mbff(x_\ell,s_\ell)} \sum_{\ell}^n s_\ell^2\,\Gamma^{s_\ell \mssg}(\hat f_i,\hat f_p)(x_\ell,s_\ell)
\\
\nonumber
&\ +\sum_i^k (\partial_i F) \tonde{\sum_\ell^n s_\ell \hat\mbff(x_\ell,s_\ell)} \sum_{\ell}^n s_\ell \,\Delta^{s_\ell \mssg,z}\restr_{z=x_\ell} \hat f_i(z,s_\ell)
\\
\label{eq:LtoLap}
=&\ (\Delta^{\mssg^{n,\mbfs},\mbfz}\restr_{\mbfz=\mbfx} \Phi_n^*u_0)(\mbfz,\mbfs)\comma
\end{align}
%
By~\eqref{eq:LtoLap} together with the time-reversal~$t\mapsto T-t$, the Cauchy problem~\eqref{eq:BackHeat} with~$u_0$ in place of~$u_{0,n}$ transforms into the Cauchy problem
\begin{equation}\label{eq:HeatProd}
\begin{aligned}
\partial_t h-\tfrac{1}{2}\Delta^{\mssg^{n,\mbfs}} h=&0\comma
\\
h(0)=&\Phi_n^*u_0\comma
\end{aligned}\qquad t\in[0,T] \fstop
\end{equation}

Since~$\man^{n,\mbfs}$ is a closed manifold, by standard results the latter Cauchy problem has a unique solution, say~$t\mapsto h(t)$, additionally satisfying~$h(t)\in\mcC^\infty(\man^{n,\mbfs})$ for all~$t\in [0,T]$. Finally, note that every function~$h\in\mcC^\infty(\man^{n,\mbfs})$ may be written as
\begin{align*}
h(\mbfx)=(U_n v)(\mbfx,\mbfs)
\end{align*}
for some~$v=G\circ\hat\mbfg^\trid\in\hTF{\infty}{1/n}$ (not necessarily in minimal form) with~$G\in \Cc^\infty(\R^{nd};\R)$. As a consequence, there exist functions~$t\mapsto G(t)\in \Cc^\infty(\R^{nd};\R)$ and~$t\mapsto \hat g_i(t)\in\Cc^\infty(\hman_{1/n})$ for~$i\leq nd$ such that~$h(t)=\Phi_n^* u(t)$, where~$u(t)\eqdef G(t)\circ\hat \mbfg(t)^\trid$.
We have thus constructed the unique solution~$t\mapsto u(t)$ of the Cauchy problem~\eqref{eq:BackHeat} with initial data~$u_0\in \hTF{\infty}{1/n}$, additionally satisfying~$u(t)\in \hTF{\infty}{1/n}$. As usual, the representation of~$u$ by~$G$ and~$\hat\mssg$ is not unique, cf.~Rmk.~\ref{r:ReprCyl}. Note that the strong solution~$h_{u(t)}$ to~\eqref{eq:HeatProd} is smooth, hence the corresponding function~$u(t)\in\hTF{\infty}{1/n}$ is a strong solution to~\eqref{eq:BackHeat} in the sense of the strong topology of~$L^2(\msP,\DF_\mssm)$ and therefore satisfies~\iref{i:p:ESA:2}. 

\paragraph{Approximations} Let~$u_{0,n}\in \hTF{\infty}{1/n}$ be given by Corollary~\ref{c:Domains}, thus satisfying~\iref{i:p:ESA:1}. Constructing solutions~$t\mapsto u_n(t)$ to the Cauchy problems~\eqref{eq:BackHeat} as above concludes the proof.
\end{proof}
\end{prop}

\subsection{A quasi-homeomorphic Dirichlet form}
In this section, assume~$d\geq 2$ whenever not explicitly stated otherwise.

\smallskip

The coincidence of the `finite-dimensional' approximating forms in Lemma~\ref{l:Intertwining} is not sufficient to establish the identification of the process~$\Eta_\bullet$.
To this end, we need to construct a Dirichlet form on~$L^2(\hM,\hvolm_\beta)$ quasi-homeomorphic to~$\mcE$ in the sense of~\cite[Dfn.~3.1]{CheMaRoe94}, for which we can characterize the associated process.
Namely, the superposition~$(\widehat\dirE,\dom{\widehat\dirE})$ of $(\dirE^\mbfs,\dom{\dirE^\mbfs})$ indexed by~$\mbfs\in(\To,\PD_\beta)$, that is, the Dirichlet form on~$L^2(\hM,\hvolm_\beta)$ with semigroup~$\widehat\Heat_\bullet$ defined as in~\eqref{eq:IntroHatHeat}.

\begin{figure}[htb!]
{\normalsize
\begin{align*}
\xymatrix@R=40pt@C=80pt@M=7pt{
{\begin{matrix}(\mcE,\dom{\mcE}) \textrm{~on~} \\ L^2(\msP,\Bo{0}(\msP),\T_\mrma,\DF_\mssm)\end{matrix}} \ar@{<=>}[r]^{\qquad\boldPhi^*}& {\begin{matrix}(\widehat\dirE,\dom{\widehat\dirE}) \textrm{~on~}\\ L^2(\hM,\T_\mrmp,\hvolm_\beta) \end{matrix}}
\\
{\begin{matrix}(\mcE^{1/n},\dom{\mcE^{1/n}}) \textrm{~on~}\\ L^2(\msP,\Bo{1/n}(\msP),\DF_\mssm) \end{matrix}}\ar[u]^{n\rar\infty} \ar@{<=>}[r]_{\qquad\Phi_n^*}& {\begin{matrix}(\hat\mssE^n,\dom{\hat\mssE^n})\\ \textrm{~on~} L^2(Q_n)\end{matrix}} \ar@{.>}[ul]_{n\rar\infty}\ar@{.>}[u]_{n\rar\infty}
}
\end{align*}
}
\caption{Quasi-homeomorphism of forms}
\label{fig:2}
\end{figure}

\begin{thm}\label{t:QuasiHomeo} The map~$\boldPhi$ in~\eqref{eq:Intro:Phi} is a quasi-homeomorphism intertwining
\begin{align*}
(\widehat\dirE,\dom{\widehat\dirE}) \text{~~on~~}
L^2(\hM,\T_\mrmp,\hvolm_\beta) \qquad \text{and}\qquad (\mcE,\dom{\mcE}) \text{~~on~~} L^2(\Ppaiso,\T_\mrma,\DF_\mssm) \fstop
\end{align*}

\begin{proof}
For~$i<j\in \N$ set~$\mbfU_{i,j,\delta}\eqdef \set{\mbfx\in\mbfM: \mssd_\mssg(x_i,x_j)<\delta}$. Note that~$\mbfU_{i,j,\delta}$ is open and that for every~$\eps>0$ there exists~$\delta(\eps)>0$ such that $\n\volm\, \mbfU_{i,j,\delta(\eps)}<\eps$. As a consequence, the set
\begin{align*}
\mbfU_m\eqdef \bigcup_{\substack{i,j \, : \, i<j}} \mbfU_{i,j,\delta(2^{-i-j}/m)}\subset \mbfM\comma \qquad m\in \N\comma
\end{align*}
is open, relatively compact, and satisfies~$\n\volm\, \mbfU_m\leq 1/m$,~$\overline{\mbfU_{m+1}}\subset \mbfU_m$ and~$\mbfU_m \downarrow_m \andi{\mbfM}^\complement$.
Finally, set~$\mbfF_m\eqdef \mbfU_m^\complement$ and note that~$\mbfF_m$ is compact and satisfies~$\mbfF_m\uparrow_m \andi{\mbfM}$.

\paragraph{A nest for~$\widehat\dirE$} The set~$\andi{\mbfM}$ is~$\dirE^\mbfs$-coexceptional by Lemma~\ref{l:Except}. Since~$\Cap_\mbfs\eqdef \Cap_{\dirE^\mbfs}$ is a Choquet capacity,
\begin{align*}
\mlim \Cap_\mbfs(\mbfF_m^\complement)\leq& \mlim \Cap_\mbfs\ttonde{\overline{\mbfF_m^\complement}} =\inf_m\Cap_\mbfs\ttonde{\overline{\mbfU_m}}
\\
=&\ \Cap_\mbfs \ttonde{\cap_m \overline{\mbfU_m}}=\Cap_\mbfs\ttonde{\cap_m \mbfU_m}=\Cap_\mbfs(\andi{\mbfM}^\complement)=0\comma
\end{align*}
hence~$\seq{\mbfF_m}_{m}$ is a nest for~$(\dirE^\mbfs,\dom{\dirE^\mbfs})$ for every~$\mbfs\in\To$.
Set now
\begin{align*}
\To_m\eqdef&\tset{\mbfs\in\To : s_\ell-s_{\ell+1}\geq 2^{-\ell-1}/m}
\end{align*}
and note that~$\To_m$ is compact (closed) and satisfies~$\To_m\uparrow_m \Tso$. Then,~$\widehat \mbfF_m\eqdef \To_m\times \mbfF_m$ is compact (closed) in~$\andihM$ (in~$\widehat\mbfM$) and satisfies~$\widehat \mbfF_m\uparrow_m \andihM$. The family~$\tseq{\widehat \mbfF_m}_m$ is a nest for~$(\widehat\dirE,\dom{\widehat\dirE})$ by~\cite[Ex.~3.12(iii)]{LzDS20}.

\paragraph{A nest for~$\mcE$}
Since~$\boldPhi\colon (\andihM,\T_\mrmp)\rar (\msP,\T_\mrma)$ is a homeomorphism onto its image~$\Ppaiso$ by Proposition~\ref{p:SuppDF}\iref{i:p:SuppDF2}, then~$G_m\eqdef \boldPhi(\widehat \mbfF_m)\subset\msP$ is itself compact in~$(\Ppaiso,\T_\mrma)$, hence compact in~$(\msP,\T_\mrma)$ and, in turn, compact in~$(\msP,\T_\mrmn)$ by Proposition~\ref{p:EthierKurtz}\iref{i:p:EthierKurtz0}; also cf.~\cite[Lem.~2.4]{EthKur94}. Set 
\begin{align*}
\ttonde{\hTF{2}{\eps}}_{G_m}\eqdef \tset{u\in\hTF{2}{\eps} : u\equiv \zero \, \DF_\mssm\textrm{-a.e. on } G_m^\complement}\subset \dom{\mcE}_{G_m}\comma \qquad \eps\in I\comma
\end{align*}
and note that~$\ttonde{\hTF{2}{0}}_{G_m}\subset \mcC(\msP,\T_\mrma)$ for every~$m\in \N$ by Remark~\ref{r:MeasContCyl}\iref{i:r:MeasContCyl:3}. Then, in order to prove that~$\seq{G_m}_{m}$ is a nest for~$(\mcE,\dom{\mcE})$ we need to show that~$\cup_m\ttonde{\hTF{2}{0}}_{G_m}$ is dense in~$\dom{\mcE}$.

We start by reducing the statement to a finite-dimensional case.
In fact, by Corollary~\ref{c:Domains}, it suffices to show that~$\mcC_n\eqdef\cup_m\ttonde{\hTF{2}{1/n}}_{G_m}$ is dense in~$\hTF{2}{1/n}$ for~$n\in\N$. To this end, fix~$n\in \N$ and let~$u\eqdef F\circ\hat\mbff^\trid$ be arbitrary in~$\hTF{2}{1/n}$.
By definition of~$\hTF{0}{1/n}$, it holds that~$\boldPhi^*u=(\Phi_n^*u)\circ\Proj{n}{}$ $\hvolm_\beta$-a.e.. Therefore, by Lemma~\ref{l:Intertwining} it suffices to establish that~$D_n'\eqdef \Phi^*\mcC_n$ is dense in~$\dom{\hat \mssE^n}$.

To this end, set~$T^n_m\eqdef \Proj{n}{}(\To_m)$,~$F^n_m\eqdef \Proj{n}{}(\mbfF_m)$ and~$\hat F^n_m\eqdef T^n_m\times F^n_m$.
Since~$\boldPhi$ is a homeomorphism onto~$G_m$ for every~$m$, one has~$\Phi_n^*\ttonde{\ttonde{\hTF{1}{1/n}}_{G_m}}=\ttonde{\Phi_n^*\ttonde{\hTF{1}{1/n}}}_{\hat F^n_m}$.
Thus, the conclusion follows by Lemma~\ref{l:NestD} with~$T^n_m$ and~$F^n_m$ as above.

\paragraph{Intertwining} It suffices to prove the intertwining property~$\widehat\dirE\circ \boldPhi^*=\mcE$ for all~$u\in\mcC$ with~$\mcC$ dense in~$\dom{\mcE}$ and~$\boldPhi^*\mcC$ dense in~$\dom{\widehat\dirE}$. 
We choose~$\mcC\eqdef \hTF{1}{0}=\cup_n \hTF{1}{1/n}$, cf.~Rmk.~\ref{r:MeasContCyl}\iref{i:r:MeasContCyl:8}. The first density requirement follows by definition of~$(\mcE,\dom{\mcE})$.

By standard topological facts and Lemma~\ref{l:Intertwining}, one has
\begin{align*}
\cl_{\widehat\dirE_1^{1/2}}\tonde{\boldPhi^*\bigcup_n \hTF{1}{1/n}}\supset&\ \bigcup_n\cl_{\widehat\dirE_1^{1/2}} \ttonde{\boldPhi^* \hTF{1}{1/n}}
=\bigcup_n\cl_{\widehat\dirE_1^{1/2}} \ttonde{\Phi_n^* \hTF{1}{1/n}\circ \Proj{n}{}}
\\
=&\ \bigcup_n\cl_{(\hat\mssE^n)_1^{1/2}}\ttonde{\Phi_n^* \hTF{1}{1/n}}\circ \Proj{n}{}
=\bigcup_n \dom{\hat\mssE^n}\circ \Proj{n}{} \fstop
\end{align*}
As a consequence, it suffices to show that~$\cup_n \dom{\hat\mssE^n}\circ \Proj{n}{}$ is dense in~$\dom{\widehat\dirE}$.
By our usual reduction argument, it suffices to show that~$\cup_n \dom{\mssE^{n,\mbfs}}\circ\Proj{n}{}$ is dense in~$\dom{\dirE^\mbfs}$ for $\PD_\beta$-a.e.~$\mbfs\in\To$.
This is however immediate by definition of~$(\dirE^\mbfs,\dom{\dirE^\mbfs})$, since~$\cup_n \dom{\mssE^{n,\mbfs}}\circ\Proj{n}{}\supset \Cyl{\infty}$ as in~\eqref{eq:Cyl}.

As for the intertwining, for all~$u\in \hTF{1}{1/n}$ it holds by~\eqref{eq:IntertwiningN} that~
\begin{align}\label{eq:t:QuasiHomeo:1}
\widehat\dirE(\boldPhi^*u)=\widehat\dirE(\Phi_n^*u \circ \Proj{}{n}) \fstop
\end{align}
Noting that~$\Phi^*_n u\circ \Proj{}{n}\in\Cyl{1}$ by definition of~$\hTF{1}{1/n}$, it follows by definition of~$\widehat\dirE$ and~\eqref{eq:ADKBSC0} that
\begin{align}\label{eq:t:QuasiHomeo:2}
\widehat\dirE(\Phi_n^*u \circ \Proj{}{n})=\hat\mssE^n(\Phi^*_n u) \fstop
\end{align}
Respectively by Lemma~\ref{l:Intertwining} and definition of~$(\mcE^{1/n},\dom{\mcE^{1/n}})$,
\begin{align}\label{eq:t:QuasiHomeo:3}
\hat\mssE^n(\Phi^*_n u)=\mcE^{1/n}(u)=\mcE(u) \fstop
\end{align}
Finally, combining~\eqref{eq:t:QuasiHomeo:1}--\eqref{eq:t:QuasiHomeo:3} concludes the proof of the intertwining property.
\end{proof}
\end{thm}

As a consequence of the regularity of~$(\widehat\dirE,\dom{\widehat\dirE})$ on~$L^2(\hM,\T_\mrmp,\hvolm_\beta)$ and of~\cite[Thm.~3.7]{CheMaRoe94} we have

\begin{cor}\label{c:QuasiHomeo} The Dirichlet form~$(\mcE,\dom{\mcE})$ on $L^2(\Ppaiso,\DF_\mssm)$ is $\T_\mrma$-quasi-regular.
\end{cor}

\begin{rem} Although~$\boldPhi^*\colon L^2(\msP,\DF_\mssm)\rar L^2(\widehat\mbfM,\widehat\volm_\beta)$ is an order isomorphism, Theorem~\ref{t:QuasiHomeo} does not follow from the general result~\cite[Thm.~3.12]{LenSchWir18}, where the intertwined quasi-regular Dirichlet forms are additionally assumed irreducible. We postpone a study of the $\mcE$-invariant sets to Theorem~\ref{t:PropertiesEta} below.
\end{rem}

\subsection{Properties of~\texorpdfstring{$\Eta_\bullet$}{Eta}}
In this section we collect some properties of the process~$\Eta_\bullet$.

\subsubsection{Group actions}\label{ss:Groups}
In order to discuss invariant sets, we start with some preliminaries about group actions on general spaces. Let~$(Y,\T)$ be a Polish space with Borel $\sigma$-algebra $\Bo{\T}$, and~$\mssn$ be a probability measure on~$(Y,\Bo{\T})$.

\paragraph{Preliminaries} Let~$G$ be a group acting on~$Y$, write~$G\acts Y$ and~$g.y\in Y$ for any~$g\in G$ and~$y\in Y$. We denote by~$\quotient{Y}{G}$ the quotient of~$Y$ by the action of~$G$, always endowed with the quotient topology and the induced Borel $\sigma$-algebra, and by $\pr^G\colon Y\rar \quotient{Y}{G}$ the projection to the quotient. 

We say that~$A\subset Y$ is $G$-invariant if~$G.A\eqdef\set{g.y: g\in G,y\in A}\subset A$ (equivalently~$G.A=A$) and that~$f\colon Y\rar \R$ is $G$-invariant if it is constant on $G$-orbits, i.e.~$f(y)=f(g.y)$ for every~$g\in G$,~$y\in Y$. We say that~$A\subset Y$ is $(G,\mssn)$-invariant if there exists a $G$-invariant~$A_1\in\Bo{\T}$ such that~$A\triangle A_1$ is $\mssn$-negligible. If~$A\in\Bo{\T}$ and~$(Y,\Bo{\T},G,\mssn)$ is a continuous dynamical system with invariant measure~$\mssn$, then $(G,\mssn)$-invariance coincides with the standard definition~\cite[Eqn.~(1.2.13)]{DaPZab96}. 
For any~$N\in \overline{\N}$ let
\begin{enumerate*}[label=$\boldsymbol{(}$\bfseries{\itshape\alph*}$\boldsymbol{)}$]
\item $Y^{\times N}\eqdef \prod^N Y$, resp.~$\mbfY\eqdef Y^{\times\infty}$, always endowed with the product topology;
and \item~\
\end{enumerate*}
\begin{align}\label{eq:TildeSpace}
\andi{Y^{\times n}}\eqdef \tset{\seq{y_i}_i^n: y_i\neq y_j \text{ for } i\neq j}\comma
\end{align}
resp.~$\andi{\mbfY}$, defined analogously, always endowed with the trace topology of~$Y^{\times n}$, resp.~$\mbfY$.
Additionally, let~$\Proj{}{n}\colon \mbfY\rar Y^{\times n}$ be defined by~$\Proj{}{n}\colon \mbfy\eqdef\seq{y_i}_i^\infty\mapsto\seq{y_i}_i^n$.

If~$G\acts Y$, then $G\acts Y^{\times n}$ and~$G\acts \andi{Y^{\times n}}$ coordinate-wise.
We say that~$G\acts Y$ is
\begin{enumerate*}[label=$\boldsymbol{(}$\bfseries{\itshape\alph*}$\boldsymbol{)}$]
\item \emph{transitive} if for every~$y_1,y_2\in Y$ there exists~$g\in G$ such that~$g.y_1=y_2$;
\item \emph{$n$-transitive} if~$G\acts \andi{Y^{\times i}}$ is transitive for every~$i\leq n$;
\item \emph{finitely transitive} if~$G\acts \andi{Y^{\times n}}$ is transitive for every finite~$n$;
\item \emph{$\sigma$-transitive} if~$G\acts \andi{\mbfY}$ is transitive.
\end{enumerate*}
Finally, for~$p\in[1,\infty]$, we denote by~$L^p_G(Y,\mssn)$ the space of classes $u\in L^p(Y,\mssn)$ such that~$u$ has a $G$-invariant $\mssn$-representative.

\begin{prop}\label{p:MeasGInvariant} Let~$G$ be a group acting on~$Y$. Then, a $G$-invariant subset~$A\subset Y$ is Borel measurable if and only if so is~$\pr^G(A)$. Furthermore, for every~$p\in[1,\infty]$, the space~$L^p_G(Y,\mssn)$ is isomorphic to the space~$L^p(\quotient{Y}{G},\pr^G_\pfwd\mssn)$.
\begin{proof}
Let~$\pr\eqdef \pr^G$.
If~$\pr(A)$ is Borel, then so is~$A=\pr^{-1}(\pr(A))$ by measurability (continuity) of~$\pr$. Vice versa, if~$A$ is Borel $G$-invariant, then so is~$A^\complement$. Furthermore,~$(\pr(G.\!\set{y}))^\complement=\pr((G.\!\set{y})^\complement)$, hence, by $G$-invariance of~$A$,~$A^\complement$,
\begin{align*}
\pr(A)^\complement=&\ \pr(G.A)^\complement=\pr\tonde{\cup_{y\in A} G.\set{y}}^\complement=\cap_{y\in A} \pr(G.\set{y})^\complement
\\
=&\ \cap_{y\in A}\pr((G.\set{y})^\complement)=\pr\tonde{\tonde{\cup_{y\in A} G.\set{y}}^\complement}=\pr((G.A)^\complement)=\pr(A^\complement) \fstop
\end{align*}
By continuity of~$\pr$, both~$\pr(A)$ and~$\pr(A)^\complement$ are analytic, and therefore Borel by \cite[Cor.~6.6.10]{Bog07}. The second assertion is a straightforward consequence.
\end{proof}
\end{prop}

\paragraph{Group actions on~$\msP$}
Recall that~$\mfG\eqdef \Diff^\infty_+(\man)$ and let~$\mfF^w\eqdef \seq{\uppsi^{w,t}}_{t\in \R}$ be the one-para\-meter subgroup of~$\mfG$ generated by~$w\in \Vect^\infty$, and~$\mfI\eqdef \Iso(\man,\Bo{\mssg})$ be the group of bijective bi-measurable transformations of~$(\man,\Bo{\mssg})$. The natural action~$G\acts \man$ of any~$G\subset \mfI$ lifts to an action on~$\msP$ as in~\eqref{eq:Intro:GAct}, denoted by~$G_\pfwd$.
\begin{prop}\label{p:TransAct} Assume~$d\geq 1$. Then,
\begin{enumerate*}[label=$\boldsymbol{(}$\bfseries{\itshape\roman*}$\boldsymbol{)}$]
\item\label{i:p:TransAct1}~$\mfI\acts \man$ is $\sigma$-transitive.
Assume~$d\geq 2$. Then,
\item\label{i:p:TransAct2} for every~$n\in \N$ and every~$\mbfx,\mbfx'\in \andi{\man^{\times n}}$ there exists~$w=w_n\in\Vect^\infty$ such that~
\end{enumerate*}
\begin{align}\label{eq:p:TransAct0}
\mbfx'=(\fl^{w,1})^{\times n}(\mbfx)\eqdef\tseq{\fl^{w,1}(x_1),\dotsc, \fl^{w,1}(x_n)} \semicolon
\end{align}
and
\begin{enumerate*}[resume*]
\item\label{i:p:TransAct3} $\mfG\acts \man$ is finitely transitive.
\end{enumerate*}

\begin{proof} \iref{i:p:TransAct1} Let~$\mbfx,\mbfx'\in \andi{\mbfM}$. The map~$g\colon \man\rar \man$ defined by~$g(x_i)\eqdef x_i'$ for all~$i\in \N$ and~$g(x)\eqdef x$ if~$x\neq x_i$ for all~$i\in \N$ is bijective since~$x_i\neq x_j$ and~$x_i'\neq x_j'$ for every~$i\neq j$ and it is straightforwardly bi-measurable.
\iref{i:p:TransAct3} is an immediate consequence of~\iref{i:p:TransAct2}.
Since~$\mbfx,\mbfx'\in\andi{\man^{\times n}}$ there exist~$\eps>0$ and smooth arcs~$\gamma^i_\emparg\colon I\rar \man$ satisfying
\begin{enumerate*}[label=$\boldsymbol{(}$\bfseries{\itshape\alph*}$\boldsymbol{)}$]
\item $\gamma^i_0 = x_i$ and~$\gamma^i_1=x_i'$ for all~$i\leq n$;
\item $\gamma^i_\emparg$ is a simple open arc and $B_\eps(\im \gamma^i_\emparg)$ is contractible for all~$i\leq n$; and 
\item\label{i:p:TransAct:3} $B_\eps(\im \gamma^i_\emparg)\cap B_\eps(\im \gamma^j_\emparg)=\emp$ for~$i\neq j\leq n$.
\end{enumerate*}
For each~$i\leq n$ define a vector field~$w^i$ on~$\im\gamma^i_\emparg$ by~$w^i_{\gamma_t^i}\eqdef \dot\gamma_t^i$ for every~$t\in I$, which is well-posed by~\iref{i:p:TransAct:3}. 
By standard techniques involving partitions of unity, each~$w^i$ may be extended to a (non-relabeled) globally defined smooth vector field vanishing outside~$B_\eps(\im \gamma^i_\emparg)$.
Let~$w\eqdef \sum_i^n w^i$. By construction, one has~$\fl^{w,t}(x_i)=\gamma^i_t$ for every~$i\leq n$ and~$t\in I$, thus~\eqref{eq:p:TransAct0} holds.
\end{proof}
\end{prop}

As a consequence of Propositions~\ref{p:MeasGInvariant},~\ref{p:SuppDF}\iref{i:p:SuppDF2} and~\ref{p:TransAct}\iref{i:p:TransAct1} we have

\begin{cor}\label{c:HomeoQuot}
The Borel spaces~$\ttonde{\quotient{(\Ppaiso,\T_\mrma)}{\mfI_\pfwd},\Proj{\mfI_\pfwd}{\pfwd}\DF_{\beta\n\mssm}}$ and~$(\Tso,\PD_\beta)$ are homeomorphic, isomorphic measure spaces.
\end{cor}

We say that~$A\subset \Ppaiso$ is \emph{$\mfF^{n,w}_\pfwd$-invariant} if
\begin{align*}
\eta=\sum_is_i\delta_{x_i}\in A\comma\qquad \orb_{\mfF^{n,w}}(\eta)\eqdef \bigcup_{t\in \R}\set{\sum_{i\leq n} s_i \delta_{\fl^{w,t}(x_i)}+\sum_{i>n} s_i\delta_{x_i}}\subset A\fstop
\end{align*}

By definition, $\mfF^{\infty,w}_\pfwd$-invariance coincides with~$\mfF^w_\pfwd$-invariance. For~$N\in \overline{\N}$ we say that~$A$ is~$\mfF^N_\pfwd$-invariant if it is $\mfF^{N,w}_\pfwd$-invariant for each~$w\in \Vect^\infty$.
Consistently, for any~$\eta\in \Ppaiso$ we set
\begin{align*}
\orb_{\mfF^N}(\eta)\eqdef \bigcup_{w\in \Vect^\infty} \orb_{\mfF^{N,w}}(\eta) \fstop
\end{align*}
Since the natural action of~$\Diff^\infty(\man)$ on~$\man$ is finitely transitive but not $\sigma$-transitive, $\mfF^n_\pfwd$-invariance and $\mfF^\infty_\pfwd$-invariance are non-comparable notions, as we now show.

\begin{prop}\label{p:GInvar}
Let~$A\subset \Ppaiso$. Then,
\begin{enumerate*}[label=$\boldsymbol{(}$\bfseries{\itshape\roman*}$\boldsymbol{)}$]
\item\label{i:p:GInvar1} $A$ is $\mfI_\pfwd$-invariant if and only if it is $\mfF^n_\pfwd$-invari\-ant for all~$n\in \N$;

\item\label{i:p:GInvar3} if~$A$ is $\mfF^n_\pfwd$-invariant for some~$n\in \N$, then it is also $\mfF^k_\pfwd$-invariant for all~$k\leq n\in\N$;

\item\label{i:p:GInvar2} there exists $A\subset \Ppaisofs$ such that~$A$ is $\mfF^\infty_\pfwd$-invariant but not $\mfI_\pfwd$-invariant.
\end{enumerate*}

\begin{proof}
The forward implication in~\iref{i:p:GInvar1} is straightforward, cf.\ the proof of Prop.~\ref{p:TransAct}\iref{i:p:TransAct1}.
For the reverse implication let~$A$ be $\mfF^n_\pfwd$-invariant for all~$n\in \N$. Set~$\Proj{\times}{n}\eqdef \Proj{}{n}\circ\boldPhi^{-1}$. By Prop.~\ref{p:TransAct}\iref{i:p:TransAct2}
\begin{align*}
\Proj{\times}{n}(A)=&\ \Proj{\times}{n}(\mfF^n.A)=\Proj{}{\Sigma^n}\circ \boldPhi^{-1}(A) \times \andi{\man^{\times n}}=\Proj{\times}{n}(\mfI.A) \fstop
\end{align*}

As a consequence,
\begin{align*}
A=\bigcap_n (\Proj{\times}{n})^{-1}\ttonde{\Proj{\times}{n}(A)}=\bigcap_n (\Proj{\times}{n})^{-1}\ttonde{\Proj{\times}{n}(\mfI.A)}=\mfI.A\comma
\end{align*}
that is,~$A$ is~$\mfI$-invariant.
\iref{i:p:GInvar3} is straightforward since~$\mfF^k.A\subset \mfF^n.A$ for all~$A\subset \Ppaiso$ and all~$k\leq n$.

In order to show~\iref{i:p:GInvar2} let~$\eta=\sum_i s_i\delta_{x_i}\in \Ppaisofs$ be arbitrary. By~\iref{i:p:GInvar1}~$\cup_n\orb_{\mfF_n}(\eta)=\orb_{\mfI}(\eta)$. Let~$\mfH\eqdef\Homeo(\man)$.
We show the stronger statement that 
\begin{align*}
\orb_\mfI(\eta)\supsetneq \orb_{\mfH}(\eta)\supset \orb_{\mfF^\infty}(\eta)\fstop
\end{align*}
Indeed, for fixed~$n\in \N$ let~$y\neq x_n$. On the one hand, the measure~$\tilde\eta\eqdef s_n \delta_y+\sum_{i:i\neq n} s_i \delta_{x_i}$ satisfies~$\tilde\eta\in \orb_{\mfI}(\eta)$ and~$\tilde\eta\neq\eta$, cf. the proof of Prop.~\ref{p:TransAct}\iref{i:p:TransAct1}.
On the other hand, argue by contradiction that there exists~$h\in \mfH$ such that~$h_\pfwd \eta=\tilde\eta$. Since for~$i\neq j$ one has~$s_i\neq s_j$ by definition of~$\Ppaisofs$, then~$h(x_i)=x_i$ for all~$i\neq n$. Again since~$\eta\in \Ppaisofs$, the set~$\set{x_i}_{i:i\neq n}$ is dense in~$\man$, hence, by continuity of~$h$ it must be~$h=\id_\man$, and~$\tilde{\eta}=\eta$, a contradiction.
\end{proof}
\end{prop}

\subsubsection{Pathwise properties}
Let~$(Y,\Bo{\T},\mssn)$ be as in \S\ref{ss:Groups} and recall the following definition.

\begin{defs}[$E$-invariance]\label{d:EInvariance}
Let~$(E,\dom{E})$ be a conservative Dirichlet form on $L^2(Y,\mssn)$. A set~$A\in\Bo{\T}$ is termed $E$-invariant if
\begin{align}\label{eq:Invar}
\car_A u\in\dom{E}\qquad \text{and} \qquad E(u)=E(\car_A u)+E(\car_{A^\complement} u) \comma \qquad u,v\in\dom{E} \fstop
\end{align}

A Borel set~$A$ is $E$-invariant if and only if so is~$A^\complement$. Since~$\car\in\dom{E}$ and~$E(\car)=0$, choosing~$u=\car$ in~\eqref{eq:Invar} yields~$E(\car_A)=0$ for every $E$-invariant~$A$.
Finally, recall that if~$(E,\dom{E})$ is additionally (quasi-)regular with properly associated Markov diffusion process~$M_\bullet$, then~$A$ is $E$-invariant iff it is~$M_\bullet$-invariant.
\end{defs}

In the next result we collect the pathwise properties of~$\Eta$.

\begin{thm}\label{t:PropertiesEta}
Assume~$d\geq 2$ and let~$\Eta_\bullet$ be defined as in~\eqref{eq:Eta}. Then,
\begin{enumerate*}[label=$\boldsymbol{(}$\bfseries{\itshape\roman*}$\boldsymbol{)}$]
\item\label{i:t:PropertiesEta:1} $\Eta_\bullet$ satisfies~\eqref{eq:Intro:Eta};

\item\label{i:t:PropertiesEta:2} $\Eta_\bullet$ is not irreducible: a measurable set~$A\subset \msP$ is $\Eta_\bullet$-invariant if and only if it is~$(\mfI_\pfwd,\DF_\mssm)$-in\-vari\-ant;

\item\label{i:t:PropertiesEta:3} $\Eta_\bullet$ is not ergodic;

\item\label{i:t:PropertiesEta:4}$\Eta_\bullet$ has a (non-relabeled) distinguished extension to all starting points in~$\Ppa$ and satisfying~\eqref{eq:ass:A};

\item\label{i:t:PropertiesEta:5} $\Eta_\bullet$ has $\T_\mrma$-continuous sample paths;

\item\label{i:t:PropertiesEta:6}$\Eta_\bullet$ is a solution to the martingale problem of~$(\mcE,\dom{\mcE})$ in the sense of~\cite[Dfn.~3.1]{AlbRoe95}. In particular, if the initial distribution of~$\Eta_0$ satisfies~$\law(\Eta_0)\ll\DF_\mssm$, then, for each~$u\in\hTF{2}{0}$, the process
\end{enumerate*}
\begin{align*}
M^u_t\eqdef u(\Eta_t)-u(\Eta_0)-\int_0^t \mbfL u(\Eta_s) \diff s
\end{align*}
is a martingale with quadratic variation process
\begin{align*}
\qvar{M^u}_t=\int_0^t \boldGamma(u)(\Eta_s) \diff s \fstop
\end{align*}

\begin{proof}
Since~$\boldPhi$ is bijective between an $\widehat \dirE$-coexceptional and an~$\mcE$-coexceptional set, Equation~\eqref{eq:Intro:Eta} is satisfied as a consequence of Theorem~\ref{t:QuasiHomeo}.

\paragraph{Invariant sets} Assume first that~$A\subset \msP$ is~$(\mfI_\pfwd,\DF_\mssm)$-invariant. Without loss of generality,~$A\subset \Ppaiso$, since~$\DF_\mssm\Ppaiso=1$ and~$\Ppaiso$ is $\mfI_\pfwd$-, hence~$(\mfI_\pfwd,\DF_\mssm)$-, invariant. By a straightforward density argument and Corollary~\ref{c:HomeoQuot}, cf.~\eqref{eq:TestF},
\begin{align*}
L^2(\Tso,\PD_\beta)\cong L^2\ttonde{\quotient{\Ppaiso}{\mfI_\pfwd}, \pr^{\mfI_\pfwd}_\pfwd \DF_\mssm}\cong L^2_{\mfI_\pfwd}(\msP,\DF_\mssm)=\cl_{L^2(\DF_\mssm)}\hTF{1}{-,0}\fstop
\end{align*}

Since~$\grad\equiv \zero$ on~$\hTF{1}{-,0}$, one has~$\mcE\equiv \zero$ on~$\cl_{\mcE^{1/2}_1}(\hTF{1}{-,0})=\cl_{L^2(\DF_\mssm)}\hTF{1}{-,0}$.
By strong locality,~$\boldGamma$ satisfies the Leibniz rule, hence
\begin{align}\label{eq:Invar2}
\car_A u\in\dom{\mcE} \quad \text{and}\quad \mcE(\car_A u,v)=\mcE(\car_A u,\car_A v)=\mcE(u,\car_A v)\comma \quad u,v\in \dom{\mcE}\comma
\end{align}
as soon as~$\car_A\in L^2_{\mfI_\pfwd}(\msP,\DF_\mssm)$ or, equivalently,~$A$ is~$(\mfI_\pfwd,\DF_\mssm)$-invariant. Thus,~\eqref{eq:Invar} follows by~\eqref{eq:Invar2} since~$\boldGamma(\car_A,u)=0$ for every~$u\in\dom{\mcE}$.
Viceversa, assume that~$A$ is $\mcE$-invariant. If~$\DF_\mssm A=0$, resp.~$1$, then~$A$ is $(\mfI_\pfwd,\DF_\mssm)$-invariant since~$\emp$, resp.~$\andihM$, is.
Assume then~$\DF_\mssm A\in (0,1)$. Without loss of generality,~$A\subset \Ppaiso$, since~$\Ppaiso$ is $\mcE$-coexceptional.
Thus,~$\boldPhi^{-1}$ is well-defined on~$A$ and~$B\eqdef\boldPhi^{-1}(A)\subset \andihM$ is $\widehat\dirE$-invariant by Theorem~\ref{t:QuasiHomeo}.
Since~$(\dirE^\mbfs,\dom{\dirE^\mbfs})$ is ergodic for every~$\mbfs\in \To$ by the discussion in~\cite[\S1]{BenSaC97}, it follows that the only $\widehat\dirE$-invariant sets are of the form~$C\times \prod_i ^\infty U_i$ where~$C\in\Bo{\mrmu}(\Tso)$ is any measurable set and~$U_i$ satisfies either~$U_i=\man_i$ or~$U_i=\emp$ for all~$i\in \N$.
Since~$\widehat\volm_\beta B\in (0,1)$, it must be~$U_i=\man$ for every~$i\in \N$ and~$\PD_\beta C\in (0,1)$, that is~$B=C\times \mbfM$. The~$(\mfI_\pfwd,\DF_\mssm)$-invariance of~$A$ follows by the~$\mfI\acts \mbfM$-invariance of~$\mbfM$ since
\begin{align*}
g_\pfwd \boldPhi(\mbfs,\mbfx)=\boldPhi(\mbfs,g^{\times \infty}(\mbfx))\comma\qquad g\in \mfI \fstop
\end{align*}

\paragraph{Lack of ergodicity} Let~$\eps\in (0,1)$. Since~$\PD_\beta$ is diffuse, by the main result in~\cite{Sie22} there exists~$A''\in\Bo{\mrmu}(\Tso)$ with~$\PD_\beta A''=\eps$.
Let~$A'$ be the corresponding subset of~$\quotient{\Ppaiso}{\mfJ_\pfwd}$ via the homeomorphism of Corollary~\ref{c:HomeoQuot}.
Then,~$A\eqdef (\pr^{\mfI_{\pfwd}})^{-1}(A')$ is $\mfI_\pfwd$-invariant by definition, $\mcB_\mrmn(\msP)$-measurable by Proposition~\ref{p:MeasGInvariant}, satisfying~$\DF_\mssm A=\eps$ by Corollary~\ref{c:HomeoQuot} and $\mcE$-invariant by the previous step.
As a consequence,~$(\mcE,\dom{\mcE})$ is \emph{not} ergodic.

\paragraph{Continuity of paths and extension} For every~$\eta_0\in \Ppaiso$ the path~$t\mapsto \Eta_t^{\eta_0}$ is $\T_\mrma$-continuous as a consequence of Corollary~\ref{c:QuasiHomeo} and the standard theory of Dirichlet forms.
Consistently with the definition of~$\BM^\mbfs_\bullet$ for~$\mbfs\in\Tso$, for~$\mbfs\in\To\setminus \Tso$ we set
\begin{align*}
\BM^{\mbfs;\mbfx_0}_t(\omega)\eqdef\ttonde{x^1_{t/s_1}, x^2_{t/s_2},\dotsc}\comma
\end{align*}
where~$\tseq{x^i_t}_{t\geq 0}$ are independent Brownian motions on~$\man$ and, conventionally,~$x^i_{t/s_i}=x^i_0$ for all~$t\geq 0$ whenever~$s_i=0$. Then, letting
\begin{align*}
\widehat\BM^{\mbfs,\mbfx_0}_\bullet\eqdef \BM^{\mbfs;\mbfx_0}_\bullet\comma\qquad \textrm{and} \qquad \eta_0\eqdef \boldPhi(\mbfs,\mbfx_0)\comma \quad \Eta_\bullet^{\eta_0}\eqdef \boldPhi\circ \widehat\BM^{\mbfs,\mbfx_0}_\bullet
\end{align*}
yields the desired extension to all starting points in~$\boldPhi(\To\times\andi{\mbfM})= \Ppa$ satisfying~\eqref{eq:ass:A}.

\paragraph{Martingale problem} Since~$(\mcE,\dom{\mcE})$ on~$L^2(\Ppaiso,\T_\mrma,\DF_\mssm)$ is not $\T_\mrma$-regular because of Remark~\ref{r:LocCompact}, a proof of~\iref{i:t:PropertiesEta:6} is not entirely standard. It follows by checking the assumptions in~\cite[Thm.~3.4]{AlbRoe95}. The $\T_\mrma$-quasi-regularity of~$(\mcE,\dom{\mcE})$ is shown in Corollary~\ref{c:QuasiHomeo}, thus it suffices to check~\cite[(A), p.~517]{AlbRoe95} for the choice~$D=\hTF{\infty}{0}$: The probability measure~$\DF_\mssm$ is fully supported by Proposition~\ref{p:SuppDF}\iref{i:p:SuppDF1}. The $\T_\mrma$-continuity of~$u\in\hTF{\infty}{0}$ is noted in Remark~\ref{r:MeasContCyl}\iref{i:r:MeasContCyl:3}.
\end{proof}
\end{thm}

As a consequence of Theorem~\ref{t:PropertiesEta}\iref{i:t:PropertiesEta:5}--\iref{i:t:PropertiesEta:6}, the process~$\Eta_\bullet$ defined as in~\eqref{eq:Eta} is, equivalently, a solution of the following martingale problem.

\begin{cor}[Martingale problem]\label{c:MartingaleP}
For every~$\hat f\in\mcC^2(\hman_0)$ as in Definition~\ref{d:CylFunc}, the process
\begin{align*}
M^{\hat f}_t\eqdef \hat f^\trid(\Eta_t)-\hat f^\trid(\Eta_0)-\int_0^t \mbfB_0\tquadre{\nabla\hat f}(\Eta_s) \diff s
\end{align*}
is a continuous martingale with quadratic variation process
\begin{align*}
\qvar{M^{\hat f}}_t=\int_0^t \Gamma(\hat f)^\trid (\Eta_s) \diff s \fstop
\end{align*}
\end{cor}

\begin{rem}[Distinguished invariant measures]
By Theorem~\ref{t:QuasiHomeo} and~\cite[Cor.~3.18]{LzDS20}, it follows from the extremality of $\widehat\BM_\bullet$-ergodic measures that~$\mbfQ\in\msP(\msP)$ is $\Eta_\bullet$-ergodic if and only if it is of the form~$\mbfQ_\mbfs\eqdef \boldPhi_\pfwd(\delta_\mbfs\otimes\n\volm)$ for some~$\mbfs\in \Tso$. It is straightforward that every such measure satisfies
\begin{align}\label{eq:ass:B}
\mbfE_{\mbfQ_\mbfs}\, \mbfE_{P_\emparg}[\Eta^\emparg_t A]=\n\mssm A\comma \qquad A\in\Bo{\mssg}\comma t\geq 0\fstop
\end{align}

More generally,~\eqref{eq:ass:B} holds for any $\Eta_\bullet$-invariant measure~$\mbfQ$, since
\begin{align*}
\mbfQ\in \overline{\conv}\set{\mbfQ_\mbfs : \mbfs\in\Tso}\fstop
\end{align*}
\end{rem}

\begin{rem}[Comparison with the MMAF]
As already noted in~\cite[\S1]{KonvRe18} for the Modified Massive Arratia Flow in the case~$d=1$, also in the case~$d\geq 2$ we do not expect the family~$\tseq{\widehat\mssW^n_\bullet}_n$ to be a compatible family of Feller semigroups in the sense of Le Jan--Raimond~\cite[Dfn.~1.1]{LeJRai04}; thus the process~$\Eta_\bullet$ is not induced by a stochastic flow.
\end{rem}

\subsubsection{The Rademacher property}
We show a weak form of the Rademacher property for $(\mcE,\dom{\mcE})$. We refer to the results in~\cite{LzDS19b}, from which a proof is adapted.
\begin{prop}\label{p:RademacherE}
Assume~$d\geq 2$. If~$u\in \Lip(\msP_2)$, then~$u\in \dom{\mcE}$ and~$\boldGamma(u)\leq \Lip[u]^2$ $\DF_\mssm$-a.e..
\begin{proof}
By~\cite[Thm.~3.9]{LzDS19a},~$\Psi^{w,t}_\pfwd \DF_\mssm=\DF_{\uppsi^{w,t}_\pfwd \mssm}$. 
Let~$\mcF$ be the set of all bounded measurable functions~$u$ on~$\msP$ for which there exists a measurable section~$\mbfD u$ of~$T^\Der\msP_2$ such that
\begin{align*}
\int_\msP \scalar{\mbfD u(\eta)}{\mbfD u(\eta)}_{\Vect_\eta} \diff\DF_\mssm(\eta) <\infty
\end{align*}
and, for all~$s\in \R$ and~$w\in\Vect^\infty$,
\begin{align*}
\frac{u\circ \Psi^{w,t}-u}{t}\xrightarrow{t\rar 0} \scalar{\grad u}{w}_{\Vect_\emparg} \quad \textrm{in}\quad L^2(\msP,\DF_{\uppsi^{w,s}_\pfwd\mssm}) \fstop
\end{align*}

By Lemma~\ref{l:DirDer},~$\hTF{\infty}{0}\subset \mcF$. Since the generator of~$(\mcE,\dom{\mcE})$ is essentially self-adjoint on~$\hTF{\infty}{0}$ by Proposition~\ref{p:ESA}, the form~$(\mcE,\dom{\mcE})$ coincides with the form~$(\mcE,\msF)$ defined in~\cite{LzDS19b} with~$\mbbP=\DF_\mssm$.
Moreover,~$(\mcE,\dom{\mcE})$ coincides with the closure of~$(\mcE,\TF{\infty})$ by Lemma~\ref{l:TFDensity}.
Thus, for~$\mbbP=\DF_\mssm$, the forms~$(\mcE,\msF_0)$,~$(\mcE,\msF_{\cont})$ and~$(\mcE,\msF)$ defined in~\cite[Thm.~2.10(1)]{LzDS19b} all coincide with $(\mcE,\dom{\mcE})$ by~\cite[Rmk.~2.13]{LzDS19b}.

As already noted in~\cite[\S5.4]{LzDS19b},~$\DF_\mssm$ satisfies Assumptions~2.6$(i)$-$(ii)$ there. Since we have the closability of~$(\mcE,\TF{\infty})$ independently of Assumption~2.6$(iii)$ there, the strategy of~\cite{LzDS19b} applies verbatim, except for~\cite[Lem.~4.8, Prop.~4.9]{LzDS19b}. We show suitable replacements respectively in Lemma~\ref{l:Rademacher} and Proposition~\ref{p:Rademacher} below.
\end{proof}
\end{prop}

Next, we discuss the short-time asymptotics of~$\Eta_\bullet$.
For~$A_1,A_2\in\mcB_\mrmn(\msP)$ of positive~$\DF_\mssm$-measure set
\begin{align*}
\mssd_{W_2}(A_1,A_2)\eqdef \DF_\mssm\textrm{-}\essinf_{\mu_i\in A_i} W_2(\mu_1,\mu_2)\comma
\\
p_t(A_1,A_2)\eqdef \int_{A_1} \int_{A_2} p_t(\mu_1,\diff\mu_2) \diff\DF_\mssm(\mu_1)
\end{align*}
and let
\begin{align}
\mssd_{\mcE}(\mu,\nu)\eqdef \sup\set{u(\mu)-u(\nu) : u\in\dom{\boldGamma}\cap\mcC(\msP)\comma \boldGamma(u)\leq 1\quad  \DF_\mssm\text{-a.e.}}
\end{align}
be the intrinsic distance of~$(\mcE,\dom{\mcE})$. Then,
\begin{cor}[Varadhan short-time asymptotics lower bound]\label{c:Varadhan}
It holds that 
\begin{align}\label{eq:c:Varadhan:0}
W_2(\emparg,\mu)\leq \mssd_\mcE(\emparg,\mu)\comma \qquad \mu\in\msP_2\comma
\end{align}
and
\begin{align}\label{eq:Varadhan}
\lim_{t\downarrow 0} t\log p_t(A_1,A_2)\leq-\tfrac{1}{2}\mssd_{W_2}(A_1,A_2)^2 \fstop
\end{align}
\begin{proof}
Equation~\eqref{eq:c:Varadhan:0} is an immediate consequence of Proposition~\ref{p:RademacherE}. For~$W_2$ is $\T_\mrmn$-continuous,~$W_2(\emparg, A)=W_2\ttonde{\emparg, \overline A}$, where~$\overline A$ denotes the $\T_\mrmn$-closure of~$A$.
By compactness of~$(\msP,\T_\mrmn)$, the set~$\overline A$ too is $\T_\mrmn$ compact.
In particular, there exists a sequence~$\seq{\mu_n}_n$ so that~$W_2\ttonde{\emparg, \overline A}=\nlim \min_{m\leq n} W_2(\emparg, \mu_m)$.
Arguing similarly to~\cite[Lem.~3.4]{LzDS20}, it follows that~$W_2(\emparg, A)\in\dom{\boldGamma}$ and~$\boldGamma \ttonde{W_2(\emparg, A)}\leq 1$ $\DF_\mssm$-a.e.; furthermore~$W_2(\emparg, A)=0$ $\DF_\mssm$-a.e.\ on~$A$.
The Varadhan-type estimate~\eqref{eq:Varadhan} now follows combining the general results~\cite[Thms.~1.1,~1.2]{HinRam03}.
\end{proof}
\end{cor}

Although most of the previous results only hold when~$d\geq 2$, it is possible to construct a regular strongly local Dirichlet form on~$\msP$ also in the case when~$d=1$, i.e. when~$\man=\mbbS^1$, as a restriction of the form~$(\mcE,\dom{\mcE})$.

\begin{prop}[Reduced form]
The form~$(\mcE,\TF{1})$ on~$L^2(\msP_2(\mbbS^1),\DF_\mssm)$ is closable. Its closure~$(\mcE^\red,\dom{\mcE^\red})$ is a $\T_\mrmn$-regular strongly local Dirichlet form.
\begin{proof}
By Lemma~\ref{l:Domains} we have~$\TF{1}\subset \dom{\mcE}$, hence the statement is well-posed. The closability of~$(\mcE,\TF{1})$ follows since~$(\mcE,\dom{\mcE})$ is closed. The Markov property and strong locality are inherited from~$(\mcE,\dom{\mcE})$. The density of~$\TF{1}$ in~$\Cb(\msP_2)$ holds as in the proof of Corollary~\ref{c:Regularity}. 
\end{proof}
\end{prop}

If~$d\geq 2$ we have~$(\mcE^\red,\dom{\mcE^\red})=(\mcE,\dom{\mcE})$ by Lemma~\ref{l:TFDensity}. However, if~$d=1$, the proof of Lemma~\ref{l:TFDensity} fails and the form~$(\mcE,\dom{\mcE})$ might be \emph{not} $\T_\mrmn$-regular.

\section{Appendix}

\subsection{Measurability properties} 

\begin{lem}\label{l:Meas}
For~$A\in\Bo{\mssg}$ and~$\mu\in\msP$ define the evaluation map $\ev_A\colon \mu\mapsto \mu A$.
Then,
\begin{enumerate*}[label=$\boldsymbol{(}$\bfseries{\itshape\roman*}$\boldsymbol{)}$]
\item\label{i:l:Meas1} the map~$\ev_A$ is $\Bo{\mrmn}(\msP)$-measurable;
\item\label{i:l:Meas2} $\Bo{\mrmn}(\msP)$ is generated by the maps~$\set{\ev_A}_{A\in \Bo{\mssg}}$;
\item\label{i:l:Meas3} the convex combination~$(r,\mu,\nu)\mapsto (1-r)\mu+r\nu$ of parameter~$r$ is jointly~$\mcB(I)\otimes \mcB_\mrmn(\msP)^{\otimes 2}$\-mea\-surable.
\end{enumerate*}

\begin{proof}
\iref{i:l:Meas1}~is a consequence of~\iref{i:l:Meas2} which is in turn~\cite[Thm.~1.5]{Kal17}. Since $\Mb(\man)\supset\msP$, endowed with the weak* topology, is a measurable vector space,~\iref{i:l:Meas3} follows by~\cite[Prop.~I.2.3, p.~16]{VakTarCho87}.
\end{proof}
\end{lem}

\begin{lem}\label{l:JointMeas}
The map~$\ev\colon (\mu,x)\mapsto \mu_x\eqdef \mu\!\set{x}$ is $\Bo{\mrmn}(\msP)\otimes \Bo{\mssg}$-measurable.
\begin{proof}
Let~$\mssh_t^*\colon \msP\rar \msP$ be the heat flow on measures. Then,
\begin{enumerate*}[label=$\boldsymbol{(}$\bfseries{\itshape\alph*}$\boldsymbol{)}$]
\item $\mssh_t^*\colon \msP\rar \msP$ is narrowly continuous for every~$t>0$; 
\item $t\mapsto \mu_t\eqdef \mssh_t^*\mu$ is narrowly continuous for every~$\mu\in\msP$;
\item $\mu_t\ll \mssm$ for every~$t>0$ and every~$\mu\in\msP$.
\end{enumerate*}
For each~$\eps>0$,~$t>0$ and~$x\in X$ the map~$\mu\mapsto \mu_t B_\eps(x)$ is measurable, since it is the composition of the continuous map~$\mu\mapsto \mu_t$ with the measurable map~$\ev_{B_\eps(x)}$ (see Lem.~\ref{l:Meas}). Moreover, it is readily seen  by Dominated Convergence that for each~$\eps>0$,~$t>0$ and~$\mu\in \msP$ the map~$x\mapsto \mu_tB_\eps(x)$ is continuous, since~$\diff \mu_t(y)=f_t(y)\diff\mssm(y)$ for some~$f_t\in L^1(\man, \mssm)$. That is,~$\ev_{\eps,t}\colon (\mu,x)\mapsto \mu_tB_\eps(x)$ is a Carath\'eodory map between Polish spaces, hence it is jointly measurable. Since the pointwise limit of (jointly-)measurable maps is  (jointly-)measurable, it suffices to show the existence of~$\lim_{\eps\downarrow 0} \lim_{t\downarrow 0} \ev_{\eps,t}=\ev$. To this end,
\begin{align*}
\liminf_{\eps\downarrow 0}\liminf_{t\downarrow 0} \mu_t B_\eps(x)\geq& \liminf_{\eps\downarrow 0} \mu B_\eps(x)=\mu_x \comma\\
\limsup_{\eps\downarrow 0}\limsup_{t\downarrow 0} \mu_t B_\eps(x)\leq& \limsup_{\eps\downarrow 0}\limsup_{t\downarrow 0} \mu_t \ttonde{\overline{B_\eps(x)}} \leq \limsup_{\eps\downarrow 0} \mu\, \ttonde{\overline{B_\eps(x)}}\\
\leq& \limsup_{\eps\downarrow 0} \mu B_{2\eps}(x)=\mu_x
\end{align*}
by the Portmanteau Theorem and the outer regularity of~$\mu$.
\end{proof}
\end{lem}

\begin{prop}\label{p:Dynkin}
Let~$\Omega$ be any non-empty set, $\mcA$ be a multiplicative system of bounded real-valued functions on~$\Omega$. Let~$\mcB$ be the $\sigma$-algebra generated by~$\mcA$ and denote by~$\mcB_b$ the space of bounded $\mcB$-measurable real-valued functions. Then, for any non-negative finite measure~$\mu$ on~$(\Omega,\mcB)$, the system~$\mcA$ is dense in~$L^2(\Omega,\mu)$.

\begin{proof}
Since~$\mu$ is finite and functions in~$\mcA$ are bounded, then~$\mcA\subset L^2(\Omega,\mu)$. Let now~$v\in \mcA^\perp\subset L^2(\Omega,\mu)$ and~$H\subset \mcB_b$ be maximal such that~$\int_\Omega v\, h \diff\mu=0$ for every~$h\in H$. It suffices to show that~$v=0$ $\mu$-a.e.. We show that~$\mcB_b\subset H$, from which the previous assertion readily follows. Observe that~$H$ is a vector space, uniformly closed in~$\R^\Omega$ and closed under monotone convergence of non-negative uniformly bounded sequences by Dominated Convergence. Since~$\mcA\subset H$ is multiplicative,~$\mcB_b\subset H$ by Dynkin's Multiplicative System Theorem~\cite[Thm.~I.2.12.9(i)]{Bog07}.
\end{proof}
\end{prop}

\subsection{Quadratic forms}
Let~$(H,\norm{\emparg}_H)$ be a real \emph{separable} Hilbert space.
\begin{defs}\label{d:QForms} By a \emph{quadratic form}~$(Q,D)$ on~$H$ we shall always mean a symmetric positive semi-definite ---~if not otherwise stated, densely defined~--- bilinear form. To~$(Q,D)$ we associate the non-relabeled functional~$Q\colon H\rar \R\cup\set{+\infty}$ defined by
\begin{align*}
Q(u)\eqdef \begin{cases}
Q(u,u) & \text{if } u\in D
\\
+\infty & \text{otherwise}
\end{cases} \comma \qquad u\in H\fstop
\end{align*}

Additionally, we set for every~$\alpha>0$
\begin{align*}
Q_\alpha(u,v)\eqdef& Q(u,v)+\alpha\scalar{u}{v}_{H} \comma \qquad u,v\in D\comma
\\
Q_\alpha(u)\eqdef& Q(u)+\alpha\norm{u}_{H}^2\comma \qquad\qquad u\in H \fstop
\end{align*}

For~$\alpha>0$, we let~$\dom{Q}_\alpha$ be the completion of~$D$, endowed with the Hilbert norm~$Q_\alpha^{1/2}$.
\end{defs}

The following result is well-known.

\begin{lem}\label{l:WLSC} Let~$(Q, D)$ be a quadratic form on~$H$. The following are equivalent:
\begin{enumerate*}[label=$\boldsymbol{(}$\bfseries{\itshape\alph*}$\boldsymbol{)}$]
\item $(Q,D)$ is closable, say, with closure~$(Q,\dom{Q})$;
\item the canonical inclusion~$\iota\colon D\rar H$ extends to a continuous injection~$\iota_\alpha\colon \dom{Q}_\alpha\rar H$ satisfying~$\norm{\iota_\alpha}\leq \alpha^{-1}$;
\item $Q$ is lower semi-continuous w.r.t. the strong topology of~$H$;
\item $Q$ is lower semi-continuous w.r.t. the weak topology of~$H$.
\end{enumerate*}
\end{lem}

To every closed quadratic form~$(Q,\dom{Q})$ we associate a non-negative self-adjoint operator~$-L$, with domain defined by the equality~$\dom{\sqrt{-L}}=\dom{Q}$, such that~$Q(u,v) = \scalar{-Lu}{v}_{H}$ for all~$u,v\in\dom{Q}$. We denote the associated semigroup by~$T_t\eqdef e^{tL}$,~$t>0$, and the associated resolvent by~$G_\alpha\eqdef (\alpha-L)^{-1}$,~$\alpha> 0$. By Hille--Yosida Theorem, e.g.~\cite[p.~27]{MaRoe92},
\begin{subequations}\label{eq:Hille--Yosida}
\begin{align}
\label{eq:Hille--YosidaG}
Q_\alpha(G_\alpha u,v)=&\scalar{u}{v}_H\comma\qquad v\in\dom{Q}\comma u\in H\comma
\\
\label{eq:Hille--YosidaT}
T_t=&H\text{-}\!\!\lim_{\alpha\rar \infty} e^{t\alpha (\alpha G_\alpha-1)} \fstop
\end{align}
\end{subequations}

\subsection{Generalized Mosco convergence of quadratic forms}\label{ss:KuwaeShioya}
We shall need K.~Kuwae and T.~Shioya's \emph{generalized Mosco convergence}~\cite{KuwShi03}. We start by recalling the simplified setting introduced by A.~Kolesnikov in~\cite[\S2]{Kol06}.

\begin{defs}[Convergences]\label{d:ConvHSpaces}
Let~$\seq{H_n}_n$ and $H$~be Hilbert spaces and set $\mcH\eqdef H\sqcup\sqcup_n H_n$. Let further~$D\subset H$ be a dense subspace and~$\seq{\Phi_n}_n$ be densely defined linear operators
\begin{align}\label{eq:ConvH}
\Phi_n\colon D\rar H_n \fstop
\end{align}

We say that~$H_n$ $\mcH$-\emph{converges} to~$H$ if
\begin{align}\label{eq:ConvHSpaces1}
\nlim \norm{\Phi_n u}_{H_n}=\norm{u}_{H} \comma \qquad u\in D \comma
\end{align}
in which case we say further that a sequence~$\seq{u_n}_n$,~$u_n\in H_n$,
\begin{enumerate*}[label=$\boldsymbol{(}$\bfseries{\itshape\alph*}$\boldsymbol{)}$]
\item $\mcH$-\emph{strongly converges} to~$u\in H$ if there exists a sequence~$\seq{\tilde u_m}_m\subset D$ such that
\end{enumerate*}
\begin{align*}
\mlim\norm{\tilde u_m-u}_{H}=0 \qquad \text{and} \qquad \mlim \nlimsup\norm{\Phi_n\tilde u_m-u_n}_{H_n}=0\semicolon
\end{align*}
\begin{enumerate*}[resume*]
\item $\mcH$-\emph{weakly converges} to~$u\in H$ if, for every sequence~$\seq{v_n}_n$,~$v_n\in H_n$, $\mcH$-strongly converging to~$v\in H$,
\end{enumerate*}
\begin{align*}
\nlim \scalar{u_n}{v_n}_{H_n}=\scalar{u}{v}_{H} \fstop
\end{align*}

Let further~$\seq{B_n}_n$ be a se\-quence of bounded operators~$B_n\in \mcB(H_n)$. We say that~$B_n$ converges $\mcH$-strongly to~$B\in\mcB(H)$ if $B_nu_n$ converges $\mcH$-strongly to~$Bu$ for all sequences~$\seq{u_n}_n$, $u_n\in H_n$, such that~$u_n$ $\mcH$-strongly converges to~$u\in H$.
\end{defs}

\begin{defs}[generalized Mosco convergence]\label{d:Kuwae--Shioya}
Let~$\tseq{(Q_n,\dom{Q_n})}_n$ be a se\-quence of closed quadratic forms, $Q_n$ on~$H_n$, and~$(Q,\dom{Q})$ be a quadratic form on~$H$. We say that~$Q_n$ \emph{Mosco converges} to~$Q$ if the following conditions hold: 
\begin{enumerate*}[label=$\boldsymbol{(}$\bfseries{\itshape\alph*}$\boldsymbol{)}$]
\item $H_n$ $\mcH$-converges to~$H$;
\item\label{i:MoscoB} (\emph{weak $\Gamma$-$\liminf$}) if~$\seq{u_n}_n$,~$u_n\in H_n$, $\mcH$-weakly converges to~$u\in H$, then
\end{enumerate*}
\begin{align}\label{eq:GammaLimInf}
Q(u)\leq \nliminf Q_n(u_n) \semicolon
\end{align}
\begin{enumerate*}[label=$\boldsymbol{(}$\bfseries{\itshape\alph*}$\boldsymbol{)}$]\setcounter{enumi}{2}
\item (\emph{strong $\Gamma$-$\limsup$}) for every~$u\in H$ there exists a sequence~$\seq{u_n}_n$,~$u_n\in H_n$, $\mcH$-strongly convergent to~$u$ and such that
\end{enumerate*}
\begin{align}\label{eq:GammaLimSup}
Q(u)=\nlim Q_n(u_n)
\end{align}

Clearly, in condition~\iref{i:MoscoB} we can additionally assume~$u_n\in\dom{Q_n}$.
\end{defs}

\begin{rem} In all the above definitions, the notion of convergence \emph{does depend} on the family of linear operators~$\seq{\Phi_n}_n$. The latter is however omitted from the notation, for it will be apparent from the context.
\end{rem}

\begin{lem}\label{l:UniEquiv}
Let~$\tseq{(Q_n,\dom{Q_n})}_n$, resp.~$\tseq{(Q_n^\sharp,\dom{Q_n^\sharp})}_n$, be a sequence of closed quad\-ratic forms~$Q_n$ on~$H_n$, resp.~$Q_n^\sharp$ on~$H_n^\sharp$, and~$(Q,\dom{Q})$ be a quadratic form on~$H$.
Further assume that there exist unitary operators~$\seq{U_n}_n$ such that~$Q_n^\sharp\circ U_n=Q_n$. 
Then,~$\tseq{(Q_n,\dom {Q_n})}_n$ Mosco converges to~$(Q,\dom{Q})$ if and only if~$\tseq{(Q_n^\sharp,\dom{Q_n^\sharp})}_n$ Mosco converges to~$(Q,\dom{Q})$.
\begin{proof}
Assume~$\tseq{(Q_n,\dom {Q_n})}_n$ Mosco converges to~$(Q,\dom{Q})$. Let~$D\subset H$ and~$\Phi_n\colon D\rar H_n$ be as in Definition~\ref{d:ConvHSpaces} for~$n\in \N$.
Set~$\Phi_n^\sharp\eqdef U_n\circ\Phi_n\colon D\rar H_n^\sharp$. Then, since~$U_n\colon H_n\rar H_n^\sharp$ is unitary and by~\eqref{eq:ConvHSpaces1},
\begin{align*}
\nlim \tnorm{\Phi_n^\sharp u}_{H_n^\sharp}=\nlim \norm{\Phi_n u}_{H_n}=\norm{u}_H\comma \qquad u\in D\comma
\end{align*}
hence~$H_n^\sharp$ $\mcH$-converges to~$H$.
Analogously, one can show that~$\seq{u_n}_n$,~$u_n\in H_n$, $\mcH$-strongly, resp.\ $\mcH$-weakly, converges to~$u\in H$ if and only if~$\tseq{u_n^\sharp}_n$, $u_n^\sharp\eqdef U_n(u_n)\in H_n^\sharp$, $\mcH$-strongly, resp. $\mcH$-weakly converges to~$u\in H$.
Thus, let~$u\in H$ and~$\seq{u_n}_n$,~$u_n\in H_n$, be as in~\eqref{eq:GammaLimSup} and note that~$\tseq{u_n^\sharp}_n$ defined as above $\mcH$-strongly converges to~$u\in H$ and
\begin{align*}
Q(u)=\nlim Q_n(u_n)= \nlim Q_n^\sharp (U_n u_n) = \nlim Q_n^\sharp(u_n^\sharp)\comma
\end{align*}
which proves the $\Gamma$-$\limsup$ condition~\eqref{eq:GammaLimSup} for~$Q^\sharp_n$.
Finally, let~$\tseq{u_n^\sharp}_n$,~$u_n^\sharp\in H_n^\sharp$, be $\mcH$-weakly converging to~$u\in H$ and set~$u_n\eqdef U_n^{-1}(u_n^\sharp)\in H_n$. Then,~$\seq{u_n}_n$ $\mcH$-weakly converges to~$u\in H$, hence, by assumption on~$Q_n$,
\begin{align*}
Q(u)\leq \nliminf Q_n(u_n)=\nliminf Q_n^\sharp (U_n u_n )= \nliminf Q_n^\sharp(u_n^\sharp)\comma
\end{align*}
which proves the $\Gamma$-$\liminf$ condition~\eqref{eq:GammaLimInf} for~$Q_n^\sharp$ and concludes the proof.
\end{proof}
\end{lem}

\begin{lem}\label{l:ConvHStandard}
Let~$H$ be a Hilbert space and~$\seq{P_n}_n$ be an increasing sequence of orthogonal projectors~$P_n\colon H\rar H$ strongly converging to~$\id_H$. Set~$H_n\eqdef \ran P_n$ and let further~$u\in H$ and~$\seq{u_n}_n$ be a sequence such that~$u_n\in H_n$. Then,
\begin{enumerate*}[label=$\boldsymbol{(}$\bfseries{\itshape\roman*}$\boldsymbol{)}$]
\item\label{i:l:ConvHStandard0} $H_n$ $\mcH$-converges to~$H$;
\item\label{i:l:ConvHStandard1}~$\seq{u_n}_n$ $\mcH$-strongly converges to~$u\in H$ if and only if it strongly converges to~$u$ in~$H$;
\item\label{i:l:ConvHStandard2}~$\seq{u_n}_n$ $\mcH$-weakly converges to~$u\in H$ if and only if it weakly converges to~$u$ in~$H$.
\end{enumerate*}
\begin{proof} 
\iref{i:l:ConvHStandard0} is an immediate consequence of the strong convergence of~$P_n$ to~$\id_H$.

\iref{i:l:ConvHStandard1} Assume~$u_n$ strongly converges to~$u$ and choose~$D=H$,~$\Phi_n\eqdef P_n$ and~$\tilde u_m=u_m$ in Definition~\ref{d:ConvHSpaces}.
By strong convergence of~$P_n$ to~$\id_H$ one has~$\mlim \norm{u_m-u}_H=0$. Furthermore,
\begin{align*}
\mlim \nlimsup \norm{P_n u_m-u_n}_{H_n}=&\mlim \nlimsup \norm{P_n(u_m-u_n)}_{H_n}
\\
\leq& \mlim\nlimsup \norm{P_n} \norm{u_m-u_n}_H=0
\end{align*}
and the conclusion follows.
Viceversa, assume that~$u_n$ $\mcH$-strongly converges to~$u$. Then,
\begin{align*}
\norm{u_n-u}_H\leq& \norm{u_n-P_n \tilde u_m}_H+\norm{P_n \tilde u_m-\tilde u_m}_H+\norm{\tilde u_m-u}_H
\\
=& \norm{u_n-P_n \tilde u_m}_{H_n}+\norm{P_n \tilde u_m-\tilde u_m}_H+\norm{\tilde u_m-u}_H \fstop
\end{align*}

Taking first the limit superior in~$n$ and, subsequently, the limit in~$m$, the above inequality readily yields the conclusion.
A proof of~\iref{i:l:ConvHStandard2} follows similarly to~\iref{i:l:ConvHStandard1} by definition of $\mcH$-weak convergence, and it is therefore omitted.
\end{proof}
\end{lem}

The main result concerning generalized Mosco convergence is the following:
\begin{thm}[Kuwae--Shioya~{\cite[Thm.~2.4]{KuwShi03}}]
Let~$\tseq{(Q_n,\dom{Q_n})}_n$ be a sequence of closed quadratic forms, $Q_n$ on~$H_n$, and~$(Q,\dom{Q})$ be a closed quadratic form on~$H$. Then, the following are equivalent:
\begin{enumerate*}[label=$\boldsymbol{(}$\bfseries{\itshape\alph*}$\boldsymbol{)}$]
\item $Q_n$ Mosco converges to~$Q$;
\item $G_{n,\alpha}$ $\mcH$-strongly converges to~$G_\alpha$ for every~$\alpha>0$;
\item $T_{n,t}$ $\mcH$-strongly converges to~$T_t$ for every~$t>0$.
\end{enumerate*}
\end{thm}

\subsection{Capacity estimates} In order to simplify the statement of the next results, set
\begin{align*}
c_{k,d}\eqdef \sqrt{\tfrac{k}{d-1}}
\comma\quad
s_{k,d}(r)\eqdef
\begin{cases}
\sin(c_{k,d}\, r) & \text{if } k>0
\\
r & \text{if } k=0
\\
\sinh(c_{k,d}\, r) & \text{if } k<0
\end{cases}
\comma \quad
V_{k,d}(r)\eqdef \int_0^r s_{k,d}(u)^{d-1}\diff u
\end{align*}
and
\begin{align}\label{eq:volBall}
v_r(x)\eqdef \mssm B^\mssg_r(x)\comma\quad x\in \man\comma r>0\comma \qquad D_\mssg\eqdef \diam_{\mssd_\mssg} \man\fstop
\end{align}

The following is well-known.
\begin{prop}[Bishop--Gromov volume comparison] Let~$(\man,\mssg)$ be satisfying Assumption~\ref{ass:RM}. Then,~$(\man,\mssd_\mssg,\mssm)$ satisfies
\begin{align}\label{eq:p:BG}
\frac{v_R(x)}{v_r(x)}\leq \frac{V_{k,d}(R)}{V_{k,d}(r)}\comma \qquad x\in \man\comma 0<r<R
\end{align}
where~$k\eqdef (d-1)\inf_{x\in \man} \Ric_x$.
\end{prop}

\begin{defs}[Packings and coverings] For fixed~$r>0$, we say that~$\seq{x_j}_j^n\subset \man$ is 
\begin{enumerate*}[label=\bfseries({\itshape\alph*})]
\item an \emph{$r$-packing} of~$\man$ if~$\tseq{B^{\mssg}_r(x_j)}_{j\leq n}$ is a disjoint family and $\cup_j^n B^{\mssg}_r(x_j)$ intersects every ball of radius~$r$ in~$\man$;
or
\item an~\emph{$r$-covering} of~$\man$ if~$\man\subset \cup_j^n B^{\mssg}_r(x_j)$.
\end{enumerate*}
The \emph{covering number} of~$\man$ is defined by
\begin{align*}
c_{\man,\mssg}(r)\eqdef \min\tset{n\in \N : \exists \seq{x_j}_j^n \ \ r\text{-covering of~} \man} \fstop
\end{align*}

We say that an $r$-covering~$\seq{x_j}_j^n$ is ($r$-)\emph{optimal} if~$n=c_{\man,\mssg}(r)$.
\end{defs}

The following is an exercise in~\cite[$E_+$,~5.31~Ex.~($b$)]{Gro06}.
We provide a proof for completeness.

\begin{lem}[Covering number of~$\man$]\label{l:CovNum}
Let~$r>0$. Then,
\begin{align*}
c_{\man,\mssg}(r)\leq \frac{V_{k,d}(D_\mssg/2)}{V_{k,d}(r/2)}\fstop
\end{align*}
\begin{proof} 
Alternatively letting either~$r\rar 0$ or~$R\rar D_\mssg/2$ in~\eqref{eq:p:BG} we have
\begin{align}\label{eq:BishopGromov}
r^d\lesssim \beta \frac{V_{k,d}(r)}{V_{k,d}(D_\mssg/2)}\leq v_r(x)\leq V_{k,d}(r) \lesssim r^d \comma \qquad r>0\comma x\in \man\fstop
\end{align}

Let~$\seq{x_j}_j^n$ be an $r/2$-packing of~$\man$. Note that it is an $r$-covering. By disjointness,
\begin{align*}
\mssm\ttonde{ \cup_j^n B^{\mssg}_{r/2}(x_j) }= \sum_j^n \mssm B^{\mssg}_{r/2}(x_j) \leq \mssm \man= \beta \comma
\end{align*}
hence~$n\leq \beta/\inf_{x\in \man} v_{r/2}(x)$ and the conclusion follows by~\eqref{eq:BishopGromov}.
\end{proof}
\end{lem}

\begin{lem}\label{l:ProdCap} For~$i=1,2$ let~$(\man_i,\mssg_i)$ be satisfying Assumption~\ref{ass:RM}, with canonical form $(\mssE^i,\dom{\mssE^i})$. Denote by~$(\man,\mssg)$ the product manifold~$(\man_1,\mssg_1)\times (\man_2,\mssg_2)$, with canonical form $(\mssE^\mssg,\dom{\mssE^\mssg})$.
For every $\mssE^i$-capacitable~$A_i,B_i$ with~$A_i\subset B_i$ let further~$u_i\eqdef u_{A_i,B_i}\in \dom{\mssE^i}$ be the equilibrium potential of~$(A_i,B_i)$.
Then, the set~$A_1\times A_2$ is $\mssE^\mssg$-capacitable and
\begin{align*}
\Cap_{\mssg}(A_1\times A_2,B_1\times B_2)\leq  \Cap_{1}(A_1,B_1)\norm{u_2}_{L^2(\n\mssm_2)}^2+ \Cap_{2}(A_2,B_2) \norm{u_1}_{L^2(\n\mssm_1)}^2 \fstop
\end{align*}

\begin{proof} Straightforward.
\end{proof}
\end{lem}

\begin{prop}\label{p:ProdCap}
For~$i=1,2$ let~$(\man,\mssg_i)$ be satisfying Assumption~\ref{ass:RM}, with canonical form $(E^i,\dom{E^i})$ and \emph{same} underlying differential manifold~$\man$. Denote by~$(\man^{\times 2},\mssg)$ the product manifold~$(\man,\mssg_1)\times (\man,\mssg_2)$, with canonical form $(E^\mssg,\dom{E^\mssg})$. Then, $\Cap_{\mssg}(\Delta \man)=0$.

\begin{proof}
Fix~$0<\eps<\delta<1$ and let~$\tseq{x_j^i}_j^{n_i}$ be an optimal $\eps$-covering for~$(\man,\mssg_i)$. Then, their union, relabeled~$\seq{y_j}_j^{n}$, is an~$\eps$-covering of both~$(\man,\mssg_i)$ and~$n\eqdef n_1+n_2\lesssim \eps^{-d}$ by Lemma~\ref{l:CovNum}. 
For simplicity of notation, for any~$r>0$ write~$B^i_r(x)\eqdef B^{\mssg_i}_r(x)$, $B^i_{r,j}\eqdef B^{\mssg_i}_r(y_j)$ and~$v_{i,r}(x)\eqdef v^{\mssg_i}_r(x)$.
Let~$u_{i,j,\eps,\delta}$ be the equilibrium potential of the pair~$(B^i_{\eps,j}, B^i_{\delta,j})$ for the form~$(\mssE^i,\dom{\mssE^i})$.
Since~$\seq{y_j}_j^n$ is a covering,~$\Delta \man\subset \cup_j^n B^1_{\eps,j}\times B^2_{\eps,j}$, thus
\begin{align*}
\Cap_{\mssg} (\Delta \man)\leq&\ \Cap_{\mssg}\tonde{ \cup_j^n \ttonde{B^1_{\eps,j}\times B^2_{\eps,j}}} \leq \sum_j^n \Cap_{\mssg}\ttonde{B^1_{\eps,j}\times B^2_{\eps,j}}
\\
\leq& \sum_j^n \Cap_{\mssg}\ttonde{B^1_{\eps,j}\times B^2_{\eps,j}, B^1_{\delta,j}\times B^2_{\delta,j}}
\\
\leq& \sum_j^n \Cap_{1}\ttonde{B^1_{\eps,j},B^1_{\delta,j}}\norm{u_{2,j,\eps,\delta}}_{L^2(\n\mssm_2)}^2+\Cap_{2}\ttonde{B^2_{\eps,j},B^2_{\delta,j}}\norm{u_{1,j,\eps,\delta}}_{L^2(\n\mssm_1)}^2
\end{align*}
by Lemma~\ref{l:ProdCap}. As a consequence, since~$0\leq u_{i,j,\eps,\delta}\leq\car_{B^i_{\delta,j}}$,
\begin{align}\label{eq:Estimate1}
\Cap_{\mssg} (\Delta \man)\lesssim &\,\eps^{-d}\, \sup_{x\in \man} \big(\Cap_{1}\ttonde{B^1_\eps(x),B^1_{\delta}(x)}\cdot v_{2,\delta}(x)
\\
\nonumber
&\qquad\qquad\qquad+\Cap_{2}\ttonde{B^2_\eps(x),B^2_\delta(x)}\cdot v_{1,\delta}(x) \big) \fstop
\end{align}

Now, if~$i,j=1,2$, $i\neq j$,
\begin{align}\label{eq:Estimate2}
\Cap_{i}\ttonde{B^i_\eps(x),B^i_\delta(x)} \cdot v_{j,\delta}(x)\leq \Cap_{i}^{(0)}\ttonde{B^i_\eps(x), B^i_\delta(x)} \cdot v_{j,\delta}(x)+ v_{1,\delta}(x)\cdot v_{2,\delta}(x) \fstop
\end{align}

By~\cite[Eqn.~(2.2)]{Stu95}, also cf.~\cite[Eqn.~(2.2)]{Gri99}, and~\eqref{eq:BishopGromov}, one has
\begin{align}
\nonumber
\sup_{x\in \man}\,& \Cap_{i}^{(0)}\ttonde{B_\eps^i(x), B_\delta^i(x)}\leq
\\
\leq& \sup_{x\in \man}\tonde{\int_\eps^\delta \frac{r-\eps}{v_{i,r}(x)-v_{i,\eps}(x)} \diff r}^{-1}
\\
\nonumber
\leq& \tonde{\int_\eps^\delta \frac{r-\eps}{\sup_{x\in \man} v_{i,r}(x)- \inf_{x\in \man}v_{i,\eps}(x)} \diff r}^{-1}
\\
\nonumber
\leq& \tonde{\int_\eps^\delta \frac{r-\eps}{V_{k,d}(r)-V_{k,d}(\eps)\, \beta^{-1} V_{k,d}(D_i/2)^{-1}} \diff r}^{-1}
\\
\label{eq:Cap0}
\lesssim&\, c_\eps\eqdef 
\begin{cases} 
\eps^{d-2} & \textrm{if~~} d\geq 3 \textrm{~~and~~} \delta\eqdef 2\eps
\\
\ttonde{\ln\tfrac{\delta+\eps}{2\eps}}^{-1} & \textrm{if~~} d=2 \textrm{~~and~~} \delta\eqdef 1\wedge D_1/2\wedge D_2/2
\end{cases}
\fstop
\end{align}

Finally, combining Equations~\eqref{eq:BishopGromov}--\eqref{eq:Cap0} yields
\begin{align*}
\Cap_{\mssg}(\Delta \man)\lesssim\, \eps^{-d} (c_\eps\,\eps^d+\eps^d\,\eps^d)
\end{align*}
and letting~$\eps$ tend to~$0$ yields the desired conclusion since~$c_\eps\rar 0$.
\end{proof}
\end{prop}

\begin{lem}\label{l:Coexcept}
Let~$\mbfs\in\Tso$ and~$\mssW^{n,\mbfs}_\bullet$ be defined as in~\eqref{eq:BMn}. Then,~$\andi{\man^{\times n}}$ is~$\mssW^{n,\mbfs}_\bullet$-co\-ex\-cep\-tion\-al.

\begin{proof}
By the standard theory of Dirichlet forms, the statement is equivalent to the set~$\andi{\man^{\times n}}$ being $\mssE^{n,\mbfs}$-coexceptional; see, e.g.~\cite[Thm.~4.1.2(i)]{FukOshTak11} or~\cite[Thm.\ IV.5.29(i)]{MaRoe92}.
Let
\begin{align*}
\man_{i,j}^n\eqdef \set{\mbfx\in \man^{\times n}: x_i=x_j \text{ for } i\neq j} \fstop
\end{align*}

Since~$\man^{\times n}\setminus\andi{\man^{\times n}}\subset \cup_{i,j: i\neq j} \man_{i,j}^n$, it suffices to show that~$\Cap_{n,\mbfs}(\man_{i,j}^n)=0$, where $\Cap_{n,\mbfs}$ denotes the capacity associated to~$(\mssE^{n,\mbfs},\dom{\mssE^{n,\mbfs}})$. Without loss of generality, we can assume~$i=1$,~$j=2$.
Set~$\mssg_i\eqdef s_i\,\mssg$,~$\mssd_{1,2}\eqdef \mssd_{\mssg_1\oplus \mssg_2}$ and~$B^{1,2}_\eps(A)\eqdef B^{\mssd_{1,2}}_\eps(A)\subset \man^{\times 2}$ be the $\eps$-neighborhood of~$A\subset \man^{\times 2}$. 
Denote by~$\Cap_{1,2}$, resp.~$\Cap_{3,\dotsc,n}$, the capacity of the form~$\ttonde{\mssE^{2,(s_1,s_2)},\dom{\mssE^{2,(s_1,s_2)}} }$, resp.~$\ttonde{\mssE^{n-2,(s_3,\dotsc,s_n)},\dom{\mssE^{n-2,(s_3,\dotsc,s_n)}} }$.
For~$0<\eps<r$, let now~$u_{1,2,\eps}$ be the equilibrium potential of~$B^{1,2}_\eps(\Delta\man)$ for the~$\Cap_{1,2}$. By Lemma~\ref{l:ProdCap},
\begin{align*}
\Cap_{n,\mbfs}(\man^n_{1,2})\leq&\ \Cap_{n,\mbfs}\ttonde{B_\eps^{1,2}(\Delta \man)\times \man^{\times n-2}, \man^{\times 2}\times \man^{\times n-2}}
\\
\leq&\ \Cap_{1,2}\ttonde{B_\eps^{1,2}(\Delta \man)}\norm{\car}_{L^2(\n\mssm^{n-2})}^2
\\
&\ +\Cap_{3,\dotsc,n}(\man^{n-2},\man^{n-2})\norm{u_{1,2,\eps}}_{L^2(\n\mssm^2)}^2
\\
\leq&\ \Cap_{1,2}\ttonde{B_\eps^{1,2}(\Delta \man)}\cdot 1+1\cdot \norm{u_{1,2,\eps}}_{L^2(\n\mssm^2)}^2
\\
\leq&\ 2\,\Cap_{1,2}\ttonde{B_\eps^{1,2}(\Delta \man)} \fstop
\end{align*}

The conclusion follows by~Proposition~\ref{p:ProdCap} letting~$\eps\rar 0$.
\end{proof}
\end{lem}

\begin{lem}\label{l:Except}
Let~$\mbfs\in\Tso$ and~$\BM^\mbfs_\bullet$ be defined as in~\eqref{eq:Intro:BM}. Then,~$\andi{\mbfM}$ is $\BM^\mbfs_\bullet$-coexceptional.

\begin{proof}
Denote by~$\tau^\mbfs$ the \emph{first touching time} of~$\andi{\mbfM}^\complement$ for~$\BM^\mbfs_\bullet$ in the sense of~\cite[\S{IV.5}, Eqn.~(5.14)]{MaRoe92}.
Since~$\andi{\mbfM}^\complement$ is measurable and~$\BM^\mbfs_\bullet$ has infinite life-time,  it suffices to show
\begin{align}\label{eq:l:Except1}
P^\mbfs_{\n\volm} \set{\tau^\mbfs<\infty}=0\comma
\end{align}
where~$P^\mbfs_{\n\volm}$ is defined analogously to~\cite[\S{IV.1}, Eq.~(1.4)]{MaRoe92}.
With slight abuse of notation, for every~$\mbfx_0\in\mbfM$ we denote both~$\mbfx_0$ and~$\Proj{n}{}(\mbfx_0)$ by~$\mbfx_0$, the distinction being apparent from the contextual index~$n$.
Note that~$\mbfM\setminus\andi{\mbfM}=\cap_n \Proj{-1}{n}\ttonde{\man^{\times n}\setminus \andi{\man^{\times n}}}$ and $\Proj{n}{}\circ \BM^{\mbfs;\mbfx_0}_\bullet=\mssW^{n,\mbfs;\mbfx_0}_\bullet$ by Theorem~\ref{t:ADKBSC}\iref{i:t:ADKBSC:6}. Moreover, since~$\man^{\times n}$ is compact,~$\Proj{n}{}$ is a closed map, hence~$\Proj{n}{}\ttonde{\overline{\BM^\mbfs_{[0,t]}}}=\overline{\mssW^{n,\mbfs}_{[0,t]}}$ by~\cite[Cor.~3.1.11]{Eng89}.
As a consequence, letting~$\tau^{n,\mbfs}$ be the first touching time of~$\man^{\times n}\setminus \andi{\man^{\times n}}$ for~$\mssW^{n,\mbfs}_\bullet$,
\begin{align}\label{eq:l:Except2}
\set{\tau^{n,\mbfs}<\infty}\subset \set{\tau^{n+1,\mbfs}<\infty} \subset \set{\tau^\mbfs<\infty} = \nlim \set{\tau^{n,\mbfs}<\infty} \fstop
\end{align}
Furthermore,
\begin{align}\label{eq:l:Except3}
P^\mbfs_{\n\volm}\set{\tau^{n,\mbfs}<\infty}= P^{n,\mbfs}_{\n\mssm^n}\set{\tau^{n,\mbfs}<\infty}=0
\end{align}
since~$\man^{\times n}\setminus\andi{\man^{\times n}}$ is~$\mssW^{n,\mbfs}_\bullet$-exceptional for every~$\mbfs\in\Tso$ by Lemma~\ref{l:Coexcept}.
Finally,~\eqref{eq:l:Except2} and~\eqref{eq:l:Except3} yield~\eqref{eq:l:Except1} by the Borel--Cantelli Lemma.
\end{proof}
\end{lem}

\subsection{Operators and domains}


\begin{lem}\label{l:PrepTFDensity} Let~$\varrho\in \mcC^\infty(I)$,~$f\in\mcC^\infty(\man)$ and~$g\in\mcC^\infty(\man^{\times 2})$. Then,
\begin{align*}
u_{g,\varrho}\colon \mu\mapsto \int_\man f(x)\cdot \varrho\tonde{\int_\man g(x,y) \diff\mu(y)} \diff\mu(x)
\end{align*}
satisfies~$u_{g,\varrho}\in\cl_{\mcE^{1/2}_1}(\TF{\infty})$ and
\begin{equation}\label{eq:l:PrepTFDensity:0}
\begin{aligned}
\ttonde{\mssg^\flat \grad u_{g,\varrho}(\eta)(x)}(\emparg)=&(\diff f)_x(\emparg_x)\cdot \varrho\ttonde{g(x,\emparg)^\trid \eta}
\\
&+f(x)\cdot\varrho'\ttonde{g(x,\emparg)^\trid\eta}\int_\man (\diff^{\otimes 2} g)_{x,y}(\emparg_x,\emparg_y)\diff\eta(y) \fstop
\end{aligned}
\end{equation}
\begin{proof}
Let~$u\eqdef u_{g,\varrho}$. Note that
\begin{align*}
\grad_w u(\eta)=&\diff_t\restr_{t=0}\int_\man (f\circ \uppsi^{w,t})(x)\cdot \varrho\tonde{\int_\man g\ttonde{\uppsi^{w,t}(x),\uppsi^{w,t}(y)} \diff\eta(y)} \diff\eta(x)
\\
=&\int_\man \tgscal{\nabla_x f}{w_x} \cdot \varrho\ttonde{g(x,\emparg)^\trid \eta} \diff\eta(x)
\\
&+\int_\man f(x) \cdot \varrho'\ttonde{g(x,\emparg)^\trid \eta}\cdot \int_\man \tscalar{\nabla^{\otimes 2}_{x,y} g}{(w_x,w_y)}_{\mssg_x\oplus\mssg_y} \diff\eta(y) \diff\eta(x)\comma
\end{align*}
whence~\eqref{eq:l:PrepTFDensity:0} follows.
By a straightforward approximation argument in the appropriate~$\mcC^1$-topologies, it suffices to show~$u_{g,\varrho}\in\cl_{\mcE^{1/2}_1}(\TF{\infty})$ when~$\varrho\in I[r]$, the space of real-valued polynomials on~$I$, and~$g$ is of the form~$g=\sum_k^n a_k\otimes b_k$, where~$k\leq n\in \N$ and~$a_k,b_k\in\mcC^\infty(\man)$.
Finally, since~$\varrho\mapsto u_{g,\varrho}$ is linear and~$\grad$ is a linear operator, it suffices to show the statement when~$\varrho(r)\eqdef r^m$ for~$m\in \N$.
For such a choice of~$g$ and~$\varrho$, one has in fact
\begin{align*}
u_{g,\varrho}(\eta)=&\sum_{\mbfj\in \N_0^n \,:\, \abs{\mbfj}=m} \tbinom{m}{\mbfj} (f\cdot \mbfa^\mbfj)^\trid\eta \cdot (\mbfb^\trid\eta)^\mbfj\in\TF{\infty}\comma\qquad \mbfa\eqdef \seq{a_k}_k^n\comma\mbfb\eqdef\seq{b_k}_k^n \fstop \qedhere
\end{align*}
\end{proof}
\end{lem}

\begin{lem}\label{l:TFDensity}
The set~$\TF{\infty}$ is dense in~$\dom{\mcE}$.
\begin{proof}
In order to prove the statement, it suffices to show that~$u\eqdef (f\otimes\varrho)^\trid\in \cl_{\mcE^{1/2}_1}(\TF{\infty})$ for all~$f\in\mcC^\infty(\man)$ and~$\varrho\in\mcC^\infty(I)$.
Denote by~$\Cap$ the (first order) capacity associated to the canonical form~$(\mssE,\dom{\mssE})$ of~$(\man,\mssg)^{\times 2}$. Let~$\seq{g_n}_n\in\dom{\mssE}$ be a minimizing sequence for~$\Cap (\Delta \man)=0$ by Prop.~\ref{p:ProdCap}. By standard arguments, we may assume that~$g_n\in\mcC^\infty(\man^{\times 2})$ additionally satisfies
\begin{equation*}
\begin{aligned}
0\leq g_n\leq& 1\comma \quad & g_n(x,x)=&1\comma
\\
\nabla^z\restr_{z=x}g_n(z,z)=&0\comma \quad & \mssE_1(g_n)\leq& 2^{-n} \comma
\end{aligned}
\qquad x\in \man\comma n\in \N \comma
\end{equation*}
and further that~$g_n(x,y)=g_n(y,x)$, so that we may write~$g^x_n\eqdef g_n(x,\emparg)$, unambiguously.

By Lemma~\ref{l:PrepTFDensity}, for every~$f\in\mcC^\infty(\man)$,~$\varrho\in I[r]$ and~$n\in \N$,
\begin{align*}
u_n(\eta)\eqdef \int_\man f(x) \cdot \varrho\ttonde{(g^x_n)^\trid \eta} \diff\eta(x) \in \cl_{\mcE^{1/2}_1}(\TF{\infty}) \fstop
\end{align*}
Thus, it suffices to show that~$\mcE^{1/2}_1$-$\nlim u_n=u$. 
By~\eqref{eq:l:PrepTFDensity:0},
\begin{align*}
\mcE(u-u_n)\leq& \norm{\abs{\nabla f}}_\infty^2\underbrace{\int_\msP\int_\man \abs{\varrho(\eta_x)-\varrho\ttonde{(g^x_n)^\trid \eta}}^2 \diff\eta(x)\diff\DF_\mssm(\eta)}_{I_{1,n}}
\\
&+\norm{f}_\infty^2\underbrace{\int_\msP \int_\man \abs{\varrho'\ttonde{(g^x_n)^\trid \eta}}^2 \int_\man \abs{\nabla^{\otimes 2}_{x,y} g_n}^2 \diff\eta(y)\diff\eta(x) \diff\DF_\mssm(\eta)}_{I_{2,n}} \fstop
\end{align*}

Concerning~$I_{1,n}$, one has, by the Mecke identity~\eqref{eq:Mecke} and properties of~$g^x_n$, that
\begin{align*}
I_{1,n}\leq& \norm{\varrho'}_\infty^2 \int_\msP \int_\man \int_I \abs{r-rg^x_n(x)-(1-r)\int_\man g^x_n(y) \diff\eta(y)}^2 \diff\Beta(r)\diff\n\mssm(x)\diff\DF_\mssm(\eta)
\\
\leq& \norm{\varrho'}_\infty^2 \int_\msP \int_\man \int_\man \abs{g^x_n(y)}^2 \diff\eta(y)\diff\n\mssm(x)\diff\DF_\mssm(\eta)
\\
=& \norm{\varrho'}_\infty^2 \int_{\man^{\times 2}} \abs{g^x_n(y)}^2 \diff\n\mssm^{\otimes 2}(x,y)
\leq 2^{-n}\norm{\varrho'}_\infty^2
\end{align*}
and therefore vanishing as~$n\rar\infty$.
A proof of the convergence~$\nlim u_n=u$ pointwise on~$\Ppa$ and in~$L^2(\msP,\DF_\mssm)$ is analogous to that for~$I_{1,n}$ and therefore it is omitted.

Concerning~$I_{2,n}$, by Cauchy--Schwarz inequality, properties of~$g_n$ and the Mecke identity~\eqref{eq:Mecke},
\begin{align*}
&I_{2,n}\leq
\\
\leq&\norm{\varrho'}_\infty^2\int_\msP \int_\man \int_I \quadre{(1-r) \int_\man \abs{\nabla^{\otimes 2}_{x,y} g_n}^2 \diff\eta(y)+r\abs{\nabla^{\otimes 2}_{x,x} g_n}^2}\diff\Beta_\beta(r)\diff\n\mssm(x)\diff\DF_\mssm(\eta)
\\
\leq& \norm{\varrho'}_\infty^2\int_\man \abs{\nabla^{\otimes 2}_{x,y} g_n}^2 \diff\n\mssm^{\otimes 2}(x,y)
\leq 2^{-n} \norm{\varrho'}_\infty^2\comma
\end{align*}
which concludes the proof by letting~$n\rar\infty$.
\end{proof}
\end{lem}

\begin{lem}\label{l:Rademacher}
For~$w\in\Vect^\infty$ let $\mcA^w$ be the form on~$\dom{\mcE}$ defined in Corollary~\ref{c:Martingale}\iref{i:c:Martingale2}. Then, for all bounded measurable~$u\colon \msP\rar \R$ and all~$v\in\TF{\infty}$,
\begin{align}\label{eq:l:Rademacher:0}
\int_\msP (u\circ \Psi^{w,t}-u) v \diff\DF_\mssm =-\int_0^1 \int_\msP \ttonde{u\circ \Psi^{w,s} \cdot \grad_w v+ \mcA^w(u\circ \Psi^{w,s},v) } \diff\DF_\mssm \diff s \fstop
\end{align}
\begin{proof}
By a monotone class argument, it suffices to show~\eqref{eq:l:Rademacher:0} for~$u\in \TF{\infty}$. Then,~$u\circ\Psi^{w,t}\in\TF{\infty}$ too.
By~\cite[Lem.~4.7]{LzDS19b},
\begin{align*}
u\circ\Psi^{w,t}-u=\int_0^t \grad_w (u\circ \Psi^{w,s}) \diff s \comma
\end{align*}
whence, integrating and applying Fubini's Theorem,
\begin{align*}
\int_\msP (u\circ\Psi^{w,t}-u)v \diff\DF_\mssm =\int_0^t \int_\msP \grad_w (u\circ \Psi^{w,s}) \cdot v \, \diff\DF_\mssm \diff s \comma
\end{align*}
hence the conclusion by properties of~$\mcA^w$.
\end{proof}
\end{lem}

\begin{prop}\label{p:Rademacher}
For~$u\in \Lip(\msP_2)$ and~$w\in\Vect^\infty$ set
\begin{align*}
\Omega^u_w\eqdef \set{\mu\in\msP : \exists\, G_w u(\mu)\eqdef \diff_t\restr_{t=0}(u\circ \Psi^{w,t})(\mu)} \fstop
\end{align*}

Let further~$\msX\subset \Vect^\infty$ be a countable $\Q$-vector space dense in~$\Vect^0$ and assume $\DF_\mssm\Omega^u_w=1$ for all~$w\in\msX$. Then~$u\in\dom{\mcE}$ and~$\boldGamma(u)\leq \Lip[u]$ $\DF_\mssm$-a.e..

\begin{proof}
It suffices to show~\cite[Eqn.~(4.20)]{LzDS19b}: the rest of the proof is identical to~\cite[Prop.~4.9]{LzDS19b}. By continuity of~$t\mapsto\mcA^w(u\circ \Psi^{w,t},v)$ for~$u,v\in \TF{\infty}$,~\cite[Eqn.~(4.19)]{LzDS19b} and Lemma~\ref{l:Rademacher} yield
\begin{align*}
\int_\msP G_wu\cdot v \, \diff\DF_\mssm=-\int_\msP u\grad_w v \diff\DF_\mssm - \mcA^w(u,v)\comma\qquad u,v\in\TF{\infty}
\end{align*}

Next, note that the map~$w\mapsto \mcA^w(u,v)$ is linear for every~$u,v\in\TF{\infty}$, since it is the limit of the linear maps~$w\mapsto \mcA^w_0(u_n,v_n)$, where~$u_n,v_n\in\hTF{\infty}{1/n}$ are the approximation of~$u$,~$v$ constructed in Corollary~\ref{c:Domains}. As a consequence, if~$w=s_1w_1+\cdots+s_kw_k$ for some~$s_i\in \R$ and~$w_i\in\msX$, then
\begin{align*}
\int_\msP G_w u\cdot v\, \diff\DF_\mssm =-\sum_i^k s_i\tonde{\int_\msP u\grad_{w_i} v\, \diff\DF_\mssm+\mcA^{w_i}(u,v)}=\sum_i^k \int_\msP G_{w_i} u\cdot v \diff\DF_\mssm
\end{align*}
and~\cite[Eqn.~(4.20)]{LzDS19b} follows by arbitrariness of~$v$.
\end{proof}
\end{prop}


{\small


\begin{thebibliography}{10}

\bibitem{AirMal06}
{Airault, H.} and {Malliavin, P.}
\newblock {Quasi-invariance of Brownian measures on the group of circle
  homeomorphisms and infinite-dimensional Riemannian geometry}.
\newblock {\em {J. Funct. Anal.}}, 241:99--142, 2006.

\bibitem{AlbDalKon97}
{Albeverio, S.}, {Daletskii, A.~Yu.}, and {Kondrat'ev, Yu.~G.}
\newblock {Infinite Systems of Stochastic Differential Equations and Some
  Lattice Models on Compact Riemannian Manifolds}.
\newblock {\em {Ukr. Math. J.}}, 49(3):360--372, 1997.

\bibitem{AlbDalKon00}
{Albeverio, S.}, {Daletskii, A.~Yu.}, and {Kondrat'ev, Yu.~G.}
\newblock {Stochastic Analysis on Product Manifolds: Dirichlet Operators on
  Differential Forms}.
\newblock {\em {J. Funct. Anal.}}, 176:280--316, 2000.

\bibitem{AlbRoe95}
{Albeverio, S.} and {R{\"{o}}ckner, M.}
\newblock {Dirichlet Form Methods for Uniqueness of Martingale Problems and
  Applications}.
\newblock In {Cranston, M.} and {Pinsky, M.}, editors, {\em {Stochastic
  Analysis -- Proceedings of the Summer Research Institute on Stochastic
  Analysis, Held at Cornell University, Ithaca, New York, July 11--30, 1993}},
  volume~57 of {\em {Proceedings of Symposia in Pure Mathematics}}, pages
  513--528. {Amer.\ Math.\ Soc.}, 1995.

\bibitem{AmbFusPal00}
{Ambrosio, L.}, {Fusco, N.}, and {Pallara, D.}
\newblock {\em {Functions of Bounded Variation and Free Discontinuity
  Problems}}.
\newblock {Oxford Mathematical Monographs}. {Oxford Science Publications},
  2000.

\bibitem{AmbGig11}
{Ambrosio, L.} and {Gigli, N.}
\newblock {A User's Guide to Optimal Transport}.
\newblock In {Ambrosio, L.}, {Bressan, A.}, {Helbing, D.}, {Klar, A.}, and
  {Zuazua, E.}, editors, {\em {Modelling and Optimisation of Flows on Networks
  -- Cetraro, Italy 2009, Editors: Benedetto Piccoli, Michel Rascle}}, volume
  {2062} of {\em {Lecture Notes in Mathematics}}, pages 1--155. {Springer},
  2013.
\newblock Throughout the present work, we refer to (the numbering of) results
  in the extended version, available at
  \url{http://cvgmt.sns.it/media/doc/paper/195/}.

\bibitem{BanRad97}
{Banakh, T.~O.} and {Radul, T.~N.}
\newblock {Topology of spaces of probability measures}.
\newblock {\em {Mat.\ Sb.}}, 188(7):23--46, 1997.

\bibitem{BenSaC97}
{Bendikov, A.} and {Saloff-Coste, L.}
\newblock {Elliptic diffusions on infinite products}.
\newblock {\em {J. reine angew. Math.}}, 493:171--220, 1997.

\bibitem{BerKon95}
{Berezansky, Yu. M.} and {Kondratiev, Yu. G.}
\newblock {\em {Spectral Methods in Infinite-Dimensional Analysis}}.
\newblock Springer Netherlands, Dordrecht, 1995.

\bibitem{Ber76}
{Berg, Ch.}
\newblock {Potential Theory on the Infinite Dimensional Torus}.
\newblock {\em {Invent. Math.}}, 32:49--100, 1976.

\bibitem{BerKlo15}
{Bertrand, J.} and {Kloeckner, B. R.}
\newblock {A Geometric Study of Wasserstein Spaces: Isometric Rigidity in
  Negative Curvature}.
\newblock {\em {International Mathematics Research Notices}}, 2015.

\bibitem{Bog07}
{Bogachev, V.~I.}
\newblock {\em {Measure Theory}}.
\newblock {Springer-Verlag}, {Berlin}, {2007}.

\bibitem{BouHir91}
{Bouleau, N.} and {Hirsch, F.}
\newblock {\em {Dirichlet forms and analysis on Wiener space}}.
\newblock {De Gruyter}, 1991.

\bibitem{Bre91}
{Brenier, Y.}
\newblock {Polar Factorization and Monotone Rearrangement of Vector-Valued
  Functions}.
\newblock {\em {Comm. Pur. Appl. Math.}}, 44:375--417, 1991.

\bibitem{CheMaRoe94}
{Chen, Z.-Q.}, {Ma, Z.-M.}, and {R\"ockner, M.}
\newblock {Quasi-homeomorphisms of Dirichlet forms}.
\newblock {\em {Nagoya Math. J.}}, 136:1--15, 1994.

\bibitem{ChoGan17}
{Chow, Y.~T.} and {Gangbo, W.}
\newblock {A partial Laplacian as an infinitesimal generator on the Wasserstein
  space}.
\newblock {\em {J.\ Differ.\ Equat.}}, 267(10):6065--6117, 2019.

\bibitem{DaPZab96}
{Da Prato, G.} and {Zabczyk, J.}
\newblock {\em {Ergodicity for infinite-dimensional systems}}, volume 229 of
  {\em {London Mathematical Society Lecture Note Series}}.
\newblock {Cambridge University Press}, 1996.

\bibitem{LzDS19b}
{Dello Schiavo, L.}
\newblock {A Rademacher-type theorem on $L^2$-Wasserstein spaces over closed
  Riemannian manifolds}.
\newblock {\em {J.\ Funct.\ Anal.}}, 278(108397), 2019.

\bibitem{LzDS19a}
{Dello Schiavo, L.}
\newblock {Characteristic functionals of Dirichlet measures}.
\newblock {\em {Electron.\ J.\ Probab.}}, 24(115):1--38, 2019.

\bibitem{LzDS20}
{Dello Schiavo, L.}
\newblock {Ergodic Decomposition of Dirichlet Forms via Direct Integrals and
  Applications}.
\newblock march 2020.
\newblock {arXiv:2003.01366}.

\bibitem{LzDSLyt17}
{Dello Schiavo, L.} and {Lytvynov, E.}
\newblock {A Mecke-type Characterization of the Dirichlet--Ferguson Measure}.
\newblock {\em arXiv 1706.07602}, 2017.

\bibitem{DonGri93}
{Donnelly, P.} and {Grimmett, G.}
\newblock {On the asymptotic distribution of large prime factors}.
\newblock {\em {J. London Math. Soc.}}, 47(3):395--404, 1993.

\bibitem{DonJoy89}
{Donnelly, P.} and {Joyce, P.}
\newblock {Continuity and weak convergence of ranked and size-biased
  permutations on the infinite simplex}.
\newblock {\em {Stoch. Proc. Appl.}}, 31(1):89--103, 1989.

\bibitem{DoeSta09}
{D\"oring, M.} and {Stannat, W.}
\newblock {The logarithmic Sobolev inequality for the Wasserstein diffusion}.
\newblock {\em {Probab.\ Theory Relat.\ Fields}}, 145:189--209, 2009.

\bibitem{Eng89}
{Engelking, R.}
\newblock {\em {General Topology}}, volume~6 of {\em {Sigma series in pure
  mathematics}}.
\newblock {Heldermann}, {Berlin}, {1989}.

\bibitem{EthKur94}
{Ethier, S. N.} and {Kurtz, T. G.}
\newblock {Convergence to Fleming--Viot process in the weak atomic topology}.
\newblock {\em {Stoch. Proc. Appl.}}, 54(1):1--27, {1994}.

\bibitem{Fer73}
{Ferguson, T. S.}
\newblock {A Bayesian analysis of some nonparametric problems}.
\newblock {\em {Ann.\ Statist.}}, 1:{209--230}, 1973.

\bibitem{FleVio79}
{Fleming, V.~H.} and {Viot, M.}
\newblock {Some Measure-Valued Markov Processes in Population Genetics Theory}.
\newblock {\em {Indiana J.\ Math.}}, 28(5):817--843, 1979.

\bibitem{FukOshTak11}
{Fukushima, M.}, {Oshima, Y.}, and {Takeda, M.}
\newblock {\em {Dirichlet forms and symmetric Markov processes}}, volume~19 of
  {\em {De Gruyter Studies in Mathematics}}.
\newblock {de Gruyter}, extended edition, 2011.

\bibitem{GanKimPac10}
{Gangbo, W.}, {Kim, H. K.}, and {Pacini, T.}
\newblock {Differential Forms on Wasserstein Space and Infinite-dimensional
  Hamiltonian Systems}.
\newblock {\em {Mem. Am. Math. Soc.}}, 211(995), 2010.

\bibitem{Gig11}
{Gigli, N.}
\newblock {On the inverse implication of Brenier--McCann theorems and the
  structure of $\tonde{\msP_2(M),W_2}$}.
\newblock {\em {Methods and Applications of Analysis}}, 18(2):127--158, 2011.

\bibitem{Gig12}
{Gigli, N.}
\newblock {Second Order Analysis on~$\tonde{\msP_2(M),W_2}$}.
\newblock {\em {Mem. Am. Math. Soc.}}, 216(1018), 2012.

\bibitem{Gri99}
{Grigor'yan, A.}
\newblock {Isoperimetric inequalities and capacities on Riemannian manifolds}.
\newblock In {Rossmann, J.}, {Tak\'a\v{c}, P.}, and {Wildenhain, G.}, editors,
  {\em {The Maz'ya Anniversary Collection -- Volume 1: On Maz'ya's work in
  functional analysis, partial differential equations and applications}},
  volume 109 of {\em {Operator Theory -- Advances and Applications}}, pages
  139--153. {Springer Basel AG}, 1999.

\bibitem{Gri09}
{Grigor'yan, A.}
\newblock {\em {Heat Kernel and Analysis on Manifolds}}, volume~47 of {\em
  {Advanced Studies in Mathematics}}.
\newblock {AMS--IP}, 2009.

\bibitem{Gro06}
{Gromov, M.}
\newblock {\em {Metric structures for Riemannian and non-Riemannian spaces}},
  volume 152 of {\em {Progress in Mathematics}}.
\newblock {Birkh\"auser}, {3} edition, {2006}.

\bibitem{Han02}
{Handa, K.}
\newblock {Quasi-invariance and reversibility in the Fleming--Viot process}.
\newblock {\em {Probab. Theory Relat. Fields}}, 122:545--566, 2002.

\bibitem{Hin09}
{Hino, M.}
\newblock {Energy measures and indices of Dirichlet forms, with applications to
  derivatives on some fractals}.
\newblock {\em {Proc. London Math. Soc.}}, 100:269--302, 2009.

\bibitem{HinRam03}
{Hino, M.} and {Ram{\'{i}}rez, J. A.}
\newblock {Small-Time Gaussian Behavior of Symmetric Diffusion Semigroups}.
\newblock {\em {Ann. Probab.}}, 31(3):1254--1295, 2003.

\bibitem{JiaDicKuo04}
{Jiang, T.~J.}, {Dickey, J.~M.}, and {Kuo, K.-L.}
\newblock {A new multivariate transform and the distribution of a random
  functional of a Ferguson--Dirichlet process}.
\newblock {\em {Stoch.\ Proc.\ Appl.}}, 111(1):77--95, 2004.

\bibitem{Kak48}
{Kakutani, S.}
\newblock {On Equivalence of Infinite Product Measures}.
\newblock {\em {Ann. Math.}}, 49(1):214--224, Jan 1948.

\bibitem{Kal17}
{Kallenberg, O.}
\newblock {\em {Random Measures, Theory and Applications}}.
\newblock {Springer}, 2017.

\bibitem{Kin75}
{Kingman, J. F. C.}
\newblock {Random Discrete Distributions}.
\newblock {\em J. Roy. Stat. Soc. B Met.}, 37(1):1--22, 1975.

\bibitem{Kol06}
{Kolesnikov, A. V.}
\newblock {Mosco convergence of Dirichlet forms in infinite dimensions with
  changing reference measures}.
\newblock {\em {J. Funct. Anal.}}, 230(2):382--418, 2006.

\bibitem{Kon17}
{Konarovskyi, V. V.}
\newblock {A System of Coalescing Heavy Diffusion Particles on the Real Line}.
\newblock {\em {Ann.\ Probab.}}, 45(5):3293--3335, 2017.

\bibitem{KonvRe17}
{Konarovskyi, V. V.} and {Renesse, M.-K. von}.
\newblock {Reversible Coalescing-Fragmentating Wasserstein Dynamics on the Real
  Line}.
\newblock {\em {arXiv:1709.02839}}, page~59, 2017.

\bibitem{KonvRe18}
{Konarovskyi, V. V.} and {Renesse, M.-K. von}.
\newblock {Modified Massive Arratia flow and Wasserstein diffusion}.
\newblock {\em {Comm. Pur. Appl. Math.}}, 2018.
\newblock (to appear).

\bibitem{KonLytVer15}
{Kondratiev, Yu.~G.}, {Lytvynov, E.~W.}, and {Vershik, A.~M.}
\newblock {Laplace operators on the cone of Radon measures}.
\newblock {\em {J.\ Funct.\ Anal.}}, 269(9):2947--2976, 2015.

\bibitem{KuwShi03}
{Kuwae, K.} and {Shioya, T.}
\newblock {Convergence of spectral structures: a functional analytic theory and
  its applications to spectral geometry}.
\newblock {\em {Comm. Anal. Geom.}}, 11(4):599--673, 2003.

\bibitem{Las18}
{Last, G.}
\newblock {An integral characterization of the Dirichlet process}.
\newblock {\em {J.\ Theor.\ Probab.}}, 49, 2019.

\bibitem{LeJRai04}
{{Le Jan, Y.} and {Raimond, O.}}
\newblock {Flows, Coalescence and Noise}.
\newblock {\em {Ann. Probab.}}, 32(2):1247--1315, 2004.

\bibitem{LenSchWir18}
{Lenz, D.}, {Schmidt, M.}, and {Wirth, M.}
\newblock {Geometric Properties of Dirichlet Forms under Order Isomorphism}.
\newblock {\em {arXiv:1801.08326}}, 2018.

\bibitem{Lot07}
{Lott, J.}
\newblock {Some Geometric Calculations on Wasserstein Space}.
\newblock {\em {Comm.\ Math.\ Phys.}}, 277(2):423--437, 2007.

\bibitem{MaRoe92}
{Ma, Z.-M.} and {R\"ockner, M.}
\newblock {\em Introduction to the Theory of (Non-Symmetric) Dirichlet Forms}.
\newblock {Graduate Studies in Mathematics}. Springer, 1992.

\bibitem{MaXia01}
{Ma, Z.-M.} and {Xiang, K.-N.}
\newblock {Superprocesses of Stochastic Flows}.
\newblock {\em {Ann. Probab.}}, 29(1):317--343, 2001.

\bibitem{McC97}
{McCann, R. J.}
\newblock {A Convexity Principle for Interacting Gases}.
\newblock {\em {Adv. Math.}}, 128, 1997.

\bibitem{Ott01}
{Otto, F.}
\newblock {The Geometry of Dissipative Evolution Equations: The Porous Medium
  Equation}.
\newblock {\em {Comm. Part. Diff. Eq.}}, 26(1-2):101--174, 2001.

\bibitem{OveRoeSch95}
{Overbeck, L.}, {R{\"o}ckner, M.}, and {Schmuland, B.}
\newblock {An analytic approach to Fleming--Viot processes with interactive
  selection}.
\newblock {\em {Ann.\ Probab.}}, 23(1):1--36, 1995.

\bibitem{ReeSim75}
{Reed, M.} and {Simon, B.}
\newblock {\em Methods of Modern Mathematical Physics II -- Fourier Analysis,
  Self-Adjointness}.
\newblock Academic Press, New York, London, 1975.

\bibitem{vReStu09}
{Renesse, M.-K. von} and {Sturm, K.-T.}
\newblock {Entropic measure and Wasserstein diffusion}.
\newblock {\em {Ann.\ Probab.}}, 37(3):1114--1191, 2009.

\bibitem{vReYorZam08}
{Renesse, M.-K. von}, {Yor, M.}, and {Zambotti, L.}
\newblock {Quasi-invariance properties of a class of subordinators}.
\newblock {\em {Stoch.\ Proc.\ Appl.}}, 118(11):2038--2057, 2008.

\bibitem{Sch97}
{Schied, A.}
\newblock {Geometric Aspects of Fleming--Viot and Dawson--Watanabe Processes}.
\newblock {\em {Ann. Probab.}}, 25(3):1160--1179, 1997.

\bibitem{Sch02}
{Schied, A.}
\newblock {Geometric Analysis for Symmetric Fleming--Viot Operators:
  Rademacher's Theorem and Exponential Families}.
\newblock {\em {Potential Anal.}}, 17:351--374, 2002.

\bibitem{Set94}
{Sethuraman, J.}
\newblock {A constructive definition of Dirichlet priors}.
\newblock {\em {Stat.\ Sinica}}, 4(2):639--650, 1994.

\bibitem{Sha11}
{Shao, J.}
\newblock {A New Probability Measure-Valued Stochastic Process with
  Ferguson-Dirichlet Process as Reversible Measure}.
\newblock {\em {Electron.\ J.\ Probab.}}, 16(9):271--292, 2011.

\bibitem{Sie22}
{Sierpi{\'n}ski, W. F.}
\newblock {Sur les fonctions d'ensemble additives et continues}.
\newblock {\em {Fund. Math.}}, 3(1):240--246, 1922.

\bibitem{Stu95}
{Sturm, K.-T.}
\newblock {Sharp estimates for capacities and applications to symmetric
  diffusions}.
\newblock {\em {Probab. Theory Relat. Fields}}, 103:73--89, 1995.

\bibitem{Stu11}
{Sturm, K.-T.}
\newblock {Entropic Measure on Multidimensional Spaces}.
\newblock In {Dalang, R.}, {Dozzi, M.}, and {Russo, F.}, editors, {\em {Seminar
  on Stochastic Analysis, Random Fields and Applications VI}}, volume~63 of
  {\em {Progr.\ Probab.}}, pages 261--277. {Birkh{\"{a}}user/Springer Basel
  AG}, Basel, 2011.

\bibitem{VakTarCho87}
{Vakhania, N.~N.}, {Tarieladze, V.~I.}, and {Chobanyan, S.~A.}
\newblock {\em {Probability Distributions on Banach Spaces}}, volume~14 of {\em
  {Mathematics and its Applications (Soviet Series)}}.
\newblock {D.~Reidel Publishing Co.}, {Dordrecht}, 1987.

\bibitem{VerGelGra75}
{Vershik, A.~M.}, {Gel'fand, I.~M.}, and {Graev, M.~I.}
\newblock {Representations of the Group of Diffeomorphisms}.
\newblock {\em {Russ.~Math.~Surv.+}}, 30(6):1--50, 1975.

\bibitem{Vil09}
{Villani, C.}
\newblock {\em {Optimal transport, old and new}}, volume 338 of {\em
  {Grundlehren der mathematischen Wissenschaften}}.
\newblock {Springer-Verlag}, 2009.

\end{thebibliography}

}
\end{document}